\documentclass[a4paper,12pt]{article}

\usepackage[utf8]{inputenc}
\usepackage[top=2cm,bottom=2cm,left=1.6cm,right=1.6cm]{geometry}
\usepackage{xcolor}
\usepackage{enumerate}
\usepackage{longtable}
\usepackage{rotating}
\usepackage{pdflscape}
\usepackage{multirow}

\usepackage{mathpazo}

\title{New families of vertex-transitive quaternionic lattices in products of trees}
\author{Justin Vast\thanks{Justin Vast is an F.R.S.-FNRS Research Fellow, and is supported in part by the FWO and the
F.R.S.-FNRS under the EOS programme (project ID 40007542).}}
\date{\today\vspace{-3ex}}

\usepackage[english]{babel}

\usepackage{amsmath,amsfonts,amssymb,amsthm}

\usepackage{hyperref} 
\usepackage[english]{cleveref}

\allowdisplaybreaks

\theoremstyle{plain}
\newtheorem{theorem}{Theorem}[section]
\newtheorem{theo}[theorem]{Theorem}
\newtheorem{proposition}[theorem]{Proposition}
\newtheorem{lemma}[theorem]{Lemma}
\newtheorem{corollary}[theorem]{Corollary}

\newtheorem*{theoremstar}{Theorem}
\theoremstyle{definition}
\newtheorem{definition}[theorem]{Definition}

\theoremstyle{remark}
\newtheorem{remark}[theorem]{Remark}
\newtheorem{notation}[theorem]{Notation}

\addtocontents{toc}{\protect\setcounter{tocdepth}{1}}

\newcommand{\ZZ}{\mathbb{Z}}
\newcommand{\QQ}{\mathbb{Q}}

\newcommand{\FF}{\mathbb{F}}
\newcommand{\TT}{\mathbb{T}}
\newcommand{\Aut}{\operatorname{Aut}}
\newcommand{\Typ}{\operatorname{Typ}}
\newcommand{\Ram}{\operatorname{Ram}}
\newcommand{\SL}{\operatorname{SL}}
\newcommand{\PSL}{\operatorname{PSL}}

\newcommand{\PGL}{\operatorname{PGL}}
\newcommand{\End}{\operatorname{End}}
\newcommand{\st}{\hspace{.1cm}|\hspace{.1cm}}
\newcommand{\Nrd}{\operatorname{Nrd}}
\newcommand{\Trd}{\operatorname{Trd}}
\newcommand{\Idl}{\operatorname{Idl}}
\newcommand{\PIdl}{\operatorname{PIdl}}

\newcommand{\Cl}{\operatorname{Cl}}
\newcommand{\ClR}{\operatorname{Cl}_\mathrm{R}}
\newcommand{\mass}{\operatorname{mass}}
\newcommand{\localiz}[2]{#1[{#2}^{-1}]}
\newcommand{\Sym}{\operatorname{Sym}}
\newcommand{\Alt}{\operatorname{Alt}}
\newcommand{\charac}{\operatorname{char}}

\usepackage{tikz}
\usetikzlibrary{calc}
\usetikzlibrary{math}
\usetikzlibrary{external}

\usepackage{graphicx}

\usepackage{natbib}
\usepackage{array}
\newcolumntype{L}[1]{>{\raggedright\arraybackslash}p{#1}}

\renewcommand{\footnotemark}{}

\begin{document}

\maketitle

\begin{abstract}
We describe infinite families of quaternionic arithmetic lattices in $\Aut(\mathcal{T}_1)\times \cdots \times \Aut(\mathcal{T}_n)$ acting transitively on the vertices of the Cartesian product $\mathcal{T}_1\times \cdots \times \mathcal{T}_n$, where $\mathcal{T}_1,\dots, \mathcal{T}_n$ are regular trees.
We present an algorithm to determine a presentation for each of these lattices.
We also introduce a method for constructing such lattices that act regularly on the set of vertices.
We illustrate our results with numerous explicit examples, including several families of known cases, as well as many new instances in arbitrary characteristic.
\end{abstract}

\renewcommand{\baselinestretch}{0.2}\normalsize
\tableofcontents
\renewcommand{\baselinestretch}{1.0}\normalsize

\section{Introduction} \label{Introduction}

\subsection{Review of known results}\label{subsec:ReviewOfKnownResults}

Let $\mathcal{T}$ and $\mathcal{T}'$ be two regular trees of degrees at least two.
A \textbf{BMW group} is a subgroup $\Gamma \leqslant \Aut(\mathcal{T}) \times \Aut(\mathcal{T}')$ acting regularly (i.e.\ acting freely and transitively) on the vertices of $\mathcal{T}\times \mathcal{T}'$, the Cartesian product of the trees.
The three letters BMW stand for the initials of M.\ Burger, S.\ Mozes and D.\ Wise,
who studied those groups extensively.
For an overview on these groups, see \cite{Caprace17}, where the name BMW was suggested.
One may generalise the definition of BMW groups to include any subgroup $\Gamma \leqslant \Aut(\mathcal{T}_1)\times \cdots \times \Aut(\mathcal{T}_n)$ acting regularly on the vertices of $\mathcal{T}_1\times \cdots \times \mathcal{T}_n$, where $\mathcal{T}_1,\dots, \mathcal{T}_n$ are regular trees and $n \geq 1$.

A. Lubotzky constructed in \cite[Section 7.4]{Lubotzky94} an arithmetic BMW group acting on a single regular tree of degree $p+1$, for each prime $p$ congruent to $1$ modulo $4$.
His goal was to give explicit examples of Ramanujan graphs coming from congruence quotients of those groups.
To do so, he considered the $\QQ$-algebraic group $G'$ such that 
\[G'(\QQ) = \PGL_1(B) = B^\times/\QQ^\times,\]
where $B$ is $(-1,-1\st \QQ)$, the Hamilton quaternion algebra over the field of rational numbers.
The ramified places of this quaternion algebra are $\infty$ and $2$.
To properly define the group $\Gamma'_S = G'(\localiz{\ZZ}{S})$ of $S$-integral points of $G'$, the usual $\QQ$-basis $1,i,j,k$ of $B$ is used. 
By \cite[Lemma 7.4.1]{Lubotzky94}, for a prime $p$ congruent to $1$ modulo $4$, the congruence subgroup $\Gamma'_p(2)$ of $\Gamma'_p$, which is by definition the kernel of the reduction map
\[\Gamma'_p \rightarrow G'(\localiz{\ZZ}{p}/(2)),\]
happens to be torsion-free and to have a regular action on the vertices of the Bruhat-Tits tree of $\PGL_2(\QQ_p)$, through the inclusion
\[\Gamma'_p(2) \hookrightarrow G'(\QQ_p) \simeq \PGL_2(\QQ_p).\]
In other words, $\Gamma'_p(2)$ is a torsion-free BMW group acting on a single regular tree of degree $p+1$.
However, if $p$ is congruent to $3$ modulo $4$, the congruence subgroup $\Gamma'_p(2)$ is not vertex-transitive on the Bruhat-Tits tree and is therefore not a BMW group.  

Lubotzky remarked in \cite[Remarks 7.4.4 (iv)]{Lubotzky94} that his construction was made possible by the fact that the Hurwitz maximal order
\[O = \ZZ \oplus \ZZ i \oplus \ZZ j \oplus \ZZ \frac{1+i+j+k}{2}\]
has class number one, that is, every right $O$-ideal is principal.
Indeed, his main ingredient to prove the transitivity is Jacobi's four-square theorem, which is a consequence of $O$ having class number $1$.
This theorem states that if $p$ is an odd prime number, there are exactly $8(p+1)$ quadruples $(a,b,c,d)$ of integers such that $p = a^2 + b^2 + c^2 + d^2$.
One can see such a quadruple as an element of the Lipschitz order
\[\ZZ \oplus \ZZ i \oplus \ZZ j \oplus \ZZ k \subseteq B \]
of reduced norm $p$.
Then, one proves that the subgroup generated by those elements is vertex-transitive on the Bruhat-Tits tree of $\PGL_2(\QQ_p)$.
Note that Lubotzky chose the basis $1,i,j,k$ of the non-maximal Lipschitz order to define the $S$-integral points of his $\QQ$-algebraic group, and not a basis of the Hurwitz maximal order.
Nevertheless, his remark suggested that one could carry out the same arguments for other quaternion algebras over $\QQ$ (or even over global fields), upon the condition that a maximal order has class number one.
This was done by P. Chiu in \cite{Chiu92}, who constructed $3$-regular Ramanujan graphs using techniques in the same vein, but with a maximal order of class number one in $(-2,-13 \st \QQ)$, which is the quaternion algebra over $\QQ$ ramifying at $\infty$ and $13$.
Another occurrence of this strategy was given by M. Morgenstern in \cite{Morgenstern94}, who constructed more Ramanujan graphs by considering a maximal order of class number one in the quaternion algebra over the global function field $\FF_q(t)$ ramifying at the infinite place and the place $t-1$.

In his paper \cite{Mozes95}, S. Mozes considered the arithmetic group
\[\Gamma'_{p,q} = G'(\ZZ[\{p,q\}^{-1}]),\]
inside the same algebraic group $G'$ as Lubotzky \cite{Lubotzky94}, with respect to the same $\QQ$-basis $1,i,j,k$ of $B$, where $p,q$ are two distinct prime numbers both congruent to $1$ modulo $4$.
This time, $\Gamma'_{p,q}$ acts on the Cartesian product of the Bruhat-Tits trees of $\PGL_2(\QQ_p)$ and $\PGL_2(\QQ_q)$ through the inclusion
\[\Gamma'_{p,q} \hookrightarrow G'(\QQ_p) \times G'(\QQ_q) \simeq \PGL_2(\QQ_p) \times \PGL_2(\QQ_q).\]
Mozes proved that the congruence subgroup $\Gamma'_{p,q}(2)$ is a torsion-free BMW group for its action on the product of the corresponding Bruhat-Tits trees.
Moreover, Mozes' groups are irreducible cocompact lattices in $\PGL_2(\QQ_p) \times \PGL_2(\QQ_q)$, and the closures of the projections respectively contain $\PSL_2(\QQ_p)$ and $\PSL_2(\QQ_q)$.

D. Rattaggi went one step further in his thesis \cite{Rattaggi04}, by generalising the ideas of Mozes.
Namely, to any pair of distinct odd primes $p$ and $q$, not necessarily congruent to $1$ modulo $4$, he constructed a torsion-free irreducible lattice in $\PGL_2(\QQ_p) \times \PGL_2(\QQ_q)$ acting regularly on the vertices, still by using the Lipschitz order in the Hamilton quaternion algebra over $\QQ$.
As Lubotzky, his main ingredient to prove the transitivity is Jacobi's four-square theorem.
He then cleverly found the right generators for his groups.
We will see in the first BMW family in \Cref{subsec:-1_-1_Q} that this choice of generators comes from a congruence condition on a larger arithmetic group.

Another step was achieved by the works of 
N. Rungtanapirom, J. Stix and A. Vdovina \cite{RungStixVdov19}, who constructed families of BMW groups, yielding explicit cubical Ramanujan complexes from congruence subgroups.
They also used the Hamilton quaternion algebra $B$ over $\QQ$, but unlike their predecessors, to define the $S$-integral points of the $\QQ$-algebraic group $G$ such that $G(\QQ) = B^\times / \QQ^\times$, they used the basis $1,i,j,(1+i+j+k)/2$ of the maximal Hurwitz order, not the basis of the Lipschitz order.
They found out that one could take any finite set $S$ of odd prime numbers and still obtain a BMW group acting on the product of the Bruhat-Tits trees of $\PGL_2(\QQ_p)$, for $p\in S$.
More precisely, they proved in \cite[Proposition 3.11]{RungStixVdov19} that the congruence subgroup $\Gamma_{S}(2)$ is BMW, where
\[\Gamma_{S} = G(\localiz{\ZZ}{S}),\]
and $\Gamma_{S}(2)$ is the kernel of
\[\Gamma_{S} \rightarrow G(\localiz{\ZZ}{S}/(2)).\]
Note that in their original work, $\Gamma_{S}(2)$ is not seen as a congruence subgroup.
Indeed, they chose to work with $\Gamma_{S\cup \{2\}}$, which does not have a congruence quotient at $2$ since $2$ is invertible in $\ZZ[(S\cup \{2\})^{-1}]$.
To overcome this, they introduced the ``$2$-adic filtration'' of $\Gamma_{S\cup \{2\}}$.
After a close look at the definitions, one observes that this filtration is exactly the congruence subgroup $\Gamma_{S}(2)$.
Their proof of transitivity is again based on Jacobi's four-square theorem.
For the free action of $\Gamma_{S}(2)$, their approach is different, though.
They computed the vertex-covolume $c = 1/24$ of $\Gamma_{S\cup \{2\}}$ using Prasad's volume formula, and proved that $\Gamma_{S}(2)$ has index $1/c = 24$: this subgroup then has a vertex-covolume of $1$.
Together with the transitive action, it must act freely, as wanted.
Moreover, they give the conditions on $S$ for $\Gamma_{S}(2)$ to be torsion-free, and they mention that if $S=\{p,q\}$ with $p,q\equiv 1 \mod 4$, then one obtains the same group as Mozes, i.e.\ $\Gamma_{S}(2) = \Gamma'_{S}(2)$.

By imitating the choice of generators made by Rattaggi and using the same algebraic group $G$ as previously, they also found another infinite family of torsion-free BMW groups in \cite[Theorem 3.15]{RungStixVdov19}, extending Rattaggi's findings.
It is not explicitly stated in their article, but as we will see in the second BMW family in \Cref{subsec:-1_-1_Q}, each of these BMW groups is a congruence subgroup of $\Gamma_S$ for some $S$. 

In \cite{RungStixVdov19} and \cite{StixVdov17}, the same authors applied the same recipe to a maximal $\FF_q[t]$-order $O$, in the quaternion algebra over $\FF_q(t)$ which ramifies at the places $t$ and $\infty$, where $q>1$ is any prime power.
In this way too, they found new BMW groups.
This construction is a generalisation of \cite{Morgenstern94}, passing from an action on one tree to an action on a product of trees.
See \Cref{subsec:ExFqtDeg1} for additional details.

Let us mention Rungtanapirom's thesis \cite[Chapter 4]{Rungtanapirom18}, where one of his main results was the construction of fake quadrics.
To do so, he used a BMW group acting on a product of two $3$-regular trees, from the quaternion algebra over $\FF_2(t)$ ramifying at the places $t+1$ and $t^2+t+1$.
We will see in \Cref{subsec:ExFqtDeg2} that this group is part of an infinite family of BMW groups.
To the best of our knowledge, the description of that family is new and original.

Concluding the state of the art, more recently in \cite{MendelYu26}, J. Mendel and J. Yu constructed families of BMW groups from orders of class number one in totally definite quaternion algebras over number fields.
This generalises all the results cited above in characteristic zero.
Using congruence subgroups, they also give new explicit examples of cubical Ramanujan complexes.
This addresses Lubotzky's remark \cite[Remarks 7.4.4 (iv)]{Lubotzky94} mentioned earlier.

The work of Mendel and Yu shares some similarities with ours, although both works were carried out at the same time, completely independently and without consultation.
Nevertheless, there are some substantial differences: they mainly consider orders of class number one, and only work over number fields, while in our work, dynamical aspects such as the construction of Ramanujan cube complexes are not addressed.

\subsection{Main results and content of the paper}

The present paper aims to provide a unified theory for arithmetic quaternionic groups acting vertex-transitively on a Cartesian product of Bruhat-Tits trees.
We then use this theory to construct finite index subgroups of those vertex-transitive quaternionic lattices that still act transitively, but also freely, and are thus BMW-groups. Note however that in general, a vertex-transitive quaternionic lattice need not contain a vertex-transitive subgroup acting freely, see \Cref{subsec:-1_-1_Qsqrt69}.
As a by-product, we will recover all the families of BMW groups cited above together with many new examples.

We start with \Cref{Preliminaries}, where we lay out the necessary preliminaries.
We first recall what the Bruhat-Tits tree of $\PGL_2(K_v)$ for a local field $K_v$ is.
The most adequate definition for our purposes will be to view the vertices as maximal orders of $M_2(K_v)$.
Afterwards, we expose some fundamental results concerning quaternion algebras, and specifically their arithmetic.

We then study arithmetic quaternionic lattices in \Cref{sec:TransQuatLat}.
Let $B$ be a division quaternion algebra over a global field $K$, and let $O$ be a maximal order in $B$.
We assume $B$ to be totally definite, that is, the subset $\Ram(B)\subseteq V^K $ of the $K$-places ramified in $B$ contains $V^K_\infty$, the set of archimedean places.
In \Cref{subsec:VKO}, we study an arithmetic invariant set  defined in \Cref{def:VKO} as
\[
V^K(O) := \{ v \in V^K\backslash \Ram(B) \st \text{all the $v$-neighbours of $O$ have the same type as $O$}\}.
\]
The main result about this set is the following.

\begin{theoremstar}[\Cref{thm:VKOnonempty}]
If the set $V^K(O)$ is non-empty, then it has an explicit description which only depends on $B$, and it contains the totally positive principal prime ideals not in $\Ram(B)$.
\end{theoremstar}

In \Cref{subsec:GammaOS}, we then study the arithmetic lattice $\Gamma_O(S)$ which has a natural action on a product of Bruhat-Tits trees indexed by $S\backslash \Ram(B)$.
We prove the following.

\begin{theoremstar}[\Cref{lemma:VertexStabGammaS}, \Cref{cor:GammaSTransIfSContainsRam}]
If $\Ram(B)\subseteq S$, then $\Gamma_O(S)$ acts vertex-transitively on the product of trees $\prod_{v\in S\backslash \Ram(B)} \mathcal{T}_v$ if and only if $S\backslash \Ram(B) \subseteq V^K(O)$.
If those conditions are satisfied, then the vertex-stabilisers are conjugate to $\Gamma_O(\Ram(B))$.
\end{theoremstar}

From this result, it is natural to examine maximal orders $O$ for which $V^K(O)$ is non-empty.
Maximal orders of type number one are the typical examples since they are the ones which satisfy $V^K(O) = V^K \backslash \Ram(B)$ by definition.
They were classified in the number field case by M. Kirschmer and D. Lorch in \cite{KirschmerLorch16} (see \Cref{thm:classiftype1} and \Cref{Appendix}).
In positive characteristic, we restate in our setting a theorem in \cite{Gekeler91} of E.-U. Gekeler on Drinfeld modules (see \Cref{thm:poschartypenum}).
It does not give a classification, but instead an explicit formula for the type number of specific quaternion algebras.

Kirschmer and Lorch also classified the maximal orders of type number $2$ in the number field case.
In \Cref{prop:VKOnonemptyType2}, we identify the ones with $V^K(O)$ non-empty.
We also enumerate some other sporadic examples we found of type number $>2$ in \Cref{subsec:TransLattTypeGt1}.
We do not know whether there only exist finitely many examples, but we expect more to exist. 
In any case, a classification seems to be out of reach at the moment.

In \Cref{subsec:FamOfVertTransSubGrps}, we address the question whether $\Gamma_O(S)$ has proper subgroups that remain vertex-transitive. 
We provide a recipe to construct such subgroups that also act freely, and are thus BMW groups. 
The conditions required by our construction need not hold in general: indeed, as mentioned above, we will construct in \Cref{subsec:-1_-1_Qsqrt69} examples where $\Gamma_O(S)$ is vertex-transitive but does not contain any subgroup acting freely and transitively.

\Cref{sec:GrpPres} is dedicated to group presentations.
The main result is \Cref{thm:PresentationGamma}, which is very similar to \cite[Algorithm 2.6]{infauthors15}, but slightly more general to apply in our setting.
It is also a special case of the more general results in \cite{Macbeath64} or \cite{Brown84}, but we take advantage of the structure of a product of trees.

Finally, in Sections \ref{sec:ExNumFields} and \ref{sec:ExFuncFields}, we give concrete examples of transitive quaternionic lattices with a main focus on BMW groups, firstly in the number field case, and secondly in the global function field case.
We revisit all the cited examples using the tools we developed, and provide many more original examples.
We made extensive use of the Magma algebra system \cite{Magma} for these two sections.

The \Cref{Appendix} contains a description of the classification of quaternion algebras of type number one over number fields coming from \cite{KirschmerLorch16}, together with relevant information for seeking vertex-transitive subgroups.

\subsection*{AI statement}

No artificial intelligence was used.

\subsection*{Acknowledgements}

I sincerely thank Pierre-Emmanuel Caprace and François Thilmany for their interest in my work and for their generous guidance throughout the completion of this project.

\section{Preliminaries} \label{Preliminaries}

\subsection{Bruhat-Tits trees}\label{subsec:BruhatTitsTrees}
We recall the construction of the \textbf{Bruhat-Tits tree} $\mathcal{T}_v$ corresponding to a non-Archimedean place $v$ in a global field $K$.
See \cite[Chapitre II §1]{Serre77}, \cite[Chap.~II~§2]{Vigneras80} or \cite[Section 23.5]{Voight21}.

Let $K_v$ be the completion of $K$ at $v$, let $R_v$ be its valuation ring, and let $k_v$ be its residue field.
The vertices of $\mathcal{T}_v$ consist of the maximal orders of $M_2(K_v)$, and two maximal orders $O_v$ and $O'_v$ are neighbours if and only if $O_v /O_v\cap O'_v$ and $O'_v / O_v\cap O'_v$ are both isomorphic to $k_v$ as $R_v$-modules.
The graph constructed in this way is a regular tree of degree $q_v+1$, where $q_v$ is the cardinality of $k_v$.

One may alternatively define $\mathcal{T}_v$ as follows.
An $R_v$-lattice $L$ is a finitely generated $R_v$-submodule of $K_v^2$.
It is said to be full if $LK_v = K_v^2$, and two lattices $L$ and $L'$ are said to be homothetic if there exists $\alpha\in K_v^\times$ such that $L' = \alpha L$.
Let $L,L'$ be full lattices.
Since $R_v$ is a discrete valuation ring and $L$ and $L'$ are torsion-free $R_v$-modules, it follows by the invariant factors theorem that they are both free of rank $2$.
Moreover, one can find a basis $\{u,v\}$ of $L$ and integers $k,l$ such that $\{\pi^k u, \pi^l v\}$ is a basis of $L'$, the element $\pi\in R_v$ being a uniformizer.
The integer $|k-l|$ depends only of the representatives of the homothety classes of $L$ and $L'$, and is called the distance between the classes.
The Bruhat-Tits tree is then the graph whose vertices are the homothety classes of full $R_v$-lattices in $K_v^2$, and two homothety classes are neighbours if and only if they are at distance $1$.

Note that $\End_{R_v}(L) = \End_{R_v}(L')$ if and only if $L$ and $L'$ are homothetic.
Moreover, the order $\End_{R_v}(L)$ is maximal for any $R_v$-lattice $L$ and conversely, given a maximal order $O$ in $M_2(K_v)$, there exists a $R_v$-lattice $L$ such that $O = \End_{R_v}(L)$.
This correspondence between maximal orders in $M_2(K_v)$ and homothety classes of lattices in $K_v^2$ gives rise to a graph isomorphism between the two constructions.
We will use the first definition here.

Importantly, the group $\PGL_2(K_v)$ acts faithfully and vertex-transitively on the Bruhat-Tits tree $\mathcal{T}_v$ by conjugating the maximal orders, and the vertex stabilisers are the conjugates of $\PGL_2(R_v)$.

\subsection{Quaternion algebras and orders}\label{subsec:QuatAlgAndOrders}
We now present the results on quaternion algebras and orders that we need for this article.
The reader is invited to consult the references \cite{Vigneras80} and \cite{Voight21} for more details.

Let $K$ be a global field and write $V^K$ for its set of places. 
If $K$ is a number field, let $R$ be its ring of integers, and let $V_\infty^K$ be the set of Archimedean places.
If $K$ has positive characteristic, fix a place $v_{\infty}\in V^K$, and let $R$ be the ring $\bigcap_{v\in V^K\backslash\{v_{\infty}\}} (R_v \cap K)$, where $R_v$ is the valuation ring of the completion $K_v$; one then sets $V_\infty^K$ to be $\{v_\infty\}$.
In any case, let $V_f^K$ be $V^K\backslash V_\infty^K$.
We will usually identify $V_f^K$ with the set of prime ideals of $R$, or even the set of prime elements in $R$ (up to a unit) if $R$ is a principal ideal domain.
If $K = \FF_q(t)$, for example, one can take $v_\infty = - \deg$; in that case, $R = \FF_q[t]$ and $V_f^K$ can be seen as the set of monic irreducible polynomials over $R$.

Let $B$ be a quaternion algebra over $K$.
If $K$ has characteristic not $2$, we will write 
\[(a,b\st K)\]
for the quaternion algebra with basis $1,i,j,k$ such that $i^2 = a$, $j^2 = b$ and $i j = k = -j i$.
If $K$ has characteristic $2$, then 
\[[a,b\st K)\]
is the quaternion algebra with basis $1,i,j,k$ such that $i^2+i = a$, $j^2 = b$ and $i j = k = j (i+1)$.

\subsection{From local to global}\label{subsec:FromLocalToGlobal}
Let us fix an embedding $\iota_v : K\hookrightarrow K_v$ for every $v\in V^K$, restricting to an embedding $\theta_v : R\hookrightarrow R_v$ for the non-Archimedean places.
The quaternion algebra $B$ is said to be \textbf{ramified at a place} $v\in V^K$ if the completion $B_v = B \otimes_{\iota_v} K_v$ is a division algebra.
It is \textbf{unramified} otherwise, that is to say, if $B_v$ is isomorphic to $M_2(K_v)$ as a $K_v$-algebra.
We denote by $\Ram(B)$ the set of ramified places for $B$.
It is a finite set of even cardinality containing no complex place.
Conversely, if $S \subseteq V^K$ is a finite subset of even cardinality containing no complex place, there exists a quaternion algebra $B$ with $\Ram(B) = S$.
Moreover, two quaternion algebras over $K$ are isomorphic if and only if they ramify at the same places.
In particular, $B$ is isomorphic to $M_2(K)$ if and only if $B_v$ is isomorphic to $M_2(K_v)$ for all $v\in V^K$: this is a particular case of the Albert–Brauer–Hasse–Noether theorem.
See \cite[Chap.~III Théorème~3.1]{Vigneras80} or \cite[Main Theorem~14.6.1, Corollary~14.6.5]{Voight21}).

There is a \textbf{local-global correspondence} between orders in $B$ and orders in the completions $B_v = B \otimes_{\iota_v} K_v$, for $v \in V_f^K$, mapping an order $O$ in $B$ to the family $(O_v)_{v\in V_f^K}$, where $O_v$ is the order $ O \otimes_{\theta_v} R_v$ in $B_v$.
The family of completions fully determine the order $O$.
If $O'$ is another order in $B$, then the families $(O_v)_{v\in V_f^K}$ and $(O'_v)_{v\in V_f^K}$ coincide on almost every place.
Conversely, given a family of local orders $(O'_v)_{v\in V_f^K}$ coinciding at almost every place with $(O_v)_{v\in V_f^K}$, there exists a (unique) order $O'$ in $B$ having those completions (cf.\ \cite[Chap.~III Proposition~5.1]{Vigneras80}, \cite[Theorem~9.4.9]{Voight21}).
Maximality is a local property, namely $O$ is a maximal order if and only if $O_v$ is a maximal order for every place $v\in V_f^K$ (cf.\ \cite[Lemma 10.4.3]{Voight21}).
Also, every order is contained in a maximal order (cf.\ \cite[Chap.~I Proposition~4.2]{Vigneras80}, \cite[10.4.2]{Voight21}).
It follows that, given any order $O$ in $B$, the completion $O_v$ is maximal for almost every $v\in V_f^K$.

\subsection{The geometry of the set of maximal orders}\label{subsec:GeomOfMaxOrders}
For a place $v\in V^K_f$ not ramified in $B$, the algebra $B_v$ is isomorphic to $M_2(K_v)$ by definition, and the maximal orders in $B_v$ are in correspondence with the vertices in the Bruhat-Tits tree at $v$ by what we explained earlier in \Cref{subsec:BruhatTitsTrees}.
On the other hand, if $v$ ramifies $B$, then $B_v$ is a division algebra and the valuation uniquely extends to $B_v$ by defining $v(x) := v(\Nrd(x))/2 \in \frac{1}{2}\ZZ$ for all $x\in B_v$, where $\Nrd$ is the reduced norm (cf.\ \cite[Chap.~II Lemma 1.4]{Vigneras80}, \cite[Lemma 13.2.1]{Voight21}).
In this case, there is a unique maximal order in $B_v$, namely $O_v = \{x\in B_v \st v(x) \ge 0\}$ (cf.\ \cite[Chap.~II Lemma 1.5]{Vigneras80}, \cite[Proposition 13.3.4]{Voight21}).
By the local-global correspondence, for a fixed maximal order $O$, there is therefore a bijection between the maximal orders in $B$ and the vertices in the restricted product
\[\sideset{}{'}\prod_{v\in V^K_f \backslash \Ram(B)} \mathcal{T}_v, \]
where a family of local orders $(O'_v)$ in the full product belongs to the restricted product if and only if $O_v = O'_v$ for almost all $v\in V^K\backslash \Ram(B)$.
\begin{notation}\label{not:ProdOfTrees}
If $T$ is a finite set of finite places not ramified in $B$, we shall abusively write
\[\prod_{v\in T} \mathcal{T}_v\]
for the set of maximal \textit{global} orders $O'$ such that $O_w = O'_w$ for all $w\in V^K_f\backslash T$.
Note that the fixed maximal order $O$ is implicit in this notation, and a different choice of maximal order may lead to a different set of maximal orders.
If $T$ is a singleton $\{v\}$, we will simply write $\mathcal{T}_v$.
Beware in this case that a maximal order $O'$ in $\mathcal{T}_v$ is an order in $B$, not in the completion $B_v$.
\end{notation}

For a place $v\in V^K_f$ not ramified in $B$, a \textbf{$v$-neighbour} of a maximal order $O$ is a maximal order $O'$ such that $O'_v$ is a neighbour of $O_v$ in the Bruhat-Tits tree $\mathcal{T}_v$ at the place $v$, and $O'_w = O_w$ for every $w\in V_f^K$ distinct from $v$.
In particular, $O'\in \mathcal{T}_v$ according to our notation.	
The \textbf{$v$-neighbourhood} of $O$ is the set of all the $v$-neighbours of $O$.
The cardinality of the $v$-neighbourhood is equal to the degree of the regular tree $\mathcal{T}_v$, that is to say $|k_v|+1$, where $k_v$ is the residue field of $K_v$.
If $\mathfrak{p}$ denotes the prime $R$-ideal corresponding to the place $v$, there is a bijective correspondence between the set of integral right $O$-ideals of reduced norm $\mathfrak{p}$ and the $v$-neighbourhood of $O$.
The bijection is given by
\[I \mapsto O_{\mathrm{L}}(I),\]
where $O_{\mathrm{L}}(I) = \{ \alpha \in B \st \alpha I \subseteq I \}$ is the left order of $I$ (cf.\ \cite[Lemma 26.4.1]{Voight21}).

\subsection{Type and class numbers}\label{subsec:TypeAndClassNumbers}
Let $O, O'$ be orders in a quaternion algebra $B$.
We say that $O$ and $O'$ belong to the same \textbf{genus} if they are locally isomorphic, that is to say if $O_v$ and $O'_v$ are isomorphic as $R_v$-algebras for all places $v \in V^K_f$.
They are of the same \textbf{type} if they are isomorphic as $R$-algebras.
Having the same type implies having the same genus, but not conversely in general.
We write $\Typ(O)$ for the set of $R$-algebra isomorphism classes of orders having the same genus as $O$.
The \textbf{type number} of $O$ is the cardinal of $\Typ(O)$ and is finite (cf.\ \cite[Chap.~III Corollaire~5.5]{Vigneras80}, \cite[Corollary~27.6.25]{Voight21}).
The type number of a maximal order $O$ counts the number of $R$-algebra isomorphism classes among maximal orders, since an order $O'$ is maximal if and only if it is locally isomorphic to $O$. 
Therefore, we shall refer to the type number of the quaternion algebra $B$ as the type number of any of its maximal orders.

The group $\PGL_1(B) = B^\times/K^\times$ naturally acts on the set of orders in $B$ by conjugation.
If $\gamma\in B^\times$ is a representative of an element $g\in \PGL_1(B)$, and if $O$ is an order, we shall write $g\cdot O$ for the order $\gamma O \gamma^{-1}$.
Note that an isomorphism of orders extends to an automorphism of $B$; by the Skolem-Noether theorem, it follows that two orders have the same type if and only if they are conjugate by an invertible element of the quaternion algebra (cf.\cite[Lemma~17.4.2]{Voight21}).
It follows that the type number of an order $O$ corresponds to the number of orbits for the action of $\PGL_1(B)$ on the set of orders in the same genus as $O$.
In particular, the type number of $B$ is the number of orbits for the action of $\PGL_1(B)$ on the vertices of the restricted product of trees $\sideset{}{'}\prod_{v\in V^K_f \backslash \Ram(B)} \mathcal{T}_v$.

Let $O$ be an order in $B$.
By definition, two right fractional $O$-ideals $I$ and $J$ are in the same \textbf{right class} if they are isomorphic as right $O$-modules.
Equivalently, they are in the same right class if and only if there exists $\alpha\in B^\times$ such that $I = \alpha J$ (cf.\ \cite[Lemma~17.3.3]{Voight21}).
We write $\Cl_\mathrm{R}(O)$ for the right class set of $O$, that is, the set of right classes of invertible right fractional $O$-ideals.
The right class set is finite and its cardinality is called the (right) \textbf{class number} of $O$ (cf.\ \cite[Main Theorem~17.7.1]{Voight21} and \cite[Chap. III Théorème 5.4]{Vigneras80}).
If $O$ is maximal, which is our case of interest in this paper, every right fractional $O$-ideal is invertible (cf.\ \cite[Proposition~16.1.2]{Voight21}).

\subsection{Two-sided ideals}\label{subsec:TwoSidedIdeals}
Let $O$ be a maximal order in $B$.
The set of \textbf{fractional two-sided ideals} of $O$ forms a free abelian group on the set of non-zero prime two-sided ideals (cf.\ \cite[Corollary~18.3.5]{Voight21}).
Let us write that group $\Idl(O)$.
A fractional two-sided ideal $I$ is said to be \textbf{principal} if there exists an element $x\in B^\times$ for which $I = x O = O x O = O x$.
The principal fractional two-sided ideals of $O$ form a subgroup of $\Idl(O)$, denoted by $\PIdl(O)$.
Similarly, let $\Idl(R)$ and $\PIdl(R)$ be the groups of fractional ideals and principal fractional ideals of $R$, respectively.
We write $\Cl(K)$ for the ideal class group of $K$, that is, the quotient $\Idl(R)/\PIdl(R)$.
From here on, we will always write ``ideal'' to mean ``\textit{fractional} ideal''.

There is a bijection given by
\begin{align*}
\{\text{Prime two-sided $O$-ideals}\} & \rightarrow  \{\text{Prime $R$-ideals}\} \\
P & \mapsto  P\cap R
\end{align*}
(cf.\ \cite[Theorem~18.3.6, (18.3.7)]{Voight21}).
The group $\Idl(R)$ of $R$-ideals embeds in $\Idl(O)$ by identifying an $R$-ideal $I$ with the two-sided $O$-ideal $I O$.
If $I = \mathfrak{p}$ is a prime $R$-ideal not dividing the $R$-\textbf{discriminant} $\mathfrak{D} = \prod_{v\in \Ram(B) \cap V^K_f} \mathfrak{p}_v$ of $B$, then $\mathfrak{p} O$ is the (unique) prime two-sided $O$-ideal $P$ such that $P \cap R = \mathfrak{p}$.
On the other hand, if the prime ideal $\mathfrak{p}$ divides $\mathfrak{D}$, there exists a (unique) prime two-sided $O$-ideal $\mathfrak{P}$ such that $\mathfrak{P}\cap R = \mathfrak{p}$ and $\mathfrak{P}^2 = \mathfrak{p} O$, and $\mathfrak{P} \notin \Idl(R)$. 
This paragraph can be summarised by the exact sequence
\[1\rightarrow \Idl(R) \rightarrow \Idl(O) \rightarrow \prod_{\mathfrak{p} | \mathfrak{D}} \ZZ/2\ZZ \rightarrow 1,\]
where the product is over the prime $R$-ideals dividing $\mathfrak{D}$ (cf.\ \cite[Theorem~18.3.6, (18.3.8)]{Voight21}).

The map 
\[\Idl(O)/\PIdl(O) \rightarrow \{ [I]\in \ClR(O) \st O_\mathrm{L}(I) \simeq O \} : I \mapsto [I]\]
is an isomorphism of sets, where $O_{\mathrm{L}}(I)$ is the left order of the ideal $I$ (cf.\ \cite[Proposition~18.5.10]{Voight21}).
Moreover, the map 
\[\ClR(O) \rightarrow \Typ(O) : [I] \mapsto \text{type of } O_{\mathrm{L}}(I)\]
is surjective (cf.\ \cite[Lemma~17.4.13]{Voight21}).
In particular, the class number of $O$ is greater than or equal to the type number of $O$.

\section{Transitive quaternionic lattices}\label{sec:TransQuatLat}
\addtocontents{toc}{\protect\setcounter{tocdepth}{2}}

As in the preliminaries, let $K$ be a global field and write $V^K$ for its set of places. 
If $K$ is a number field, let $R$ be its ring of integers, and let $V_\infty^K$ be the set of Archimedean places.
Otherwise, if $K$ has positive characteristic, fix a place $v_{\infty}\in V^K$, and let $R$ be the ring $\bigcap_{v\in V^K\backslash\{v_{\infty}\}} (R_v \cap K)$, where $R_v$ is the valuation ring of the completion $K_v$; one then sets $V_\infty^K$ to be $\{v_\infty\}$.
In both cases, let $V_f^K$ be $V^K\backslash V_\infty^K$.

Let $B$ be a division quaternion algebra over $K$, and assume that $V_\infty^K \subseteq \Ram(B)$.
If $K$ is a number field, we then say that $B$ is \textbf{totally definite}.
It follows in this case that $K$ is totally real.
Let $O$ be a maximal ($R$-)order in $B$.
Let $S$ be a finite set of places containing $V^K_\infty$, and define $S_0$ to be the intersection $S\cap \Ram(B)$.

We write $\localiz{R}{S}$ for the ring $\bigcap_{v\in V^K \backslash S} (R_v \cap K)$, that is the localization of $R$ inverting the (finite) places in $S$.
Similarly, the notation $\localiz{O}{S}$ stands for the (maximal) $\localiz{R}{S}$-order $\bigcap_{v\in V^K \backslash S} (O_v \cap B) \simeq O \otimes_R \localiz{R}{S}$ in $B$.
We will also write $\PGL_1(B)$ for the group $B^\times/K^\times$.

\subsection{The set \texorpdfstring{$V^K(O)$}{VK(O)}}\label{subsec:VKO}

\begin{definition}\label{def:VKO}
Let us define the set
\[
V^K(O) := \{ v \in V^K\backslash \Ram(B) \st \text{all the $v$-neighbours of $O$ have the same type as $O$}\}.
\]
\end{definition}

\begin{lemma}\label{lemma:inVKOiffTvSameTypeAsO}
A place $v\in V^K\backslash S_0$ belongs to $V^K(O)$ if and only if every order $O'\in \mathcal{T}_v$ has the same type as $O$.
\end{lemma}

\begin{proof}
A $v$-neighbour $O'$ of $O$ satisfies $O_w = O'_w$ for all place $w\in V^K\backslash (S_0 \cup \{v\})$, so the implication from right to left is obvious.

For the other implication, let $v\in V^K(O)$, and let $O_1,\dots,O_n$ be the $v$-neighbours of $O$.
By the Skolem-Noether theorem, there exist elements $g_1,\dots,g_n\in \PGL_1(B)$ such that $g_i\cdot O = O_i$ for all $1\leq i \leq n$.
It is then clear that the subgroup $G$ generated by those elements acts vertex-transitively on the tree $\mathcal{T}_v$, which proves that any $O'\in \mathcal{T}_v$ has the same type as $O$.
\end{proof}

\begin{remark}\label{rem:VKOFullIFFTypeNumberOne}
Recall \Cref{subsec:GeomOfMaxOrders}, where we saw that the set of maximal orders in $B$ is the vertex-set of the restricted product $\sideset{}{'}\prod_{v\in V^K \backslash \Ram(B)} \mathcal{T}_v$.
By \Cref{lemma:inVKOiffTvSameTypeAsO}, a maximal order $O$ then satisfies $V^K(O) = V^K\backslash \Ram(B)$ if and only if every maximal order is isomorphic to $O$, if and only if $O$ has type number $1$ (cf.\ \Cref{subsec:TypeAndClassNumbers}).
\end{remark}

\begin{lemma} \label{lemma:ordersisosameVKO}
Let $O$ and $O'$ be maximal orders of $B$ of the same type.
Then $V^K(O) = V^K(O')$.
\end{lemma}
\begin{proof}
The orders $O$ and $O'$ are conjugate by the Skolem-Noether theorem, so let $\alpha \in B^\times$ be such that $O' = \alpha O \alpha^{-1}$.
Let $v\in V^K(O)$, and let $N_v$ be the $v$-neighbourhood of $O$.
Then $\alpha N_v \alpha^{-1} = \{ \alpha U \alpha^{-1} \st U \in N_v \}$ is the $v$-neighbourhood of $O'$.
The $v$-neighbours are then all isomorphic to a neighbour of $O$, so isomorphic to $O$ since $v\in V^K(O)$, and isomorphic to $O'$ since $O$ and $O'$ are isomorphic by assumption.
Therefore $v\in V^K(O')$, meaning that $V^K(O)\subseteq V^K(O')$. 
By symmetry, $V^K(O) = V^K(O')$.
\end{proof}

Let $\Cl^+(R)$ be the narrow class group of $R$.
This is the quotient $\Idl(R)/\PIdl^+(R)$, where $\PIdl^+(R)$ is the subgroup of non-zero principal fractional ideals generated by a totally positive element.
Recall that an element $a\in K^\times$ is totally positive if $v(a) > 0$ for all real embeddings $v\in V^K_\infty$, and the group of totally positive elements is written $K^\times_{>0}$.
In particular, $\Cl^+(R)$ is the usual class group $\Cl(R)$ if $K$ is a function field.

Let $\mathcal{I}^+_B(R)$ be the subgroup of $\Idl(R)$ generated by $\PIdl^+(R)$, the prime ideals in $\Ram(B)\cap V^K_f$, and the squares of the prime ideals in $V^K\backslash \Ram(B)$.
Obviously, $\Idl(R)/\mathcal{I}^+_B(R)$ is an elementary abelian $2$-group and a quotient of the narrow class group $\Cl^+(R)$.

Alternatively, one can give an idèlic description of this generalised class group $\Idl(R)/\mathcal{I}^+_B(R)$ as follows.
One first view the narrow class group $\Cl^+(R)$ as 
\[\widehat{K}^\times /(K^\times_{>0} \widehat{R}^\times),\]
where $\widehat{K}^\times = \sideset{}{'}\prod_{v\in V^K_f} K_v^\times$ is the finite idèle group of $K$, and $\widehat{R}^\times$ is the subgroup $ \sideset{}{}\prod_{v\in V^K_f} R_v^\times$.
The group $\Idl(R)/\mathcal{I}^+_B(R)$ is then the further quotient
\[\widehat{K}^\times /(K^\times_{>0} \widehat{R}^\times (\widehat{K}^\times)^2 \prod_{v\in \Ram(B)\cap V^K_f} K_v^\times ).\]
Write $\widehat{B}^\times$ for the finite idèle group of $B$, defined as the restricted product $\sideset{}{'}\prod_{v\in V^K_f} B_v^\times$ with respect to the subgroups $O_v^\times\leqslant B_v^\times$:
\[\widehat{B}^\times = \{(x_v)_v \in \sideset{}{}\prod_{v\in V^K_f} B_v^\times \st x_v \in O_v^\times \text{ for all but finitely many $v\in V^K_f$}\}.\]
Now one observes that an element in $\widehat{B}^\times$ normalising $\widehat{O} = \sideset{}{}\prod_{v\in V^K_f} O_v$ has reduced norm in $\widehat{R}^\times (\widehat{K}^\times)^2 \prod_{v\in \Ram(B)\cap V^K_f} K_v^\times \leqslant \widehat{K}^\times$, and reciprocally, for a given element $a$ in this subgroup, there exists a element in $\widehat{B}^\times$ normalising $\widehat{O}$ of reduced norm $a$.
In other words, 
\[\widehat{R}^\times(\widehat{K}^\times)^2 \prod_{v\in \Ram(B)\cap V^K_f} K_v^\times = \Nrd(N_{\widehat{B}^\times}(\widehat{O})),\]
and it follows that
\[\Idl(R)/\mathcal{I}^+_B(R) = \widehat{K}^\times /(K^\times_{>0} \Nrd(N_{\widehat{B}^\times}(\widehat{O}))).\] 
This is exactly the generalised class group $\Cl_{GN(O)} R$ described by Voight in \cite[28.5.9]{Voight21}.

\begin{theo}\label{thm:VKOnonempty}
Assume that $V^K(O)$ is non-empty.
Then
\[V^K(O) = (V^K\backslash \Ram(B)) \cap \mathcal{I}^+_B(R).\]
\end{theo}

\begin{proof}
Let $\mathfrak{p}\in V^K(O)$ and let $I$ be an integral right $O$-ideal of reduced norm $\mathfrak{p}$.
The left order of $I$ is a $\mathfrak{p}$-neighbour of $O$, so they have the same type by hypothesis on $\mathfrak{p}$.
One can then find an element $g\in B^\times$ such that $g\cdot O$ is the left order of $I$.
The right $O$-ideal $g^{-1} I$ is consequently a two-sided $O$-ideal and has a decomposition as a product of prime two-sided $O$-ideals.
From the classification of prime two-sided $O$-ideals, its reduced norm is then a product of prime $R$-ideals in $\Ram(B)\cap V^K_f$ and squares of prime $R$-ideals in $V^K\backslash \Ram(B)$.
In particular, $\Nrd(g^{-1} I) \in \mathcal{I}^+_B(R)$.
On the other hand, the reduced norm of $g$ is totally positive, so that $\Nrd(g)R \in \mathcal{I}^+_B(R)$.
One deduces that 
\[\mathfrak{p} = \Nrd(I) = \Nrd(g) \Nrd(g^{-1} I) \in \mathcal{I}^+_B(R),\]
and this proves the inclusion from left to right.

Now let $\mathfrak{q}\in (V^K\backslash \Ram(B)) \cap \mathcal{I}^+_B(R)$.
We want to show that $\mathfrak{q}$ belongs to $V^K(O)$. 
As $\mathfrak{p}$ belongs to $V^K(O)$ by hypothesis, let us assume that $\mathfrak{q}$ is different from $\mathfrak{p}$.
Let $O'$ be a $\mathfrak{q}$-neighbour of $O$.
There exists an integral right $O$-ideal of reduced norm $\mathfrak{q}$ whose left order is $O'$.
One can multiply it on the right by a well-chosen two-sided $ O$-ideal to get a right $O$-ideal $J$ with a reduced norm in $\PIdl^+(R)$. 
Note that the left order of $J$ is also $O'$.
By \cite[Chap.~III Théorème~5.7]{Vigneras80} (or \cite[Corollary 28.5.5]{Voight21}), for $S = V^K_\infty \cup \{\mathfrak{p}\}$, the right $\localiz{O}{S}$-ideal $J \localiz{O}{S}$ is principal, generated by some element $x\in B^\times$.
Locally, since $\localiz{O'}{S}$ is the left $\localiz{R}{S}$-order of $J \localiz{O}{S} = x \localiz{O}{S}$, one has $x\cdot O_{\mathfrak{r}} = x\cdot (\localiz{O}{S})_{\mathfrak{r}} = (\localiz{O'}{S})_{\mathfrak{r}} =  O'_{\mathfrak{r}}$ for all $\mathfrak{r}\in V^K_f$ different from $\mathfrak{p}$.
Consider $U$, the maximal order such that $U_{\mathfrak{r}} = O_{\mathfrak{r}}$ for all $\mathfrak{r}\in V^K_f\backslash \{\mathfrak{p}\}$ and $U_{\mathfrak{p}} = x^{-1} \cdot O_{\mathfrak{p}}$.
Clearly, $U\in \mathcal{T}_{\mathfrak{p}}$, so by \Cref{lemma:inVKOiffTvSameTypeAsO}, there exists $y\in B^\times$ such that $U = y\cdot O$.
Finally, by looking locally again, one observes that $(xy)\cdot O_{\mathfrak{p}} = x\cdot U_{\mathfrak{p}} = O_{\mathfrak{p}} = O'_{\mathfrak{p}}$, and $(xy)\cdot O_{\mathfrak{r}} = x\cdot U_{\mathfrak{r}} = x\cdot O_{\mathfrak{r}} = O'_{\mathfrak{r}}$ for all $\mathfrak{r}\in V^K_f\backslash \{\mathfrak{p}\}$, that is to say $(xy)\cdot O = O'$.
The orders $O$ and $O'$ have therefore the same type.
As this is true for all $\mathfrak{q}$-neighbour $O'$ of $O$, one has $\mathfrak{q}\in V^K(O)$, and the inclusion from right to left is proven.
\end{proof}

\begin{corollary}\label{cor:NarrowClGrpTrivType1}
If the narrow class group $\Cl^+(\localiz{R}{\Ram(B)})$ is trivial, then $V^K(O)$ is non-empty if and only if $O$ has type number $1$.
\end{corollary}

\begin{proof}
First, if $O$ has type number $1$, then certainly $V^K(O)$ is non-empty by definition, without any assumption on the narrow class group.

Conversely, $\Cl^+(\localiz{R}{\Ram(B)})$ being trivial can be reformulated as $\Idl(R)$ being equal to its subgroup generated by $\PIdl^+(R)$ and the prime ideals in $\Ram(B)\cap V^K_f$.
In particular, this implies that $\mathcal{I}^+_B(R)$ is also equal to $\Idl(R)$.
If $V^K(O)$ is non-empty, then one obtains $V^K(O) = V^K\backslash \Ram(B)$ by \Cref{thm:VKOnonempty}, i.e.\ $O$ has type number $1$.
\end{proof}

\subsection{The lattice \texorpdfstring{$\Gamma_O(S)$}{GammaO(S)}}\label{subsec:GammaOS}

For each place $v\in V^K$, the completion map $\iota_v : K \hookrightarrow K_v$ induces a group embedding $\PGL_1(B) := B^\times/K^\times \hookrightarrow \PGL_1(B_v) := B_v^\times/K_v^\times$, where $B_v := B \otimes_{\iota_v} K_v$.
The group $\PGL_1(B)$ then embeds diagonally in the \textbf{adèle group}, defined as the restricted product
\[\PGL_{1,O}(\mathcal{A}) := \sideset{}{'}\prod_{v\in V^K} \PGL_1(B_v),\]
with respect to the subgroups $\PGL_1(O_v) := O_v^\times / R_v^\times$ for $v\in V^K_f$, and $\PGL_1(B_v)$ for $v\in V^K_\infty$.
In other words, an element
\[(g_v)_v \in \prod_{v\in V^K} \PGL_1(B_v)\]
in the (full) product belongs to the adèle group if and only if $g_v \in \PGL_1(O_v)$ for almost all $v\in V^K_f$.
The \textbf{group of $S$-integral adèles} is the subgroup
\[\PGL_{1,O}(\mathcal{A}(S)) := \prod_{v\in S} \PGL_1(B_v) \times \prod_{v\in V^K\backslash S} \PGL_1(O_v),\]
and from there, we define the \textbf{$S$-arithmetic group}
\[\Gamma_{O}(S) := \PGL_1(B) \cap \PGL_{1,O}(\mathcal{A}(S)). \]
Written differently,
\[\Gamma_{O}(S)  = \{g\in \PGL_1(B) \st \forall v\in V^K\backslash S,  g\in \PGL_1(O_v)\}.\]

\begin{remark}\label{rem:SintegralAdelesActOnProdOfTrees}
The group of $S$-integral adèles acts on the product of trees $\prod_{S\backslash S_0} \mathcal{T}_v$ (we recall \Cref{not:ProdOfTrees}: the vertices of the product are maximal global orders, not families of maximal local orders).
It naturally restricts to an action of $\Gamma_{O}(S)$ on the product of trees.
\end{remark}

\begin{lemma}\label{lemma:GammaS0finite}
The group $\Gamma_O(S_0)$ is finite.
\end{lemma}

\begin{proof}
Since $K$ is discrete in its adèle ring $\mathcal{A}$, so is $\PGL_1(B)$ in $\PGL_{1,O}(\mathcal{A})$, and therefore $\Gamma_O(S_0)$ is discrete in $\PGL_{1,O}(\mathcal{A}(S_0))$.
Since $S_0\subseteq \Ram(B)$, the group $\PGL_1(B_v)$ is compact for every $v$ in $S_0$; the group $\PGL_{1,O}(\mathcal{A}(S_0))$ is thus compact by Tychonoff's theorem.
The group $\Gamma_O(S_0)$, being discrete in a compact group, is consequently finite.
\end{proof}

\begin{lemma}\label{lemma:GammaSCocompactLattice}
If $S_0$ is strictly contained in $S$, then the group $\Gamma_O(S)$ is a cocompact lattice in the product $\prod_{v\in S\backslash S_0} \PGL_2(K_v)$.
\end{lemma}

\begin{proof}
Since $B$ is a division algebra, $\PGL_1(B)$ is anisotropic as an algebraic $K$-group.
By the generalisation of the Borel--Harish-Chandra theorem to any characteristic (see for instance \cite[Theorem 3.2.4]{Margulis91}), the group $\Gamma_O(S)$ is a cocompact lattice in $\prod_{v\in S} \PGL_1(B_v)$.
As $S_0 = S \cap \Ram(B)$ by definition, the group $\prod_{v\in S_0} \PGL_1(B_v)$ is compact, and $\Gamma_O(S)$ is then a cocompact lattice in $\prod_{v\in S\backslash S_0} \PGL_1(B_v)$ since $S_0$ is strictly contained in $S$ by hypothesis.
Moreover, for all $v\in S\backslash S_0$, one has $v\notin \Ram(B)$, so that $B_v \simeq M_2(K_v)$ and $\prod_{v\in S\backslash S_0} \PGL_1(B_v) \simeq \prod_{v\in S\backslash S_0} \PGL_2(K_v)$.
\end{proof}

\begin{lemma}\label{lemma:VertexStabGammaS}
There are finitely many orbits for the action of $\Gamma_{O}(S)$ on the vertices of $\prod_{v\in S\backslash S_0} \mathcal{T}_v$, and the stabiliser of the vertex $O$ is $\Gamma_O(S_0)$.
\end{lemma}

\begin{proof}
The lemma is trivial if $S_0 = S$, so assume $S_0$ to be strictly contained in $S$.

As $\Gamma_{O}(S)$ is a cocompact lattice in $\prod_{v\in S\backslash S_0} \PGL_2(K_v)$ by \Cref{lemma:GammaSCocompactLattice}, and as the latter group acts cocompactly on $\prod_{v\in S\backslash S_0} \mathcal{T}_v$, the group $\Gamma_{O}(S)$ also acts cocompactly on the product of trees.
In particular, there are finitely many orbits of vertices.

By \Cref{prop:ExactSeqGammaScupRam}, the stabiliser of $O$ for the action of $\PGL_1(B)$ on the product of trees is $\Gamma_O(\Ram(B))$.
The stabiliser of $O$ for the action of $\Gamma_O(S)$ is then the intersection of $\Gamma_O(S)$ and $\Gamma_O(\Ram(B))$, that is $\Gamma_O(S\cap \Ram(B)) = \Gamma_O(S_0)$.
\end{proof}

\begin{proposition}\label{prop:ExactSeqGammaScupRam}
One has
\[\Gamma_O(S\cup \Ram(B))  = \{g\in \PGL_1(B) \st g\cdot O \in \prod_{v\in S\backslash S_0} \mathcal{T}_v\},\]
and one has an exact sequence of the following form:
\[1 \rightarrow \Gamma_O(S) \rightarrow \Gamma_O(S\cup \Ram(B)) \rightarrow \prod_{v\in \Ram(B) \backslash S_0} \ZZ/2\ZZ.\]
In particular, the group $\Gamma_O(S)$ is a normal subgroup of $\Gamma_O(S\cup \Ram(B))$ of index a power of $2$.
Moreover, the map on the right-hand side is surjective if $\Idl(\localiz{O}{S_0}) = \PIdl(\localiz{O}{S_0})$.
\end{proposition}

\begin{proof}
We first prove the equality $\Gamma_O(S\cup \Ram(B))  = \{g\in \PGL_1(B) \st g\cdot O \in \prod_{v\in S\backslash S_0} \mathcal{T}_v\}$.
The inclusion from left to right is clear, since $\Gamma_O(S\cup \Ram(B))$ acts on the product of trees by \Cref{rem:SintegralAdelesActOnProdOfTrees}.
Reciprocally, let $g\in \PGL_1(B)$ such that $g \cdot O \in \prod_{v\in S\backslash S_0} \mathcal{T}_v$.
For all $v\in V^K\backslash (S\cup \Ram(B))$, one has $g\cdot O_v = O_v$.
Since $v$ does not ramify the quaternion algebra, the stabiliser of $O_v$ for the action of $\PGL_1(B_v)\simeq \PGL_2(K_v)$ on the Bruhat-Tits tree $\mathcal{T}_v$ is $\PGL_1(O_v)$, so $g\in \PGL_1(O_v)$ and we get the other inclusion.

We now describe the exact sequence.
For $v\in \Ram(B)\cap V^K_f$, the valuation extends uniquely to a surjective morphism $v : B_v^\times \rightarrow \frac{1}{2}\ZZ$.
It then quotients to a surjective morphism $\overline{v} : \PGL_1(B_v) \rightarrow (\frac{1}{2}\ZZ) / \ZZ \simeq \ZZ/2\ZZ$ whose kernel is exactly $\PGL_1(O_v)$.
The group $\PGL_1(O_v)$ is therefore normal of index $2$ in $\PGL_1(B_v)$, and one has the short exact sequence
\[1 \rightarrow \PGL_{1,O}(\mathcal{A}(S)) \rightarrow \PGL_{1,O}(\mathcal{A}(S\cup \Ram(B))) \rightarrow \prod_{v\in \Ram(B)\backslash S_0} \PGL_1(B_v)/\PGL_1(O_v) \rightarrow 1.\]
By intersecting the $S$- and $(S\cup \Ram(B))$-integral adèles by $\PGL_1(B)$, one gets the desired exact sequence.

Assume now that $\Idl(\localiz{O}{S_0}) = \PIdl(\localiz{O}{S_0})$.
Let $\mathfrak{p}\in \Ram(B)\backslash S_0$ a prime ideal corresponding to a valuation $v_\mathfrak{p}$, and let $P_\mathfrak{p}$ be the (unique) prime two-sided $O$-ideal such that $P_\mathfrak{p}\cap R = \mathfrak{p}$.
By the assumption, the $\localiz{O}{S_0}$-ideal $P_\mathfrak{p}\localiz{O}{S_0}$ is principal.
In other words, one can find $I$, a product of prime two-sided $O$-ideals among $\{P_\mathfrak{q} : \mathfrak{q} \in S_0\cap V^K_f\}$, such that $P_\mathfrak{p} I$ is principal; say that it is equal to $aO = Oa$ for some $a\in B^\times$.
It is clear that $a\in \Gamma_O(S\cup \Ram(B))$, and $v_\mathfrak{p}(a) = \frac{1}{2}$ whereas $v_\mathfrak{q}(a) = 0$ for all $\mathfrak{q}\in \Ram(B)\backslash S_0$ different from $\mathfrak{p}$.
The map is consequently surjective.
\end{proof}

Recall from \Cref{subsec:TwoSidedIdeals} that the map $I \mapsto [I]$, sending a two-sided $O$-ideal to its right class, induces a bijection
\[\Idl(O)/\PIdl(O) \rightarrow \{ [I]\in \ClR(O) \st O_\mathrm{L}(I) \simeq O \}.\]
Composing its inverse with the natural quotient map, one gets a surjective map
\[\psi : \{I \in \Idl_{\mathrm{R}}(O) \st O_{\mathrm{L}}(I) \simeq O \} \rightarrow \Idl(O)/\PIdl(O),\]
where $\Idl_{\mathrm{R}}(O)$ is the set of (fractional) invertible right $O$-ideals.

Note that if $I$ is a right $O$-ideal with left order isomorphic to $O$, and if $J$ is a two-sided $O$-ideal, then $IJ$ is right $O$-ideal with the same left order as $I$, and one has
\[\psi(IJ) = \psi(I) \psi(J).\]
Also, $\psi(I)$ is trivial if and only if $I$ is principal.

Let $\mathfrak{p}\in V^K(O)$ and let $N(O,\mathfrak{p})$ be the set of integral right $O$-ideals of reduced norm $\mathfrak{p}$.
As explained in \Cref{subsec:GeomOfMaxOrders} in the preliminaries, the function sending a right $O$-ideal to its left order induces a bijection between $N(O,\mathfrak{p})$ and the set of $\mathfrak{p}$-neighbours of $O$.
Since $\mathfrak{p}\in V^K(O)$ by hypothesis, for all $I$ in $N(O,\mathfrak{p})$, the left order of $I$ has the same type as $O$.
One can therefore consider the subgroup
\[\mathcal{G}(O,\mathfrak{p}) := \langle \psi(I) : I \in N(O,\mathfrak{p}) \rangle \leqslant \Idl(O)/\PIdl(O).\]
For a subset $T\subseteq V^K(O)$, we define the subgroup
\[\mathcal{G}(O, T) := \langle \mathcal{G}(O, \mathfrak{p}) : \mathfrak{p} \in T\rangle \leqslant \Idl(O)/\PIdl(O). \]

Recall now the map $P \mapsto P\cap R$ between the prime two-sided $O$-ideals and the prime $R$-ideals, which is a bijection.
For $\mathfrak{p}$ a prime $R$-ideal, write $P_{\mathfrak{p}}$ for the unique prime two-sided $O$-ideal such that $P_{\mathfrak{p}} \cap R = \mathfrak{p}$.
One has $P_{\mathfrak{p}}^2 = \mathfrak{p} O$ if $\mathfrak{p}\in \Ram(B)$, whereas $P_{\mathfrak{p}} = \mathfrak{p} O$ otherwise.

\begin{theo}\label{thm:GammaSTrans}
The group $\Gamma_O(S)$ acts vertex-transitively on the product of trees $\prod_{v\in S\backslash S_0} \mathcal{T}_v$ if and only if 
\begin{enumerate}
\item[(1)] $S\backslash S_0 \subseteq V^K(O)$ and
\item[(2)] $\mathcal{G}(O,S\backslash S_0) \leqslant \langle \PIdl(O), \Idl(R), \langle P_\mathfrak{p} : \mathfrak{p} \in S_0 \cap V^K_f \rangle\rangle / \PIdl(O)$.
\end{enumerate}
\end{theo}

\begin{proof}
Write $H$ for the group on the right-hand side, and write $v_\mathfrak{p}$ for the valuation corresponding to a prime $R$-ideal $\mathfrak{p}$.
Note that the subgroup of $\Idl(O)$ generated by $\Idl(R)$ and $\langle P_\mathfrak{p} : \mathfrak{p} \in S_0 \cap V^K_f \rangle$ is the group of two-sided $O$-ideals whose $P_{\mathfrak{q}}$-valuation is even for all $\mathfrak{q}\in \Ram(B)\backslash S_0$.

Assume first that $\Gamma_O(S)$ acts vertex-transitively.
For every $\mathfrak{p}\in S\backslash S_0$ and every $\mathfrak{p}$-neighbour $O'$ of $O$, there exists an element in $\Gamma_O(S)$ sending $O$ to $O'$.
It follows that $O$ and $O'$ are conjugate and have the same type, and from there, $S\backslash S_0 \subseteq V^K(O)$.

Next, we prove the inclusion $\mathcal{G}(O,S\backslash S_0) \leqslant H$.
Let $\mathfrak{p}\in S\backslash S_0$, and let $I\in N(O,\mathfrak{p})$; we have to show that $\psi(I)\in H$.
The left order of $I$ is a $\mathfrak{p}$-neighbour $O'$ of $O$ so, by hypothesis, there exists $g\in \Gamma_O(S)$ such that $g\cdot O = O'$.
If $a\in B^\times$ is a representative of $g$, then $O'$ is the left order of the right $O$-ideal $a O$.
Since $I$ and $a O$ share the same left order, the product $J := (a O)^{-1} I$ is well defined and is a two-sided $O$-ideal.

For all prime $\mathfrak{q}\in \Ram(B)\backslash S_0$, one has $g \in \PGL_1(O_\mathfrak{q})$ by hypothesis, so there exists a scalar $\lambda_{\mathfrak{q}} \in K_{\mathfrak{q}}^\times$ such that $\lambda_{\mathfrak{q}} a\in O_\mathfrak{q}^\times$.
It follows that the two-sided $O_\mathfrak{q}$-ideal $J O_\mathfrak{q}$ is equal to $\lambda_\mathfrak{q} O_\mathfrak{q}$ and is therefore a power of $\mathfrak{q} O_\mathfrak{q} = (P_\mathfrak{q} O_\mathfrak{q})^2$, so that $P_\mathfrak{q}$-valuation of $J$ is even.
The ideal $J$ then belongs to $\langle\Idl(R), \langle P_\mathfrak{q} : \mathfrak{q} \in S_0 \cap V^K_f \rangle\rangle$.
One concludes that 
\[\psi(I) = \psi(aO) \psi(J) = \psi(J)\in H.\]

Reciprocally, assume that $S\backslash S_0\subseteq V^K(O)$ and $\mathcal{G}(O,S\backslash S_0) \leqslant H$.
Let $O'$ be a $\mathfrak{p}$-neighbour of $O$ with $\mathfrak{p}\in S\backslash S_0$, and let $I$ be the ideal in $N(O,\mathfrak{p})$ whose left order is $O'$.
By hypothesis, $\psi(I) \in H$.
It follows that there exists a two-sided $O$-ideal $J \in \langle \Idl(R), \langle P_\mathfrak{q} : \mathfrak{q} \in S_0 \cap V^K_f \rangle\rangle$ such that $\psi(I) = \psi(J)$.
The product $I J^{-1}$ is well defined, and is a right $O$-ideal whose left order is $O'$.
Moreover, $\psi(I J^{-1})$ is trivial, meaning that $I J^{-1}$ is principal, say $I J^{-1} = a O$ for some $a\in B^\times$.
We write $g \in \PGL_1(B)$ for the class of $a$ modulo $K^\times$.
One has $g\cdot O = a O a^{-1} = O'$, so it remains to show that $g\in \Gamma_O(S)$.

Locally for any prime $R$-ideal $\mathfrak{q}$ different from $\mathfrak{p}$, one has $g\cdot O_{\mathfrak{q}} = O_{\mathfrak{q}}$.
If $q\notin \Ram(B)$, then $O_{\mathfrak{q}}$ is a vertex of the Bruhat-Tits tree at $\mathfrak{q}$ and its stabiliser for the action of $\PGL_1(B_\mathfrak{q})$ is $\PGL_1(O_\mathfrak{q})$, so $g\in \PGL_1(O_\mathfrak{q})$.
Now assume that $q\in \Ram(B)\backslash S_0$.
The ideal $a O_\mathfrak{q}$ is two-sided and equal to $J^{-1} O_\mathfrak{q}$, and it has an even $P_\mathfrak{q} O_\mathfrak{q}$-valuation; in other words, $a O_\mathfrak{q}$ is a power of $\mathfrak{q} O_\mathfrak{q}$.
As $\mathfrak{q} R_{\mathfrak{q}}$ is a principal $R_{\mathfrak{q}}$-ideal, $a$ is then an element of $O_\mathfrak{q}^\times$ up to a scalar multiple in $K_\mathfrak{q}^\times$, so $g\in \PGL_1(O_{\mathfrak{q}})$.
In conclusion, one has $g\in \Gamma_O(S_0 \cup \{\mathfrak{p}\}) \leqslant \Gamma_O(S)$, as wanted.
\end{proof}

\begin{corollary}\label{cor:GammaSTransIfClass1}
If $\Idl(\localiz{O}{S_0}) = \PIdl(\localiz{O}{S_0})$, then $\Gamma_O(S)$ acts vertex-transitively on the product of trees $\prod_{v\in S\backslash S_0} \mathcal{T}_v$ if and only if $S\backslash S_0 \subseteq V^K(O)$.
\end{corollary}

\begin{proof}
The condition $\Idl(\localiz{O}{S_0}) = \PIdl(\localiz{O}{S_0})$ is equivalent to the condition that $\Idl(O)$ is generated by $\PIdl(O)$ and $\langle P_\mathfrak{p} : \mathfrak{p} \in S_0 \cap V^K_f \rangle$.
The group on the right-hand side of (2) in \Cref{thm:GammaSTrans} is then the full group $\Idl(O)/\PIdl(O)$, and $\mathcal{G}(O,S\backslash S_0)$ is naturally a subgroup of it, whenever $S\backslash S_0 \subseteq V^K(O)$.
The result immediately follows.
\end{proof}

\begin{corollary}\label{cor:GammaSTransIfSContainsRam}
If $\Ram(B)\subseteq S$ (i.e.\ if $S_0 = \Ram(B)$), then $\Gamma_O(S)$ acts vertex-transitively on the product of trees $\prod_{v\in S\backslash S_0} \mathcal{T}_v$ if and only if $S\backslash S_0 \subseteq V^K(O)$.
\end{corollary}

\begin{proof}
Indeed, in this case, the group generated by $\Idl(R)$ and $\langle P_\mathfrak{p} : \mathfrak{p} \in S_0 \cap V^K_f \rangle$ is $\Idl(O)$.
The subgroup condition (2) in \Cref{thm:GammaSTrans} then becomes superfluous and the result follows.
\end{proof}

One may be interested in the difference between the lattice $\Gamma_S(O)$ and the ``naive'' quotient group \[\PGL_1(\localiz{O}{S}):= \localiz{O}{S}^\times/\localiz{R}{S}^\times.\]
The following proposition gives an element of answer.

\begin{proposition}\label{prop:ExactSeqPGL1GammaS}
One has an exact sequence
\[1 \rightarrow \PGL_1(\localiz{O}{S}) \rightarrow \Gamma_{O}(S) \rightarrow \PIdl(\localiz{O}{S})/\PIdl(\localiz{R}{S}) \rightarrow 1.\]
In particular, $\PGL_1(\localiz{O}{S})$ is an arithmetic group and a cocompact lattice in $\prod_{v\in S\backslash S_0} \PGL_2(K_v)$ if $S_0$ is strictly contained in $S$.
\end{proposition}

\begin{proof}
If $g$ is an element of $\Gamma_{O}(S)$, it stabilises $\localiz{O}{S}$ and $a \localiz{O}{S}$ is then a two-sided $\localiz{O}{S}$-ideal for any representative $a\in B^\times$ of $g$.
The class of $a \localiz{O}{S}$ modulo the principal $\localiz{R}{S}$-ideals does not depend on the representative $a$ of $g$.
Also, any principal two-sided $\localiz{O}{S}$-ideal is of the form $a\localiz{O}{S}$ for a certain $a$ normalising $\localiz{O}{S}$.
We therefore have a surjective group morphism
\[\phi : \Gamma_{O}(S) \rightarrow \PIdl(\localiz{O}{S})/\PIdl(\localiz{R}{S}).\]
It remains to show that its kernel is $\PGL_1(\localiz{O}{S})$.
If $a\in \localiz{O}{S}^\times$, then $a\localiz{O}{S} = \localiz{O}{S}$, showing that $\PGL_1(\localiz{O}{S})$ is in the kernel.
Reciprocally, if $a\in B^\times$ is a representative of an element $g\in \Gamma_{O}(S)$ in the kernel, there exists a principal $\localiz{R}{S}$-ideal $(\lambda)$ such that $a \localiz{O}{S} = \lambda \localiz{O}{S}$, i.e.\ $\lambda^{-1} a \in \localiz{O}{S}^\times$, as wanted.

By the exact sequence, the group $\PGL_1(\localiz{O}{S})$ is a finite-index subgroup of $\Gamma_O(S)$.
Since $\Gamma_O(S)$ is an arithmetic group and a cocompact lattice in $\prod_{v\in S\backslash S_0} \PGL_2(K_v)$ by \Cref{lemma:GammaSCocompactLattice}, the same is true for $\PGL_1(\localiz{O}{S})$.
\end{proof}

\begin{remark}\label{rem:PIdlOPIdlRsubgrpofClR}
If $\Ram(B) \subseteq S$, then $\Idl(\localiz{O}{S})$ and $\Idl(\localiz{R}{S})$ are the same group, and
\[\Cl(\localiz{R}{S})= \Idl(\localiz{O}{S})/\PIdl(\localiz{R}{S}).\]
In particular, the quotient $\PIdl(\localiz{O}{S})/\PIdl(\localiz{R}{S})$ in \Cref{prop:ExactSeqPGL1GammaS} is a subgroup of the class group of $\localiz{R}{S}$.
\end{remark}

Except for a certain maximal order in $(-1,-1 \st \QQ(\sqrt{15}))$, the ring $\localiz{R}{\Ram(B)}$ is a principal ideal domain for all the other currently known examples of maximal orders $O$ where $V^K(O)$ is non-empty, which are described in the two next subsections.
The following result is therefore useful. 

\begin{corollary}\label{cor:RS0PrincipalDomainPGL1eqGammS}
Assume that $S_0 = \Ram(B)$ and that $\localiz{R}{S_0}$ is a principal ideal domain.
Then $\PGL_1(\localiz{O}{S}) = \Gamma_O(S)$ and $\PGL_1(\localiz{O}{S_0}) = \Gamma_{O}(S_0)$.

In particular, the group $\PGL_1(\localiz{O}{S})$ acts vertex-transitively on the product of trees $\prod_{v\in S\backslash S_0} \mathcal{T}_v$, with vertex-stabilisers conjugate to $\PGL_1(\localiz{O}{S_0})$, if and only if $S\backslash S_0\subseteq V^K(O)$.
\end{corollary}

\begin{proof}
Indeed, if $\localiz{R}{S_0}$ is a principal ideal domain, so is $\localiz{R}{S}$, and their class groups are trivial.
By the exact sequence in \Cref{prop:ExactSeqPGL1GammaS} and by \Cref{rem:PIdlOPIdlRsubgrpofClR}, the groups $\Gamma_O(S)$ and $\PGL_1(\localiz{O}{S})$ are the same, and similarly for and $\Gamma_O(S_0)$ and $\PGL_1(\localiz{O}{S_0})$.
One concludes using \Cref{cor:GammaSTransIfSContainsRam}.
\end{proof}

\subsection{Transitive lattices from orders of type number 1}\label{subsec:TransLattType1}

Maximal orders of type number $1$ are of main interest since by definition, they are the ones for which $V^K(O) = V^K \backslash \Ram(B)$ by \Cref{rem:VKOFullIFFTypeNumberOne}, so that for any choice of $S$ containing $\Ram(B)$, the group $\Gamma_O(S)$ is vertex-transitive on the product of trees $\prod_{S\backslash S_0} \mathcal{T}_v$ by \Cref{cor:GammaSTransIfSContainsRam}.
A complete classification of such orders with $K$ a number field is available.

\begin{theo}[Kirschmer–Lorch \cite{KirschmerLorch16}]\label{thm:classiftype1}
Up to isomorphism, there are exactly 105 totally definite quaternion algebras over number fields of type number $1$.
\end{theo}

One can find in the Appendix the list of those quaternion algebras, together with relevant information for finding vertex-transitive subgroups. 
By checking all of them, one obtains the following fact.

\begin{corollary}\label{cor:Btymenum1thenPrincipalIdealDom}
If $B$ has type number one and $K$ is a number field, then $\localiz{R}{\Ram(B)}$ is a principal ideal domain.

In particular, for any $S$ with $\Ram(B) \subseteq S$, the groups $\Gamma_O(S)$ and $\PGL_1(\localiz{O}{S})$ are equal and act vertex-transitively on the product of trees $\prod_{S\backslash S_0} \mathcal{T}_v$, with vertex-stabilisers conjugate to $\Gamma_O(S_0) = \PGL_1(\localiz{O}{S_0})$ by \Cref{cor:RS0PrincipalDomainPGL1eqGammS}.
\end{corollary}

Actually, the ring of integers $R$ is principal too for all the examples except one, namely $B = (-75-19\sqrt{15}, -4-\sqrt{15} \st \QQ(\sqrt{15}))$.
Indeed, in this case $R= \ZZ[\sqrt{15}]$ has a class group of order $2$.

In the positive characteristic case, Gekeler used the theory of Drinfeld modules to give an explicit formula for the type number $|\Typ(O)|$ for certain quaternion algebras with exactly two ramified places.

\begin{theo}[Gekeler \cite{Gekeler91}] \label{thm:poschartypenum}
Assume that $K = \FF_q(t)$ and $v_\infty = -\deg$. 
Let $f\in \FF_q[t]$ be an irreducible monic polynomial.
Then the quaternion algebra $B$ over $K$ with $\Ram(B) = \{v_\infty, v_f\}$ has type number one if and only if $\deg(f)\in \{1,2\}$.
\end{theo}

\begin{proof}
Let $d = \deg(f)$ and let $O$ be a maximal order in $B$.
By \cite[(4.7)]{Gekeler91}, it is equal to $1$ if $d\in \{1,2\}$. 
In general, we have the following inequalities:
\[
|\Typ(O)|  > 
\begin{cases} \frac{1}{2}\left( \frac{q^d-q}{q^2-1} + 1 \right) & \text{if $d$ is odd,} \\
 \frac{1}{2} \frac{q^d-1}{q^2-1} & \text{if $d$ is even.} 
\end{cases}\]
It follows that $|\Typ(O)| > 1$ whenever $d > 2$.
\end{proof}

\subsection{Transitive lattices from orders of type number \texorpdfstring{$>1$}{>1}}\label{subsec:TransLattTypeGt1}

As we will see, maximal orders of type number $1$ are not the only ones for which $V^K(O)$ is non-empty.
We start with those of type number $2$.

\begin{lemma} \label{lemma:typenum2sametypO}
Let $O$ and $O'$ be non-isomorphic maximal orders in a quaternion algebra $B$ of type number $2$. 
Then $V^K(O) = V^K(O')$.
\end{lemma}
\begin{proof}
Let $v\in V^K\backslash \Ram(B)$ not belonging to $V^K(O')$.
By definition, there exists a $v$-neighbour $U$ of $O'$ not isomorphic to $O'$.
Since $B$ has type number two, $U$ is isomorphic to $O$, implying that $V^K(U) = V^K(O)$ by \Cref{lemma:ordersisosameVKO}.
Since $O'$ is a $v$-neighbour of $U$ not isomorphic to $U$, the place $v$ does not belong to $V^K(U)$, so $V^K(O)\subseteq V^K(O')$.
By symmetry, $V^K(O) = V^K(O')$.
\end{proof}

In the above-mentioned paper \cite{KirschmerLorch16}, Kirschmer and Lorch also give a classification of orders in quaternion algebras over number fields having type number two.

\begin{theo}[Kirschmer–Lorch \cite{KirschmerLorch16}]\label{thm:classiftype2}
Up to isomorphism, there are exactly 246 totally definite quaternion algebras over number fields of type number $2$.
\end{theo}

From that classification, one gets the following result.

\begin{proposition}\label{prop:VKOnonemptyType2}
There are exactly $3$ totally definite quaternion algebras over number fields of type number $2$ for which $V^K(O)$ is non-empty for one (hence all) of its maximal orders.
They are
\[(-1,-1 \st \QQ(\sqrt{3})), (-1,-1 \st \QQ(\sqrt{21})) \text{ and } (-1,-1 \st \QQ(\zeta_{15}+\zeta_{15}^{-1})),\] where $\zeta_{15}$ is a $15$th primitive root of unity.
\end{proposition}

\begin{proof}
By \Cref{thm:VKOnonempty}, one has $V^K(O) = (V^K\backslash \Ram(B)) \cap \mathcal{I}^+_B(R)$ if $V^K(O)$ is non-empty.
One enumerates the $246$ quaternion algebras in \Cref{thm:classiftype2}, and take any (small) prime ideal $\mathfrak{p}$ in $(V^K\backslash \Ram(B)) \cap \mathcal{I}^+_B(R)$.
One then checks whether all the $\mathfrak{p}$-neighbours of $O$ are isomorphic to $O$ or not.
It turns out that the 3 quaternions algebras above are the only ones for which $V^K(O)$ is non-empty among the whole list.
This was checked using Magma \cite{Magma}.
\end{proof}

One has $\Ram(B) = V^K_\infty$ for all those examples, so $S_0 = V^K_\infty$.
Moreover, the class groups of $\QQ(\sqrt{3})$, $\QQ(\sqrt{21})$ and $\QQ(\zeta_{15}+\zeta_{15}^{-1})$ are all trivial.
By \Cref{cor:RS0PrincipalDomainPGL1eqGammS}, for each of the six orders, the group $\Gamma_O(S) = \PGL_1(\localiz{O}{S})$ acts vertex-transitively on the product of trees if and only if $S\backslash S_0 \subseteq V^K(O)$, with vertex-stabilisers conjugate to $\Gamma_O(S_0) = \PGL_1(O)$.
Moreover, $\Gamma_O(S_0)$ is isomorphic to
\[D_{12}\text{ and }\Sym(4) \text{ for the orders of } (-1,-1 \st \QQ(\sqrt{3})),\]
\[D_{6}\text{ and }\Alt(4) \text{ for the orders of } (-1,-1 \st \QQ(\sqrt{21})),\]
\[D_{30}\text{ and }\Alt(5) \text{ for the orders of } (-1,-1 \st \QQ(\zeta_{15}+\zeta_{15}^{-1})).\]

Note that $\Idl(O) = \PIdl(O)$ for all those orders, and the narrow class groups of the base fields all have order $2$.

\vspace{5mm}

To our knowledge, no classification is available yet for orders of type number $>2$ in the number field case.
Numerical experiments make us think that those maximal orders of type number $>2$ and $V^K(O)$ non-empty are relatively rare: we only found $9$ examples, but we expect more to exist.
We do not know whether the number of such orders is finite for each characteristic.

An example of obstruction to the existence of such orders is given by \Cref{cor:NarrowClGrpTrivType1}.
Indeed, if the field has a trivial narrow class group (e.g.\  $K = \QQ$), then a maximal quaternion order $O$ with $V^K(O)$ non-empty is necessarily of type number $1$.

Using Magma \cite{Magma}, we tested quaternion algebras of the form 
\[(-1,-1\st K)\]
for several totally real number fields, and we found other interesting examples.
Namely, the fields $K = \QQ[X]/(P(X))$ for $P(X)$ an irreducible polynomial over $\QQ$ in the following list:
\[X^2-6, X^2-7, X^2-15, X^2-69,X^4-5X^2+5, X^4-4X^2+1,\] 
\[X^4-6X^2-4X+2, X^4-2X^3-3X^2+4X+1, X^4+2X^3-7X^2-8X+1.\]
The type numbers of the respective quaternion algebras are $3,3,8,5,3,3,4,4,4$.
For each of those, one maximal order $O$ (and only one up to conjugation) is such that $V^K(O)$ is non-empty, and one has $\Ram(B) = V^K_\infty$ and $\Idl(O) = \PIdl(O)$.
Note that \Cref{lemma:typenum2sametypO} is consequently no longer true if $B$ is of type number $>2$.

By \Cref{cor:GammaSTransIfClass1}, the group $\Gamma_O(S)$ acts vertex-transitively on the product of trees if and only if $S\backslash S_0 \subseteq V^K(O)$, with stabilisers conjugate to $\Gamma_O(S_0)$ by \Cref{lemma:VertexStabGammaS}.
Moreover, the stabilisers are respectively isomorphic to the following groups:
\[D_4, \Sym(3), C_2\times C_2, C_2, \Sym(4), D_8,\]
\[\Sym(3),\Sym(3),D_5.\]

The class groups of the fields are all trivial, except for $\QQ(\sqrt{15})$, where it has order $2$.
The narrow class groups of the fields all have order $2$, except for $\QQ(\sqrt{15})$, where it is isomorphic to $C_2\times C_2$.
A description of the set $V^K(O)$ is given as follows.

\begin{corollary}\label{cor:VKOtypeGt1}
If $O$ is a maximal order in one of the quaternion algebras of type number $2$ in \Cref{prop:VKOnonemptyType2}, or one of the $9$ maximal orders of type number $>2$ given just above, then
\[V^K(O) = V^K_f \cap \PIdl^+(R)
= \{ \mathfrak{p} \in V^K_f \st  \mathfrak{p} \text{ is a totally positive principal $R$-ideal}\}.\]
\end{corollary}

\begin{proof}
For each of those quaternion algebras, $\Ram(B) = V^K_\infty$.
By definition, the group $\mathcal{I}_B^+(R)$ in \Cref{thm:VKOnonempty} is then the subgroup of $\Idl(R)$ generated by $\PIdl^+(R)$ and the squares of the prime ideals in $V^K_f$.
But as the narrow class group of $K$ are all elementary abelian $2$-groups (either $C_2$ or $C_2\times C_2$), the square of any prime ideal is in $\PIdl^+(R)$, so that $\mathcal{I}_B^+(R) = \PIdl^+(R)$ and \Cref{thm:VKOnonempty} concludes.
\end{proof}

\subsection{Constructing families of vertex-transitive subgroups}\label{subsec:FamOfVertTransSubGrps}

The goal of this subsection is to describe how to find families of vertex-transitive subgroups of $\Gamma_O(S)$, with a main focus on BMW subgroups, i.e.\ those acting regularly.

We start by presenting two elementary lemmas about transitive group actions.

\begin{lemma} \label{lemma:GtransHtrans}
Let $G$ be a group acting transitively on a non-empty set $X$ and let $G_0$ be the stabiliser of an element $x\in X$.
Let $H$ be a finite index subgroup of $G$, and let $H_0$ be the stabiliser of $x$ for the restricted action on $X$.
Then $H$ is transitive on $X$ if and only if the index of $H$ in $G$ coincide with the index of $H_0$ in $G_0$.
\end{lemma}

\begin{proof}
We have
\[|G/H| \geq |(G_0 H)/H| = |G_0/(G_0\cap H)| = |G_0/H_0|, \]
and equality holds if and only if $G = G_0 H$, if and only if $H$ is transitive on $X$.
\end{proof}

\begin{lemma}\label{lemma:Htransfromquotient}
Let $G$ be a group acting transitively on a non-empty set $X$ and let $G_0$ be the stabiliser of an element $x\in X$.
Let $f : G \rightarrow Q$ be a surjective group morphism and let $A$ be a subgroup of $Q$.
If 
\[[Q : A] = [f(G_0) : A\cap f(G_0)],\]
 then $f^{-1}(A)$ is transitive on $X$.
\end{lemma}

\begin{proof}
Let $n$ be the index $[Q : A] = [f(G_0) : A\cap f(G_0)]$.
The stabiliser of $x$ for the action of $f^{-1}(A)$ on $X$ is $f^{-1}(A)\cap G_0$.
Since $f$ is surjective, the subgroup $f^{-1}(A)$ has index $n$ in $G$.
Note that, by basic set theory, one has $f(f^{-1}(A) \cap G_0) = A \cap f(G_0)$.

On one hand, 
\[[G_0 : f^{-1}(A)\cap G_0] \leq [G : f^{-1}(A)] = n.\]
On the other hand,
\begin{align*}
[G_0 : f^{-1}(A)\cap G_0] &\geq [f(G_0) : f(f^{-1}(A)\cap G_0)] \\
& = [f(G_0) : A \cap f(G_0)] = n.
\end{align*}
We therefore have $[G_0 : f^{-1}(A)\cap G_0] = n$, showing the transitivity of $f^{-1}(A)$ on $X$ with help of \Cref{lemma:GtransHtrans}. 
\end{proof}

As mentioned earlier, in all the currently known examples of maximal quaternionic orders $O$ with $V^K(O)$ non-empty, we have $\Idl(\localiz{O}{\Ram(B)}) = \PIdl(\localiz{O}{\Ram(B)})$.
By \Cref{cor:GammaSTransIfSContainsRam}, $\Gamma_O(S)$ then acts vertex-transitively with stabilisers conjugate to $\Gamma_O(S_0)$ if and only if $S\backslash S_0 \subseteq V^K(O)$, upon the additional condition that $S_0 = \Ram(B)$.

For all those known examples, except for one maximal order of $(-1,-1 \st \QQ(\sqrt{15}))$, the ring $\localiz{R}{\Ram(B)}$ is a principal ideal domain; the groups $\PGL_1(\localiz{O}{S})$ and $\Gamma_O(S)$ are therefore the same by \Cref{cor:RS0PrincipalDomainPGL1eqGammS}, and the respective vertex-stabilisers $\PGL_1(\localiz{O}{S_0})$ and $\Gamma_O(S_0)$ are also equal.
One might then prefer to work with the definition of $\PGL_1(\localiz{O}{S})$ instead of the adélic description of $\Gamma_O(S)$.
For $(-1,-1 \st \QQ(\sqrt{15}))$, one has $\Ram(B) = V^K_\infty = S_0$ and one can see from the exact sequence in \Cref{prop:ExactSeqPGL1GammaS}, that $\PGL_1(\localiz{O}{S})$ and $\PGL_1(O)$ have index $2$ in $\Gamma_O(S)$ and $\Gamma_O(S_0)$, respectively.

The two groups involved in the subgroup condition (2) in \Cref{thm:GammaSTrans} may vary in function of $S$.
It is therefore useful for proving that $\Gamma_O(S)$ is transitive for a given $S$, but not practical for finding families where $S$ varies and $S_0 = S\cap \Ram(B)$ is fixed.
Because of that, we will usually assume that
\[\Idl(\localiz{O}{S_0}) = \PIdl(\localiz{O}{S_0}),\]
which is a condition that depends only $S_0$, and ensures the vertex-transitivity of $\Gamma_O(S)$ if $S\backslash S_0 \subseteq V^K(O)$ by \Cref{cor:GammaSTransIfClass1}.

The next idea is to find congruence subgroups that are also vertex-transitive.
This is where \Cref{lemma:Htransfromquotient} comes into play.
We will assume from now that $S$ strictly contains $S_0$.

Let $\mathfrak{m}$ be a non-zero ideal of $R$, coprime with the prime ideals in $S\cap V^K_f$.
Say that 
\[\mathfrak{m} = \prod_{i=1}^n \mathfrak{p}_i^{e_i},\]
where $\mathfrak{p}_i$ are pairwise coprime prime ideals.
Let us write
\[\pi_\mathfrak{m} : \PGL_{1,O}(\mathcal{A}(S)) \rightarrow \PGL_1(O/\mathfrak{m}O) \simeq \prod_{i = 1}^n \PGL_1(O_{\mathfrak{p}_i}/\mathfrak{p}_i^{e_i} O_{\mathfrak{p}_i})\]
for the reduction map modulo $\mathfrak{m}$ of the $S$-integral adèle group.
We will use the notation
\[\Gamma_{O/\mathfrak{m}}(S) := \pi_{\mathfrak{m}}(\Gamma_{O}(S)),\]
As a corollary of the Strong Approximation Theorem for $\SL_1(B)$, since $S$ includes a place not ramified in $B$, the reduction $\Gamma_{O/\mathfrak{m}}(S)$ contains at least
\[\PSL_1(O/\mathfrak{m}O) \simeq \prod_{i = 1}^n \PSL_1(O_{\mathfrak{p}_i}/\mathfrak{p}_i^{e_i} O_{\mathfrak{p}_i});\]
see \cite[Corollary 28.5.14]{Voight21}.

Then, one looks for a subgroup $A$ of $\Gamma_{O/\mathfrak{m}}(S)$ such that 
\[[\Gamma_{O/\mathfrak{m}}(S) : A] = [\Gamma_{O/\mathfrak{m}}(S_0) : A\cap \Gamma_{O/\mathfrak{m}}(S_0)].\]
\Cref{lemma:Htransfromquotient} ensures that the action of $\pi_\mathfrak{m}^{-1}(A) \cap \Gamma_O(S)$ is transitive, with stabilisers conjugate to $\pi_\mathfrak{m}^{-1}(A)\cap \Gamma_O(S_0)$.

If one is interested in a regular action, one can try to apply the given recipe with the additional restriction that $[\Gamma_{O/\mathfrak{m}}(S) : A]$ is also the order of $\Gamma_{O}(S_0)$.
In other words, the morphism $\pi_\mathfrak{m}$ is injective on $\Gamma_{O}(S_0)$, i.e.\ $\Gamma_{O}(S_0) \simeq \Gamma_{O/\mathfrak{m}}(S_0)$, and $A\cap \Gamma_{O/\mathfrak{m}}(S_0) = 1$.

\begin{remark} \label{rem:strgapproxthmimagecontainsPSL}
If the ideal $\mathfrak{m}$ is prime, the algebra $O/\mathfrak{m} O$ is isomorphic to the matrix algebra $M_2(k_\mathfrak{m})$, where $k_\mathfrak{m}$ is the residue field of $K_\mathfrak{m}$.
One then has $\PSL_1(O/\mathfrak{m} O) \simeq \PSL_2(k_\mathfrak{m})$ and $\PGL_1(O/\mathfrak{m} O) \simeq \PGL_2(k_\mathfrak{m})$, as well as the exact sequence
\[1 \rightarrow \PSL_1(O/\mathfrak{m} O) \rightarrow \PGL_1(O/\mathfrak{m} O) \xrightarrow{\overline{\Nrd}} k_\mathfrak{m}^\times/(k_\mathfrak{m}^\times)^2 \rightarrow 1.\]
It follows that $\Gamma_{O/\mathfrak{m}}$ is $\PSL_1(O/\mathfrak{m} O)$ if for all $x\in \Gamma_O(S)$, $\Nrd(x)$ is a square modulo $\mathfrak{m}$, and $\PGL_1(O/\mathfrak{m} O)$ otherwise.

\end{remark}

\section{Group presentations}\label{sec:GrpPres}

In this section, we describe a method for finding a finite presentation of a group $\Gamma \leqslant \prod_{k=1}^n \Aut(\mathcal{T}_k)$ acting transitively (on the left) on the vertices of the Cartesian product of locally finite trees $\TT = \prod_{k=1}^n \mathcal{T}_k$, with finite vertex-stabilisers.
Our results are similar to \cite[Algorithm 2.6]{infauthors15}, but slightly more general to apply in our setting.
It is also a special case of the more general results in \cite{Macbeath64} or \cite{Brown84} which apply to groups acting by isometries on path-connected and simply connected spaces, but we take advantage of the structure of a product of trees.

As a graph, a product of trees is not simply-connected, except in trivial cases.
However, a tree can also be seen as a $1$-dimensional cube complex, and a product of $n$ trees can then be seen as an $n$-dimensional cube complex, which is simply connected and path connected.
The action of $\Gamma$ on the product of trees seen as a graph extends to a homeomorphic action on the product of trees seen as a cube complex.

Note that for all $1\leq k \leq n$, since the group $\Gamma$ acts transitively on the tree $\mathcal{T}_k$ through the projection on the $k$th coordinate, the tree is necessarily regular, say of degree $d_k$.
To avoid any triviality, we will assume that $d_k$ is at least $1$.

Fix a vertex $x = (x_k)_{1\leq k \leq n}$, and let $\Gamma_0$ be the stabiliser of $x$ for the action of $\Gamma$.
For all $1\leq i\leq n$, we define $\Gamma_i$ to be the subgroup of $\Gamma$ fixing $x_k$ for all $1\leq k \leq n$ different from $i$. 
Obviously, $\Gamma_0\leqslant \Gamma_i$.

For all $1\leq i\leq n$ and $1\leq j\leq d_i$, we write $x^{i,j} = (x^{i,j}_k)_{1\leq k \leq n} \in \TT$ for the vertex such that $x_k^{i,j} = x_k$ if $k\neq i$ and $x_i^{i,j}$ is the $j$-th neighbour of $x_i$ in $\mathcal{T}_i$.
Since the action of $\Gamma$ is transitive, there exists $g^{i,j}\in \Gamma$ for which $g^{i,j}\cdot x = x^{i,j}$.
By definition, one has $g^{i,j}\in \Gamma_i$.
Let $S_i = \{g^{i,j} : 1\leq j \leq d_i\}$.

The following theorem, due to A. M. MacBeath in \cite{Macbeath64}, gives a presentation of the group $\Gamma$.
This result does not take advantage of the geometrical structure of $\TT$, and the given presentation is therefore not very suitable.
Nevertheless, it will be of great help in our quest of a more natural presentation (cf.\ \Cref{thm:PresentationGamma}).

\begin{theo}[A. M. Macbeath, 1964]\label{thm:Macbeath}
Let $X$ be a non-empty path-connected and simply-connected topological space.
Let $\Gamma$ be a group acting on $X$ by homeomorphisms, and let $U\subseteq X$ be a path-connected open subset such that $X = \Gamma \cdot U$.
Define 
\[S = \{\gamma \in \Gamma \st \gamma\cdot U \cap U \neq \emptyset \},\]
and let $\overline{S}$ be a set of symbols together with a bijection $S\rightarrow \overline{S} : s \mapsto \overline{s}$.
Then $\Gamma$ is generated by $S$, and a presentation of $\Gamma$ is given by $\langle \overline{S} \st \mathcal{Q} \rangle$, where
\[\mathcal{Q} = \{\overline{s} \overline{t} (\overline{st})^{-1} : s,t \in S \st st\in S \} \subseteq \langle \overline{S} \rangle.\]
\end{theo}

In our setting, the space $X$ will be the product of trees $\TT$ seen as a cube complex.
For each $1\leq k \leq n$, put the usual metric $d$ on the geometric realisation of the tree $\mathcal{T}_k$, such that any two neighbouring vertices are at distance $1$ from each other.
Then put the $\infty$-metric $d = d_\infty$ on the product of trees.
In other words, if $u = (u_k)_{1\leq k \leq n}$ and $v = (v_k)_{1\leq k \leq n}$ are points in the cube complex $\TT$, then
\[d(u,v) = \max_{1\leq k \leq n} d(u_k,v_k).\]
Note that $\Gamma$ acts isometrically on $\TT$.
Let $U$ be the open ball in $\TT$ centred at $x$ and of radius $2/3$.
Since any point in $\TT$ is at distance $\leqslant 1/2$ from a nearest vertex, the transitivity of $\Gamma$ on the vertices implies that $\TT = \Gamma \cdot U$.
As in \Cref{thm:Macbeath}, write $S$ for the set $\{\gamma \in \Gamma \st \gamma\cdot U \cap U \neq \emptyset \}$.

\begin{lemma}\label{lemma:SetSElemDistLe1}
One has
\[S = \{\gamma \in \Gamma \st d(x, \gamma \cdot x) \leq 1\},\]
and the cardinal of $S$ is $|\Gamma_0|\cdot \prod_{k=1}^{n} (d_k + 1)$.
\end{lemma}

\begin{proof}
If $\gamma \in S$, let $u\in \gamma\cdot U \cap U $.
Then $d(x, \gamma\cdot x) \leqslant d(x, u) + d(u, \gamma\cdot x) < 4/3 $.
Since $d(x, \gamma\cdot x)$ is an integer, it is $\leq 1$.

Reciprocally, if $\gamma$ is such that $d(x, \gamma\cdot x) \leq 1$, then $d(x_k, \gamma \cdot x_k) \leq 1$ for all $1\leq k \leq n$.
Let $u_k\in T_k$ be the midpoint between the vertices $x_k$ and $g\cdot x_k$, which is at distance $\leq 1/2$ for each of them.
Then the point $u = (u_k)_{1\leq k \leq n}$ belongs both to $U$ and $\gamma \cdot U$, implying that their intersection is non-empty and $\gamma\in S$.

For each $1\leq k \leq n$, the number of vertices in $\mathcal{T}_k$ that are at distance $\leq 1$ from $x_k$ is $d_k + 1$.
It follows that the number of vertices at distance $\leq 1$ from $x$ is the product $\prod_{k=1}^{n} (d_k + 1)$.
Moreover, since $\Gamma$ is transitive on the vertex set and $\Gamma_0$ is the stabiliser of $x$, for each vertex $y$, there are exactly $|\Gamma_0|$ elements in $\Gamma$ sending $x$ to $y$, concluding the proof.
\end{proof}

\begin{proposition}\label{prop:SDecomposedAsProd}
The function
\[(S_1 \cup \{e\}) \times \cdots \times (S_n \cup \{e\}) \times \Gamma_0 \rightarrow S,\]
sending $(s_1, \dots, s_n, \gamma_0)$ to the product $s_1 \cdots s_n \gamma_0$ is a bijection.
\end{proposition}

\begin{proof}
Let $s\in S$.
By \Cref{lemma:SetSElemDistLe1}, one has $d(x, s\cdot x) \leq 1$, so $d(x_k, s\cdot x_k) \leq 1$ for all $1\leq k \leq n$.
Write $t_0$ for the element $s$.
Now assume that $t_i$ is defined for some $0\leq i < n$, and has the property that $t_i x_k = x_k$ for all $1\leq k \leq i$, and $d(x_k, t_i\cdot x_k) \leq 1$ for all $i+1\leq k \leq n$.
There is a unique $s_{i+1}\in S_{i+1}\cup \{e\}$ such that $s_{i+1} \cdot x_{i+1} = t_i \cdot x_{i+1}$.
Define $t_{i+1}$ to be $s_{i+1}^{-1} t_i$.
One has $t_{i+1} \cdot x_{i+1} = x_{i+1}$ by construction, then $t_{i+1} \cdot x_k = s_{i+1}^{-1} \cdot x_k = x_k$ for all $1\leq k \leq i$, and $d(x_k, t_{i+1}\cdot x_k) = d(s_{i+1} \cdot x_k, t_i\cdot x_k) = d(x_k, t_i\cdot x_k) \leq 1$ for all $i+2\leq k \leq n$.

The element $t_n$ is of the form $s_n^{-1} \cdots s_1^{-1} s$, where $s_k\in S_k\cup \{e\}$ for all $1\leq k \leq n$, and it stabilises $x$, so it equals an element $\gamma_0 \in \Gamma_0$.
One then has $s = s_1 \cdots s_n \gamma_0$.
Since the cardinals of the sets $(S_1 \cup \{e\}) \times \cdots \times (S_n \cup \{e\}) \times \Gamma_0 $ and $S$ coincide by \Cref{lemma:SetSElemDistLe1}, one concludes that the function in the statement is well defined and is a bijection.
\end{proof}

\begin{corollary}\label{cor:SwapGamma0Si}
Let $1 \leq i \leq n$.
There are unique functions $\sigma_i : \Gamma_0 \times S_i \rightarrow S_i$ and $\alpha_i : \Gamma_0 \times S_i \rightarrow \Gamma_0$ such that for all $\gamma_0\in \Gamma_0$ and $s_i \in S_i$,
\[\gamma_0 s_i = \sigma_i(\gamma_0, s_i) \alpha_i(\gamma_0, s_i). \]
\end{corollary}

\begin{proof}
Indeed, $\gamma_0 s_i \cdot x_k = x_k$ for all $1\leq k \leq n$ different from $i$, and $\gamma_0 s_i \cdot x_i$ is a neighbour of $x_i$, so that $\gamma_0 s_i \in S$.
By \Cref{prop:SDecomposedAsProd}, the element $\gamma_0 s_i$ is then uniquely writable in the form $s'_1\cdots s'_n \gamma'_0$ as in the proposition, and it is clear that $s'_k = e$ for all $k$ different from $i$.
\end{proof}

\begin{corollary}\label{cor:SwapSiSj}
Let $1\leq i,j \leq n$ with $i > j$.
There are unique functions $\lambda_{i,j} : S_i\times S_j \rightarrow S_j$,  $\mu_{i,j} : S_i\times S_j \rightarrow S_i$ and $\beta_{i,j} : S_i \times S_j \rightarrow \Gamma_0$ such that, for all $s_i\in S_i$ and $s_j\in S_j$,
\[s_i s_j = \lambda_{i,j}(s_i,s_j) \mu_{i,j}(s_i,s_j) \beta_{i,j}(s_i,s_j)\]
\end{corollary}

\begin{proof}
The proof is similar to that of \Cref{cor:SwapGamma0Si}, once you note that $s_i s_j \in S$.
\end{proof}

\begin{lemma}\label{lemma:InverseFuncSi}
Let $1\leq i \leq n$.
There are unique functions $\iota_i : S_i \rightarrow S_i$ and $\delta_i : S_i \rightarrow \Gamma_0$ such that, for all $s_i\in S_i$,
\[s_i \iota_i(s_i) = \delta_i(s_i).\]
Moreover, for all $s_i,t_i\in S_i$, one has $d(x, s_it_i \cdot x) \leq 1$ if and only if $t_i = \iota_i(s_i)$.
\end{lemma}

\begin{proof}
Let $s_i\in S_i$.
The set $\{s_i t_i \cdot x_i : t_i \in S_i\}$ is the set of neighbours of $s_i \cdot x_i$ in the tree $\mathcal{T}_i$.
Therefore, there is only one $t_i \in S_i$ such that $s_i t_i$ stabilises $x_i$, i.e.\ such that $s_i t_i \in \Gamma_0$, proving the existence and the uniqueness of $\iota_i$ and $\delta_i$.
Also, if $t_i\in S_i$ is different from $\iota_i(s_i)$, then $s_it_i\cdot x_i$ is at distance $2$ from $x_i$, proving that $d(x, s_it_i \cdot x) \leq 1$ if and only if $t_i = \iota_i(s_i)$. 
\end{proof}

\begin{remark}
Since the sets $S_1,\dots, S_n, \Gamma_0$ are pairwise disjoint, we can get rid of the subscripts of the function $\alpha_i, \beta_{i,j}, \delta_i, \sigma_i, \lambda_{i,j}, \mu_{i,j}, \iota_i$, and simply write $\alpha, \beta, \delta, \sigma, \lambda, \mu, \iota$ without any confusion.
\end{remark}

Let $\overline{\Gamma_0}$ and $\overline{S_k}$ for $1\leq k \leq n$ be sets of symbols, together with bijections
\[\Gamma_0 \rightarrow \overline{\Gamma_0} : \gamma_0  \mapsto \overline{\gamma_0},\]
\[S_k \rightarrow \overline{S_k} : s_k \mapsto \overline{s_k}.\]
Let $\overline{\alpha}, \overline{\beta}, \overline{\delta}, \overline{\sigma}, \overline{\lambda}, \overline{\mu}, \overline{\iota}$ be the functions defined over those new sets, in correspondence with $\alpha, \beta, \delta, \sigma, \lambda, \mu, \iota$, making the diagrams commute.

Let also 
\[\langle \overline{S_0} \st \mathcal{R}_0 \rangle\]
be a finite presentation of the finite group $\Gamma_0$.
Let 
\[\omega : \overline{\Gamma_0} \rightarrow \langle \overline{S_0} \rangle\]
be a function associating an element $\overline{\gamma_0}$ in $\overline{\Gamma_0}$ to a word in the free group generated by $\overline{S_0}$ representing $\gamma_0$ for that presentation.

Define
\[
\overline{A} = \bigcup_{i=0}^n \overline{S_i},
\]
\[
\mathcal{R}_1 = \bigcup_{i = 1}^n \{\overline{s}_0 \overline{s}_i [\overline{\sigma}(\overline{s}_0, \overline{s}_i) \omega(\overline{\alpha}(\overline{s}_0, \overline{s}_i))]^{-1} : \overline{s_0} \in \overline{S}_0, \overline{s_i} \in \overline{S}_i\} \subseteq \langle \overline{A} \rangle, 
\]
\[
\mathcal{R}_2 = \bigcup_{i,j = 1, i > j}^n \{\overline{s}_i \overline{s}_j [\overline{\lambda}(\overline{s}_i, \overline{s}_j) \overline{\mu}(\overline{s}_i, \overline{s}_j) \omega(\overline{\beta}(\overline{s}_i, \overline{s}_j))]^{-1} : \overline{s}_i \in \overline{S}_i, \overline{s}_j \in \overline{S}_j \} \subseteq \langle \overline{A} \rangle,
\]
\[
\mathcal{R}_3 = \bigcup_{i = 1}^n \{\overline{s}_i \overline{\iota}(\overline{s}_i)\omega(\overline{\delta}(\overline{s}_i))^{-1} : \overline{s}_i \in \overline{S}_i\} \subseteq \langle \overline{A} \rangle,
\]
\[
\mathcal{R} := \bigcup_{i=0}^3 \mathcal{R}_i \subseteq \langle \overline{A} \rangle.
\]

\begin{remark}
The sets $\mathcal{R}_1$, $\mathcal{R}_2$ and $\mathcal{R}_3$ encode the relations in $\Gamma$ observed in \Cref{cor:SwapGamma0Si}, \Cref{cor:SwapSiSj} and \Cref{lemma:InverseFuncSi}, respectively.

Informally speaking, in the group $\langle \overline{A} \st \mathcal{R} \rangle$, the relations in $\mathcal{R}_1$ allow to exchange generators in $\overline{S_0}$ with generators in $\overline{S_i}$.,
the relations in $\mathcal{R}_2$ allow to exchange generators in $\overline{S_i}$ with generators in $\overline{S_j}$,
and the relations in $\mathcal{R}_3$ allow some simplifications for products of generators in $\overline{S_i}$.
\end{remark}

\begin{remark}
If $\Gamma$ acts regularly on the product of trees, that is to say if $\Gamma_0$ is the trivial group, then the relations in $\mathcal{R}$ become simpler.
Indeed, one can take $\overline{S_0}$ and $\mathcal{R}_0$ to be empty, the relations in $\mathcal{R}_1$ become tautological, the ones in $\mathcal{R}_2$ become 
\[\overline{s}_i \overline{s}_j (\overline{\lambda}(\overline{s}_i, \overline{s}_j) \overline{\mu}(\overline{s}_i, \overline{s}_j))^{-1},\]
and the ones in $\mathcal{R}_3$ become 
\[\overline{s}_i \overline{\iota}(\overline{s}_i).\]
Moreover, the set $S_i$ is uniquely determined for all $1\leq i \leq n$.
\end{remark}

\begin{theo}\label{thm:PresentationGamma}
A presentation of $\Gamma$ is given by
\[
\langle \overline{A} \st \mathcal{R} \rangle.
\]
\end{theo}

\begin{proof}
By \Cref{cor:SwapGamma0Si}, \Cref{cor:SwapSiSj} and \Cref{lemma:InverseFuncSi}, it is clear that the morphism
\[\langle \overline{A} \rangle \rightarrow \Gamma : \overline{s}_i \in \overline{A} \mapsto s_i \]
induces a morphism 
\[\pi : \langle \overline{A} \st \mathcal{R} \rangle \rightarrow \Gamma.\]
By \Cref{thm:Macbeath}, the set $S$ generates $\Gamma$, and it follows that $\pi$ is surjective by \Cref{prop:SDecomposedAsProd}.
From the same theorem, a presentation of $\Gamma$ is given as $\langle \overline{S} \st \mathcal{Q} \rangle$, where \[\mathcal{Q} = \{\overline{s} \overline{t} (\overline{st})^{-1} : s,t \in S \st st\in S \} \subseteq \langle \overline{S} \rangle,\]
the isomorphism being induced by the morphism $\langle \overline{S}\rangle \rightarrow \Gamma : \overline{s} \in \overline{S} \mapsto s$.
After composing, one gets a surjective morphism
\[\phi : \langle \overline{A} \st \mathcal{R} \rangle \rightarrow \langle \overline{S} \st \mathcal{Q} \rangle.\]
By replacing each symbol $\overline{s}$ in $\overline{S}$ occurring in the relations $\mathcal{Q}$ by a word $w(\overline{s})$ in the free group $\langle \overline{A} \rangle$ such that $\phi(w(\overline{s})) \equiv \overline{s}$, one gets a new presentation $\langle \overline{A} \st \mathcal{Q}' \rangle$ of $\Gamma$, and a surjective morphism
\[\psi : \langle \overline{A} \st \mathcal{R} \rangle \rightarrow \langle \overline{A} \st \mathcal{Q}' \rangle.\]
A good choice for $w(\overline{s})$ is given by \Cref{prop:SDecomposedAsProd}.
Namely, if $s = s_1\cdots s_n \gamma_0$ with $s_k\in S_k\cup \{e\}$ and $\gamma_0\in \Gamma_0$, one defines $w(\overline{s}) = \overline{s_1}\cdots \overline{s_n} \omega(\overline{\gamma_0}) \in \langle \overline{A} \rangle$.

Our goal is to prove that the morphism
\[\psi' : \langle \overline{A} \rangle \rightarrow \langle \overline{A} \st \mathcal{R} \rangle\]
induces a morphism $\langle \overline{A} \st \mathcal{Q}' \rangle \rightarrow \langle \overline{A} \st \mathcal{R} \rangle$, which will be necessarily the inverse of $\psi$.

Let $s,t\in S$ with $s t \in S$.
We need to show that $\psi'(w(s)w(t)) = \psi'(w(st))$.
Say $w(\overline{s}) = \overline{s_1}\cdots \overline{s_n} \omega(\overline{g_0})$ and $w(\overline{t}) = \overline{t_1}\cdots \overline{t_n} \omega(\overline{h_0})$ and $w(\overline{st}) = \overline{a_1}\cdots \overline{a_n} \omega(\overline{i_0})$ with $\overline{s_k},\overline{t_k},\overline{a_k}\in \overline{S_k}\cup \{e\}$ and $\overline{g_0},\overline{h_0},\overline{i_0}\in \overline{\Gamma_0}$.

Using relations in $\mathcal{R}_0$, $\mathcal{R}_1$ and $\mathcal{R}_2$, one sees that $\psi'(w(s)w(t))$ can be written in the form $(\overline{u_1} \overline{v_1}) \cdots (\overline{u_n} \overline{v_n}) \omega(\overline{j_0})$ with $\overline{u_k},\overline{v_k}\in \overline{S_k}\cup \{e\}$ and $\overline{j_0} \in \overline{\Gamma_0}$.
After applying $\pi$ to that element, one gets $st = (u_1 v_1) \cdots (u_n v_n) j_0$ in $\Gamma$.

One has $d(x, st\cdot x) = \max_{1\leq k\leq n} d(x_k, u_k v_k \cdot x_k)$, and this number is $\leq 1$ by the hypothesis that $st \in S$ and by \Cref{lemma:SetSElemDistLe1}.
By \Cref{lemma:InverseFuncSi}, one then has $u_k = e$, or $v_k = e$, or ($u_k\neq e \neq v_k$ and $v_k = \iota(u_k)$), for all $1\leq k \leq n$.
In any case, by the relations in $\mathcal{R}_0$, $\mathcal{R}_1$ and $\mathcal{R}_3$, one can rewrite $\psi'(w(s)w(t))$ in the form $\overline{a'_1}\cdots \overline{a'_n} w(\overline{i'_0})$ with $\overline{a'_k}\in \overline{S_k}\cup \{e\}$ and $\overline{i'_0} \in \overline{\Gamma_0}$.

Applying $\pi$ again to that element, one gets $st = a'_1\cdots a'_n i'_0$.
As an element of $S$, $st$ is uniquely written in that form by \Cref{prop:SDecomposedAsProd}.
In conclusion, it follows that  $a_1 = a'_1,\cdots,a_n = a'_n$ and $i_0 = i'_0$, so that $\psi'(w(s)w(t)) = \psi'(w(st))$, as wanted
\end{proof}

We now turn back to our case of interest.
Let $\Gamma \leqslant \Gamma_O(S)$ be a subgroup acting vertex-transitively on the product of trees $\TT = \prod_{v\in S\backslash S_0} \mathcal{T}_v$, and assume that $S_0 = \Ram(B)$.
In order to find a presentation of $\Gamma$, we need to identify the stabiliser $\Gamma_0$ of a vertex $x$, as well as the $g^{i,j}$'s, which are elements in $\Gamma$ sending $x$ to its neighbours.

The most natural choice for $x$ is the family of completions of $O$, that is $(O_v)_{v\in S\backslash S_0}$.
By \Cref{lemma:VertexStabGammaS}, the stabiliser of $x$ is the finite group $\Gamma_O(S_0) \cap \Gamma$, which can be easily identified.
In Magma, $\Gamma_O(S_0)$ is computed using the method $\mathtt{Normalizer(O)}$.

Now recall that for any $\mathfrak{p}\in S\backslash S_0$, the set $N(O,\mathfrak{p})$ of integral right $O$-ideals of reduced norm $\mathfrak{p}$ is in bijection with the set of $\mathfrak{p}$-neighbours of $O$.
One computes $N(O,\mathfrak{p})$ for those prime ideals.
This task can be done in Magma, for example, with the method $\mathtt{MaximalRightIdeals(O, p)}$.

Then, one exhibits a set $X$ of two-sided $O$-ideals that are representatives of the elements of $\Idl(O)/\PIdl(O)$.
The method $\mathtt{TwoSidedIdealClasses(O})$ is used in Magma for that.
If one has $\Idl(\localiz{O}{S_0}) = \PIdl(\localiz{O}{S_0})$, then one can take $X$ to be a set of products of the prime two-sided $O$-ideals $P_\mathfrak{p}$, for $\mathfrak{p}\in S_0\cap V^K_f$.
Recall that $P_\mathfrak{p}$ is the unique prime two-sided $O$-ideal such that $P_\mathfrak{p}\cap R = \mathfrak{p}$.
In this case, one can specify the optional argument $\mathtt{Support}$ to be $S_0\cap V^K_f$ in the Magma method, that is to say, evaluate $\mathtt{TwoSidedIdealClasses(O : Support := S_0\cap V^K_f})$.

Now for all $I\in N(O,\mathfrak{p})$, there exists a unique $J\in X$ such that $IJ$ is a principal right $O$-ideal, generated by some $a_I\in O$.
It is clear that $a_I\in \Gamma_O(S)$ since $a_I\cdot O = a_I O a_I^{-1}$ is the left order of $I$, which is a $\mathfrak{p}$-neighbour of $O$.
This generating element is computed through the method $\mathtt{IsPrincipal(I*J)}$ in Magma.
Moreover, if $\gamma \in B^\times$ is a representative of an element $\overline{\gamma}\in \Gamma_O(S_0)$, then $a_I \gamma$ is also a generator of $I J$, and there exists at least one $\overline{\gamma}$ such that $a_I \overline{\gamma}\in \Gamma$.
Up to replacing each $a_I$ by $\overline{\gamma} a_I$, assume $\overline{\gamma}$ to be trivial.

The elements $a_I$'s, for $I\in N(O,\mathfrak{p})$, for $\mathfrak{p}\in S\backslash S_0$, are the sought $g^{i,j}$'s, and one can therefore compute a presentation of $\Gamma$ thanks to \Cref{thm:PresentationGamma}.

\section{Examples over number fields}\label{sec:ExNumFields}

In this section and the next one, we give examples of families of quaternionic arithmetic lattices acting regularly on products of trees, i.e.\ BMW groups, together with some presentations.
We recover all the examples presented in \Cref{subsec:ReviewOfKnownResults} and we also introduce new families.

First, we state the following results, which will be useful for proving that some groups we construct are torsion-free.

\begin{lemma}\label{lemma:FiniteOrderInBMWEq2}
Let $\Gamma\leqslant \prod_{k=1}^n \Aut(\mathcal{T}_k)$ be a BMW group.
If $\gamma\in \Gamma$ has finite order, then $\gamma^2 = 1$.
\end{lemma}

\begin{proof}
The group $\Gamma$ acts on each tree individually through the projections.
If $\gamma \in \Gamma$ has finite order, then $\gamma$ stabilises a vertex or an edge in each tree.
In any case, $\gamma^2$ stabilises a vertex in each tree, so it stabilises a vertex in the product.
Since the action of $\Gamma$ is free on the vertices of the product, the element $\gamma^2$ is trivial.
\end{proof}

\begin{lemma}\label{lemma:Order2Trace0}
Let $x\in B^\times$, and assume that $\charac(K)\neq 2$.
Then $x$ has order $2$ in $\PGL_1(B)$ if and only if its reduced trace is zero.
\end{lemma}

\begin{proof}
If $x \in K^\times$, then $x$ has order $1$ in $\PGL_1(B)$ and $\Trd(x) = 2 x \neq 0$, so the equivalence is satisfied.
Assume now that $x\notin K^\times$, so that $x$ has order $2$ in $\PGL_1(B)$ if and only if $x^2 \in K^\times$.
One has
\[x^2 - \Trd(x) x + \Nrd(x) = 0 \text{, i.e.\ } x^2 =  \Trd(x) x - \Nrd(x).\]
It is now clear that $x^2 \in K^\times$ if and only if $\Trd(x) = 0$, as wanted.
\end{proof}

\subsection{\texorpdfstring{$(-1,-1 \st \QQ)$}{(-1,-1 | Q)}}\label{subsec:-1_-1_Q}

Take $K = \QQ$, $R = \ZZ$, $B = (-1,-1 \st \QQ)$.
This is probably the most studied quaternion algebra for constructing BMW groups in the literature.
For instance, it appeared in \cite[Section~7.4]{Lubotzky94}, \cite{Mozes95}, \cite{Rattaggi04} and \cite{RungStixVdov19}.
Those papers studied the first and second families described just below, but we identify more BMW families.

The class group of $\QQ$ is trivial, one has $\Ram(B) = \{\infty, 2\}$,  $B$ has type number~$1$ so that $V^K(O) = V^K\backslash \Ram(B)$, and the maximal Hurwitz order $O$ has class number~$1$:
\[O = \ZZ \oplus \ZZ i \oplus \ZZ j \oplus \ZZ \frac{1+i+j+k}{2}.\]
We will write $[a,b,c,d]$, with $a,b,c,d\in \ZZ$, to represent the element 
\[a+ b i + c j +d \frac{1+i+j+k}{2}\in O.\]

For all finite set of places $S$ containing the place at infinity, the group $\Gamma_O(S)$ acts vertex-transitively on the product of Bruhat-Tits trees $\prod_{v\in S\backslash S_0} \mathcal{T}_{v}$ with stabilisers conjugate to $\Gamma_O(S_0)$ by \Cref{lemma:VertexStabGammaS} and \Cref{cor:GammaSTransIfClass1}, where $S_0 = S\cap \Ram(B)$.
Note that each finite place $v$ corresponds to a prime number $p$, so that $\mathcal{T}_{v}$ is regular of degree $p+1$.

\Cref{prop:ExactSeqGammaScupRam} is also applicable, and one gets the exact sequence
\[1 \rightarrow \Gamma_O(V^K_\infty) \rightarrow \Gamma_O(\Ram(B)) \rightarrow \ZZ/2\ZZ \rightarrow 1.\]

It is well known and not hard to see that $\Gamma_O(V^K_\infty) = \PGL_1(O)$ is isomorphic to the alternating group $\Alt(4)$, which is of order $12$.
A set of representatives in $O$ of the elements in $\Gamma_O(V^K_\infty)$ is
\[U = \{[[1,0,0,0],[0,1,0,0],[0,0,1,0],[-1,-1,-1,2]\} \cup  \{[\epsilon_1,\epsilon_2,\epsilon_3,1] : \epsilon_i \in \{0,-1\}\}\]
Written in another way,
\[U = \{1,i,j,k,\frac{\pm 1\pm i\pm j + k}{2}\}.\] 
Then, $\Gamma_O(\Ram(B))$ is isomorphic to $\Sym(4)$, and a set of representatives is given by
\[U \cup (1+i) U.\]

\paragraph{First family.}

Take any finite set $S$ such that $S_0 = S \cap \Ram(B) = V^K_\infty$.
The modulo-2 reduction group $\PGL_1(O/2O)$ has order $12$.
The reduction morphism $\pi_2 : \Gamma_O(S) \rightarrow \PGL_1(O/2O)$ restricted to $\Gamma_O(S_0)$ is injective, and is consequently an isomorphism since the orders of the groups coincide.
Therefore, by \Cref{lemma:Htransfromquotient}, the kernel $\pi_2^{-1}(1)$ is transitive on the product of trees with trivial stabilisers.
In other words, it acts regularly on the product of trees.

\begin{remark}
This construction is equivalent to the one described by Rungtanapirom Stix and Vdovina in \cite[Proposition 3.11]{RungStixVdov19}.
They noticed the fact that when $S = \{\infty, p,q\}$, where $p$ and $q$ are distinct odd prime numbers, it gives the lattices constructed by Mozes in \cite[Section 3]{Mozes95} if $p\equiv 1 \equiv q \mod 4$.
One also recovers Rattaggi's groups $\Gamma_{p,l,e_0}$ in \cite[Theorem 3.30, Theorem 3.35]{Rattaggi04} if $p\equiv \pm 1$ and  $q\equiv \pm 1 \mod 8$.
\end{remark}

We now give samples of group generators for some places in the following table.
Prime ideals $\mathfrak{p}$ are written in the first column of the table.
In the second one, representatives in $O$ of the generators of the lattice sending $O$ to each of its $\mathfrak{p}$-neighbours are given.

Afterwards, we give presentations of some lattices acting on products of two trees.
We write $\Gamma_{\mathfrak{p},\mathfrak{q}}$, where $\mathfrak{p}$ and $\mathfrak{q}$ are distinct prime ideals, for the lattice acting on the product of the Bruhat-Tits corresponding to those prime ideals.

The same explanations are valid for the other families below.
Beware that two lattices from different families might have the same name without being equal.
Compare $\Gamma_{3,5}$ for this family and $\Gamma_{3,5}$ for the next family.

\vspace{5mm}

\begin{longtable}{r|l}
$3$ & \begin{tabular}{l}
$i-j-k,\hspace{1mm} -i+j-k,\hspace{1mm} i+j+k,\hspace{1mm} i+j-k$
\end{tabular}\\
\hline
$5$ & \begin{tabular}{l}
$1-2k,\hspace{1mm} -1-2j,\hspace{1mm} 1+2i,\hspace{1mm} 1-2i,\hspace{1mm} -1+2j,\hspace{1mm} 1+2k$
\end{tabular}\\
\hline
$7$ & \begin{tabular}{l}
$2+i+j+k,\hspace{1mm} 2-i-j+k,\hspace{1mm} -2-i+j-k,\hspace{1mm} -2+i-j-k,\hspace{1mm} -2-i+j+k,$\\
$-2+i-j+k,\hspace{1mm} 2+i+j-k,\hspace{1mm} 2-i-j-k$
\end{tabular}\\
\hline
$11$ & \begin{tabular}{l}
$-3i+j+k,\hspace{1mm} i+3j-k,\hspace{1mm} -3i+j-k,\hspace{1mm} i+3j+k,\hspace{1mm} -i+3j+k,\hspace{1mm} i+j-3k,$\\
$3i+j+k,\hspace{1mm} -i+3j-k,\hspace{1mm} i-j-3k,\hspace{1mm} i+j+3k,\hspace{1mm} 3i+j-k,\hspace{1mm} -i+j-3k$
\end{tabular}\\
\hline
$13$ & \begin{tabular}{l}
$-1+2i+2j-2k,\hspace{1mm} -3+2k,\hspace{1mm} -3-2j,\hspace{1mm} -1+2i-2j-2k,\hspace{1mm} -3+2i,$\\
$1-2i+2j-2k,\hspace{1mm} 1+2i+2j+2k,\hspace{1mm} 1-2i-2j-2k,\hspace{1mm} 1+2i-2j+2k,$\\
$-3-2i,\hspace{1mm} -1-2i+2j+2k,\hspace{1mm} -3+2j,\hspace{1mm} -3-2k,\hspace{1mm} -1-2i-2j+2k$
\end{tabular}\\
\hline
$17$ & \begin{tabular}{l}
$1-4k,\hspace{1mm} 1-4i,\hspace{1mm} -3+2j-2k,\hspace{1mm} 3+2i-2k,\hspace{1mm} -3+2j+2k,\hspace{1mm} 3+2i-2j,$\\
$-3+2i+2k,\hspace{1mm} -3+2i+2j,\hspace{1mm} 1+4j,\hspace{1mm} 1-4j,\hspace{1mm} -3-2i-2j,\hspace{1mm} -3-2i-2k,$\\
$3-2i+2j,\hspace{1mm} -3-2j-2k,\hspace{1mm} 3-2i+2k,\hspace{1mm} -3-2j+2k,\hspace{1mm} 1+4i,\hspace{1mm} 1+4k$
\end{tabular}\\
\end{longtable}

\begin{align*}
\Gamma_{3,5} =& \langle a_{1}, a_{2}, a_{3}, a_{4},\hspace{1mm} b_{1}, b_{2}, b_{3}\st a_{1}^2,\hspace{1mm} a_{2}^2,\hspace{1mm} a_{3}^2,\hspace{1mm} a_{4}^2,\hspace{1mm} a_{1}b_{1}a_{3}b_{2}^{-1}, \hspace{1mm} a_{1}b_{3}a_{4}b_{1}, \hspace{1mm} a_{1}b_{3}^{-1}a_{2}b_{2}, \hspace{1mm} a_{2}b_{1}a_{4}b_{2}, \hspace{1mm} \\*
 &\hspace{5mm} a_{2}b_{3}^{-1}a_{3}b_{1}, \hspace{1mm} a_{3}b_{2}a_{4}b_{3}\rangle
\\
\Gamma_{3,7} =& \langle a_{1}, a_{2}, a_{3}, a_{4},\hspace{1mm} b_{1}, b_{2}, b_{3}, b_{4}\st a_{1}^2,\hspace{1mm} a_{2}^2,\hspace{1mm} a_{3}^2,\hspace{1mm} a_{4}^2,\hspace{1mm} a_{1}b_{1}a_{2}b_{2}, \hspace{1mm} a_{1}b_{2}a_{3}b_{3}^{-1}, \hspace{1mm} a_{1}b_{4}a_{1}b_{4}^{-1}, \hspace{1mm} \\*
 &\hspace{5mm} a_{1}b_{3}^{-1}a_{4}b_{1}, \hspace{1mm} a_{2}b_{1}a_{4}b_{4}^{-1}, \hspace{1mm} a_{2}b_{3}a_{2}b_{3}^{-1}, \hspace{1mm} a_{2}b_{4}^{-1}a_{3}b_{2}, \hspace{1mm} a_{3}b_{1}a_{3}b_{1}^{-1}, \hspace{1mm} a_{3}b_{3}a_{4}b_{4}, \hspace{1mm} a_{4}b_{2}a_{4}b_{2}^{-1}\rangle
\\
\Gamma_{3,11} =& \langle a_{1}, a_{2}, a_{3}, a_{4},\hspace{1mm} b_{1}, b_{2}, b_{3}, b_{4}, b_{5}, b_{6}, b_{7}, b_{8}, b_{9}, b_{10}, b_{11}, b_{12}\st a_{1}^2,\hspace{1mm} a_{2}^2,\hspace{1mm} a_{3}^2,\hspace{1mm} a_{4}^2,\hspace{1mm} b_{1}^2,\hspace{1mm} b_{2}^2,\hspace{1mm} \\*
 &\hspace{5mm} b_{3}^2,\hspace{1mm} b_{4}^2,\hspace{1mm} b_{5}^2,\hspace{1mm} b_{6}^2,\hspace{1mm} b_{7}^2,\hspace{1mm} b_{8}^2,\hspace{1mm} b_{9}^2,\hspace{1mm} b_{10}^2,\hspace{1mm} b_{11}^2,\hspace{1mm} b_{12}^2,\hspace{1mm} a_{1}b_{1}a_{3}b_{7}, \hspace{1mm} a_{1}b_{2}a_{4}b_{5}, \hspace{1mm} a_{1}b_{3}a_{1}b_{10}, \hspace{1mm} \\*
 &\hspace{5mm} a_{1}b_{4}a_{1}b_{11}, \hspace{1mm} a_{1}b_{6}a_{1}b_{8}, \hspace{1mm} a_{1}b_{9}a_{2}b_{12}, \hspace{1mm} a_{2}b_{1}a_{2}b_{6}, \hspace{1mm} a_{2}b_{2}a_{2}b_{7}, \hspace{1mm} a_{2}b_{3}a_{4}b_{11}, \hspace{1mm} a_{2}b_{4}a_{3}b_{8}, \hspace{1mm} \\*
 &\hspace{5mm} a_{2}b_{5}a_{2}b_{10}, \hspace{1mm} a_{3}b_{2}a_{3}b_{12}, \hspace{1mm} a_{3}b_{3}a_{3}b_{5}, \hspace{1mm} a_{3}b_{6}a_{4}b_{10}, \hspace{1mm} a_{3}b_{9}a_{3}b_{11}, \hspace{1mm} a_{4}b_{1}a_{4}b_{8}, \hspace{1mm} a_{4}b_{4}a_{4}b_{9}, \hspace{1mm} \\*
 &\hspace{5mm} a_{4}b_{7}a_{4}b_{12}\rangle
\\
\Gamma_{5,7} =& \langle a_{1}, a_{2}, a_{3},\hspace{1mm} b_{1}, b_{2}, b_{3}, b_{4}\st a_{1}b_{1}a_{3}^{-1}b_{3}, \hspace{1mm} a_{1}b_{2}a_{3}b_{4}, \hspace{1mm} a_{1}b_{3}a_{2}b_{2}, \hspace{1mm} a_{1}b_{4}a_{2}^{-1}b_{1}, \hspace{1mm} a_{1}b_{4}^{-1}a_{1}b_{4}^{-1}, \hspace{1mm} \\*
 &\hspace{5mm} a_{1}b_{3}^{-1}a_{1}b_{3}^{-1}, \hspace{1mm} a_{1}b_{2}^{-1}a_{1}b_{2}^{-1}, \hspace{1mm} a_{1}b_{1}^{-1}a_{1}b_{1}^{-1}, \hspace{1mm} a_{2}b_{1}a_{2}b_{1}, \hspace{1mm} a_{2}b_{3}a_{3}^{-1}b_{4}^{-1}, \hspace{1mm} a_{2}b_{4}a_{2}b_{4}, \hspace{1mm} \\*
 &\hspace{5mm} a_{2}b_{3}^{-1}a_{2}b_{3}^{-1}, \hspace{1mm} a_{2}b_{2}^{-1}a_{2}b_{2}^{-1}, \hspace{1mm} a_{2}b_{1}^{-1}a_{3}b_{2}, \hspace{1mm} a_{3}b_{1}a_{3}b_{1}, \hspace{1mm} a_{3}b_{3}a_{3}b_{3}, \hspace{1mm} a_{3}b_{4}^{-1}a_{3}b_{4}^{-1}, \hspace{1mm} \\*
 &\hspace{5mm} a_{3}b_{2}^{-1}a_{3}b_{2}^{-1}\rangle
\end{align*}

\paragraph{Second family.}
Take again any $S$ such that $S_0 = S \cap \Ram(B) = V^K_\infty$.
Let us go one step further by looking at the congruence quotient $\PGL_1(O/4O)$.
The natural morphism $\PGL_1(O/4O) \rightarrow \PGL_1(O/2O)$ is surjective with a kernel of order $8$ isomorphic to $C_2^3$; it follows that $\PGL_1(O/4O)$ has order $96$.
This group has a non-normal subgroup $H$ of order $8$ isomorphic to $C_2 \times C_4$ generated by $[1,0,2,0]$ and $[0,1,0,2]$.
The elements of $H$ are 
\[ [1,0,0,0],[0,1,0,2],[1,2,0,0],[2,1,0,2],[1,0,2,0],[2,1,2,2],[1,2,2,0],[0,1,2,2].\]
One also observes that $H$ has a trivial intersection with the (injective) image of $\Gamma_O(S_0) = \PGL_1(O)$ in $\PGL_1(O/4O)$.
It immediately follows from \Cref{lemma:Htransfromquotient} that $\pi_4^{-1}(H)$ acts regularly on the product of trees, where $\pi_4 : \Gamma_O(S) \rightarrow \PGL_1(O/4O)$ is the reduction map.

\begin{proposition}
The group $\pi_4^{-1}(H)$ is torsion-free.
\end{proposition}

\begin{proof}
Let $\Gamma$ be $\pi_4^{-1}(H)$.
By \Cref{lemma:FiniteOrderInBMWEq2}, it suffices to show there is no element of order $2$ in $\Gamma$, that is to say, every element in $\Gamma$ has a representative in $B^\times$ with a non-zero reduced trace by \Cref{lemma:Order2Trace0}.
Let $\gamma = [a,b,c,d]\in \Gamma$.
One has
\[ \gamma = (a+d/2) + (b+d/2) i + (c+d/2)j + d/2 k. \]
Note that the fourth coordinate of each element of $H$ is even, so $d$ is itself even, say $d = 2d'$.
Moreover, $a$ and $d'$ have different parity, so that $a+d/2 = a+d'$ is an odd integer.
It follows that $\gamma$ has a non-zero reduced trace and is of infinite order.
\end{proof}

\begin{remark}
This is equivalent to the construction of Rungtanapirom Stix and Vdovina in \cite[Theorem 3.15]{RungStixVdov19}.
If $S = \{\infty, p,l\}$, where $p$ and $l$ are distinct odd prime numbers with $p$ or $l$ congruent to $3$ modulo $4$, the corresponding lattice is $\Gamma_{p,l}$ constructed by Rattaggi in \cite[Theorem 3.30, Theorem 3.35]{Rattaggi04}.
\end{remark}

\vspace{5mm}

\begin{longtable}{r|l}
$3$ & \begin{tabular}{l}
$-1+j-k,\hspace{1mm} -1-j-k,\hspace{1mm} -1+j+k,\hspace{1mm} -1-j+k$
\end{tabular}\\
\hline
$5$ & \begin{tabular}{l}
$-1-2j,\hspace{1mm} -1-2k,\hspace{1mm} 1+2i,\hspace{1mm} 1-2i,\hspace{1mm} -1+2k,\hspace{1mm} -1+2j$
\end{tabular}\\
\hline
$7$ & \begin{tabular}{l}
$-1+2i+j-k,\hspace{1mm} 1-2i-j-k,\hspace{1mm} -1+2i-j-k,\hspace{1mm} -1-2i+j-k,\hspace{1mm} -1+2i-j+k,$\\
$-1-2i+j+k,\hspace{1mm} 1+2i+j+k,\hspace{1mm} -1-2i-j+k$
\end{tabular}\\
\hline
$11$ & \begin{tabular}{l}
$-3-j+k,\hspace{1mm} -1-3j+k,\hspace{1mm} 1+j-3k,\hspace{1mm} 3-j-k,\hspace{1mm} 1-3j-k,\hspace{1mm} 1-j-3k,$\\
$1+j+3k,\hspace{1mm} 1+3j+k,\hspace{1mm} 3+j+k,\hspace{1mm} 1-j+3k,\hspace{1mm} -1+3j-k,\hspace{1mm} -3+j-k$
\end{tabular}\\
\hline
$13$ & \begin{tabular}{l}
$-1+2i-2j-2k,\hspace{1mm} 1+2i+2j+2k,\hspace{1mm} -3+2j,\hspace{1mm} -3+2k,\hspace{1mm} 1-2i+2j-2k,$\\
$-3+2i,\hspace{1mm} 1+2i+2j-2k,\hspace{1mm} 1-2i-2j+2k,\hspace{1mm} -3-2i,\hspace{1mm} 1+2i-2j+2k,$\\
$-3-2k,\hspace{1mm} -3-2j,\hspace{1mm} 1-2i-2j-2k,\hspace{1mm} -1-2i+2j+2k$
\end{tabular}\\
\hline
$17$ & \begin{tabular}{l}
$3+2i-2j,\hspace{1mm} 3+2i+2j,\hspace{1mm} -3-2j+2k,\hspace{1mm} 3+2j+2k,\hspace{1mm} 1-4k,\hspace{1mm} 1-4i,$\\
$3+2i-2k,\hspace{1mm} -1+4j,\hspace{1mm} 3+2i+2k,\hspace{1mm} 3-2i-2k,\hspace{1mm} -1-4j,\hspace{1mm} 3-2i+2k,$\\
$1+4i,\hspace{1mm} 1+4k,\hspace{1mm} 3-2j-2k,\hspace{1mm} -3+2j-2k,\hspace{1mm} 3-2i-2j,\hspace{1mm} 3-2i+2j$
\end{tabular}\\
\end{longtable}

\begin{align*}
\Gamma_{3,5} =& \langle a_{1}, a_{2},\hspace{1mm} b_{1}, b_{2}, b_{3}\st a_{1}b_{1}a_{1}b_{2}^{-1}, \hspace{1mm} a_{1}b_{2}a_{2}b_{3}^{-1}, \hspace{1mm} a_{1}b_{3}a_{2}b_{2}, \hspace{1mm} a_{1}b_{3}^{-1}a_{2}^{-1}b_{1}^{-1}, \hspace{1mm} a_{1}b_{1}^{-1}a_{2}^{-1}b_{3}, \hspace{1mm} \\*
 &\hspace{5mm} a_{2}b_{2}^{-1}a_{2}b_{1}^{-1}\rangle
\\
\Gamma_{3,7} =& \langle a_{1}, a_{2},\hspace{1mm} b_{1}, b_{2}, b_{3}, b_{4}\st a_{1}b_{1}a_{2}b_{2}^{-1}, \hspace{1mm} a_{1}b_{2}a_{1}^{-1}b_{3}, \hspace{1mm} a_{1}b_{3}a_{2}b_{4}, \hspace{1mm} a_{1}b_{4}a_{2}^{-1}b_{2}, \hspace{1mm} a_{1}b_{4}^{-1}a_{1}b_{1}^{-1}, \hspace{1mm} \\*
 &\hspace{5mm} a_{1}b_{3}^{-1}a_{2}^{-1}b_{1}, \hspace{1mm} a_{2}b_{1}a_{2}^{-1}b_{4}^{-1}, \hspace{1mm} a_{2}b_{2}a_{2}b_{3}^{-1}\rangle 
\\
\Gamma_{3,11} =& \langle a_{1}, a_{2},\hspace{1mm} b_{1}, b_{2}, b_{3}, b_{4}, b_{5}, b_{6}\st a_{1}b_{1}a_{1}^{-1}b_{1}^{-1}, \hspace{1mm} a_{1}b_{2}a_{2}b_{5}, \hspace{1mm} a_{1}b_{3}a_{2}^{-1}b_{6}^{-1}, \hspace{1mm} a_{1}b_{4}a_{1}b_{2}^{-1}, \hspace{1mm} a_{1}b_{5}a_{2}b_{2}, \hspace{1mm} \\*
 &\hspace{5mm} a_{1}b_{6}a_{1}b_{5}^{-1}, \hspace{1mm} a_{1}b_{6}^{-1}a_{2}^{-1}b_{3}, \hspace{1mm} a_{1}b_{4}^{-1}a_{1}b_{3}^{-1}, \hspace{1mm} a_{2}b_{1}a_{2}b_{5}^{-1}, \hspace{1mm} a_{2}b_{3}a_{2}b_{2}^{-1}, \hspace{1mm} a_{2}b_{4}a_{2}^{-1}b_{4}^{-1}, \hspace{1mm} \\*
 &\hspace{5mm} a_{2}b_{6}^{-1}a_{2}b_{1}^{-1}\rangle
\\
\Gamma_{5,7} =& \langle a_{1}, a_{2}, a_{3},\hspace{1mm} b_{1}, b_{2}, b_{3}, b_{4}\st a_{1}b_{1}a_{1}^{-1}b_{4}^{-1}, \hspace{1mm} a_{1}b_{2}a_{3}^{-1}b_{1}, \hspace{1mm} a_{1}b_{3}a_{1}^{-1}b_{2}, \hspace{1mm} a_{1}b_{4}a_{3}b_{3}^{-1}, \hspace{1mm} a_{1}b_{4}^{-1}a_{2}^{-1}b_{1}^{-1}, \hspace{1mm} \\*
 &\hspace{5mm} a_{1}b_{2}^{-1}a_{2}b_{3}, \hspace{1mm} a_{2}b_{1}a_{2}^{-1}b_{4}^{-1}, \hspace{1mm} a_{2}b_{2}a_{2}^{-1}b_{3}, \hspace{1mm} a_{2}b_{4}^{-1}a_{3}^{-1}b_{2}, \hspace{1mm} a_{2}b_{3}^{-1}a_{3}b_{1}^{-1}, \hspace{1mm} a_{3}b_{1}a_{3}b_{4}^{-1}, \hspace{1mm} \\*
 &\hspace{5mm} a_{3}b_{2}a_{3}b_{3}\rangle
\end{align*}

\paragraph{Third family.}
To our knowledge, there is no reference in the literature that mentions this third family nor the following ones.
Take this time any $S$ with $\Ram(B) \subseteq S$ and $3\notin S$.
Since $3\notin \Ram(B)$, the algebra $O/3O$ is isomorphic to $M_2(\FF_3)$, and so $\PGL_1(O/3O)$ is isomorphic to $\PGL_2(\FF_3)$, which has order $24$.
The group $\Gamma_O(S_0)$ also has order $24$, and the reduction map 
$\pi_3:\Gamma_O(S) \rightarrow \PGL_1(O/3O)$ restricted to $\Gamma_O(S_0)$ is injective, and consequently bijective.
Again, by \Cref{lemma:Htransfromquotient}, the kernel $\pi_3^{-1}(1)$ acts regularly on the product of trees.

\begin{proposition}\label{prop:HamiltonOverQGamma3TorsionFree}
The group $\pi_3^{-1}(1)$ is torsion-free.
\end{proposition}

\begin{proof}
Write $\Gamma$ for $\pi_3^{-1}(1)$.
For all $n > 1$, one has the reduction morphism
\[ \Gamma_O(S) \rightarrow \PGL_1(O/(3^n)O), \]
and the image of $\Gamma$ is contained in 
\[ \PGL_1((1 + 3O)/(3^n)O),\]
which is a $3$-group.
It follows that $\Gamma$ is a residually $3$-group, so there cannot be an element of order $2$ in $\Gamma$, and \Cref{lemma:FiniteOrderInBMWEq2} applies.
\end{proof}

\vspace{5mm}

\begin{longtable}{r|l}
$5$ & \begin{tabular}{l}
$1-3j,\hspace{1mm} -1-3k,\hspace{1mm} 1-3i,\hspace{1mm} 1+3i,\hspace{1mm} -1+3k,\hspace{1mm} 1+3j$
\end{tabular}\\
\hline
$7$ & \begin{tabular}{l}
$1/2+3/2i-3/2j+3/2k,\hspace{1mm} -1/2+3/2i-3/2j-3/2k,$\\
$1/2+3/2i+3/2j-3/2k,\hspace{1mm} 1/2+3/2i+3/2j+3/2k,$\\
$1/2-3/2i-3/2j-3/2k,\hspace{1mm} 1/2-3/2i-3/2j+3/2k,$\\
$-1/2-3/2i+3/2j+3/2k,\hspace{1mm} 1/2-3/2i+3/2j-3/2k$
\end{tabular}\\
\hline
$11$ & \begin{tabular}{l}
$2+3j+3k,\hspace{1mm} -2+3i-3j,\hspace{1mm} -2+3i-3k,\hspace{1mm} 2+3i+3j,\hspace{1mm} -2-3i-3k,$\\
$2+3j-3k,\hspace{1mm} 2-3j+3k,\hspace{1mm} -2+3i+3k,\hspace{1mm} 2-3i-3j,\hspace{1mm} -2-3i+3k,$\\
$-2-3i+3j,\hspace{1mm} 2-3j-3k$
\end{tabular}\\
\hline
$13$ & \begin{tabular}{l}
$-2+3i,\hspace{1mm} -5/2+3/2i-3/2j-3/2k,\hspace{1mm} -5/2+3/2i-3/2j+3/2k,$\\
$-2+3k,\hspace{1mm} -5/2-3/2i-3/2j+3/2k,\hspace{1mm} -2+3j,\hspace{1mm} -5/2+3/2i+3/2j+3/2k,$\\
$-5/2-3/2i-3/2j-3/2k,\hspace{1mm} -2-3j,\hspace{1mm} -5/2+3/2i+3/2j-3/2k,$\\
$-2-3k,\hspace{1mm} -5/2-3/2i+3/2j-3/2k,\hspace{1mm} -5/2-3/2i+3/2j+3/2k,$\\
$-2-3i$
\end{tabular}\\
\hline
$17$ & \begin{tabular}{l}
$-4+3i-3k,\hspace{1mm} 4-3i-3j,\hspace{1mm} -4+3j-3k,\hspace{1mm} 5+3i,\hspace{1mm} -5-3k,$\\
$4+3i-3j,\hspace{1mm} -4+3i+3k,\hspace{1mm} 4+3j+3k,\hspace{1mm} -5+3j,\hspace{1mm} -5-3j,$\\
$4-3j-3k,\hspace{1mm} -4-3i-3k,\hspace{1mm} 4-3i+3j,\hspace{1mm} -5+3k,\hspace{1mm} 5-3i,\hspace{1mm} -4-3j+3k,$\\
$4+3i+3j,\hspace{1mm} -4-3i+3k$
\end{tabular}\\
\end{longtable}

\begin{align*}
\Gamma_{5,7} =& \langle a_{1}, a_{2}, a_{3},\hspace{1mm} b_{1}, b_{2}, b_{3}, b_{4}\st a_{1}b_{1}a_{1}b_{3}^{-1}, \hspace{1mm} a_{1}b_{2}a_{2}^{-1}b_{1}^{-1}, \hspace{1mm} a_{1}b_{3}a_{3}b_{2}, \hspace{1mm} a_{1}b_{4}a_{1}b_{2}^{-1}, \hspace{1mm} a_{1}b_{4}^{-1}a_{2}b_{3}, \hspace{1mm} \\*
 &\hspace{5mm} a_{1}b_{1}^{-1}a_{3}^{-1}b_{4}^{-1}, \hspace{1mm} a_{2}b_{1}a_{3}b_{3}, \hspace{1mm} a_{2}b_{2}a_{2}b_{3}^{-1}, \hspace{1mm} a_{2}b_{4}a_{2}b_{1}^{-1}, \hspace{1mm} a_{2}b_{4}^{-1}a_{3}^{-1}b_{2}^{-1}, \hspace{1mm} a_{3}b_{4}^{-1}a_{3}b_{3}^{-1}, \hspace{1mm} \\*
 &\hspace{5mm} a_{3}b_{2}^{-1}a_{3}b_{1}^{-1}\rangle
\\
\Gamma_{5,11} =& \langle a_{1}, a_{2}, a_{3},\hspace{1mm} b_{1}, b_{2}, b_{3}, b_{4}, b_{5}, b_{6}\st a_{1}b_{1}a_{3}^{-1}b_{3}^{-1}, \hspace{1mm} a_{1}b_{2}a_{2}b_{5}^{-1}, \hspace{1mm} a_{1}b_{3}a_{2}^{-1}b_{2}, \hspace{1mm} a_{1}b_{4}a_{2}b_{4}, \hspace{1mm} \\*
 &\hspace{5mm} a_{1}b_{5}a_{3}^{-1}b_{4}^{-1}, \hspace{1mm} a_{1}b_{6}a_{3}b_{6}, \hspace{1mm} a_{1}b_{6}^{-1}a_{2}^{-1}b_{5}, \hspace{1mm} a_{1}b_{5}^{-1}a_{3}b_{1}, \hspace{1mm} a_{1}b_{4}^{-1}a_{3}b_{3}, \hspace{1mm} a_{1}b_{3}^{-1}a_{2}b_{6}^{-1}, \hspace{1mm} \\*
 &\hspace{5mm} a_{1}b_{2}^{-1}a_{3}^{-1}b_{2}^{-1}, \hspace{1mm} a_{1}b_{1}^{-1}a_{2}^{-1}b_{1}^{-1}, \hspace{1mm} a_{2}b_{2}a_{3}^{-1}b_{1}^{-1}, \hspace{1mm} a_{2}b_{3}a_{3}^{-1}b_{3}, \hspace{1mm} a_{2}b_{5}a_{3}b_{5}, \hspace{1mm} a_{2}b_{4}^{-1}a_{3}^{-1}b_{6}, \hspace{1mm} \\*
 &\hspace{5mm} a_{2}b_{2}^{-1}a_{3}b_{4}^{-1}, \hspace{1mm} a_{2}b_{1}^{-1}a_{3}b_{6}^{-1}\rangle
\\
\Gamma_{7,11} =& \langle a_{1}, a_{2}, a_{3}, a_{4},\hspace{1mm} b_{1}, b_{2}, b_{3}, b_{4}, b_{5}, b_{6}\st a_{1}b_{1}a_{3}b_{6}, \hspace{1mm} a_{1}b_{2}a_{2}^{-1}b_{3}, \hspace{1mm} a_{1}b_{3}a_{2}b_{5}^{-1}, \hspace{1mm} a_{1}b_{4}a_{3}b_{2}, \hspace{1mm} \\*
 &\hspace{5mm} a_{1}b_{5}a_{4}b_{4}^{-1}, \hspace{1mm} a_{1}b_{6}a_{3}^{-1}b_{2}^{-1}, \hspace{1mm} a_{1}b_{6}^{-1}a_{2}^{-1}b_{5}, \hspace{1mm} a_{1}b_{5}^{-1}a_{4}^{-1}b_{3}^{-1}, \hspace{1mm} a_{1}b_{4}^{-1}a_{4}^{-1}b_{1}^{-1}, \hspace{1mm} a_{1}b_{3}^{-1}a_{4}b_{1}, \hspace{1mm} \\*
 &\hspace{5mm} a_{1}b_{2}^{-1}a_{2}b_{6}^{-1}, \hspace{1mm} a_{1}b_{1}^{-1}a_{3}^{-1}b_{4}, \hspace{1mm} a_{2}b_{1}a_{4}^{-1}b_{4}^{-1}, \hspace{1mm} a_{2}b_{2}a_{4}^{-1}b_{6}^{-1}, \hspace{1mm} a_{2}b_{3}a_{3}^{-1}b_{5}, \hspace{1mm} a_{2}b_{4}a_{3}^{-1}b_{1}, \hspace{1mm} \\*
 &\hspace{5mm} a_{2}b_{5}a_{3}b_{1}^{-1}, \hspace{1mm} a_{2}b_{4}^{-1}a_{4}b_{6}, \hspace{1mm} a_{2}b_{3}^{-1}a_{3}b_{4}, \hspace{1mm} a_{2}b_{1}^{-1}a_{4}b_{2}, \hspace{1mm} a_{3}b_{3}a_{4}^{-1}b_{6}, \hspace{1mm} a_{3}b_{5}a_{4}^{-1}b_{2}^{-1}, \hspace{1mm} \\*
 &\hspace{5mm} a_{3}b_{5}^{-1}a_{4}b_{3}, \hspace{1mm} a_{3}b_{2}^{-1}a_{4}b_{6}^{-1}\rangle
\end{align*}

\paragraph{Fourth family.}

Now take $S$ with $\Ram(B) \subseteq S$ and $5\notin S$.
Since $5\notin \Ram(B)$, the algebra $O/5O$ is isomorphic to $M_2(\FF_5)$, and so $\PGL_1(O/5O)$ is isomorphic to $ \PGL_2(\FF_5)$, which has order $120$.
Let $\pi_5:\Gamma_O(S) \rightarrow \PGL_1(O/5O)$ be the reduction morphism.
It is surjective by \Cref{rem:strgapproxthmimagecontainsPSL} because $\Nrd(1+i) = 2$ is not a square modulo $5$.
When $\pi_5$ is restricted to $\Gamma_O(S_0)$, it is injective with an image of index $5$.
Take $A$ any of the six $5$-Sylow subgroups of $\PGL_1(O/5O)$.
Then $\pi_5^{-1}(A)$ acts regularly on the product of trees.
Actually, since the $5$-Sylow subgroups are conjugate in $\PGL_1(O/5O)$, the corresponding lattices are also conjugate in $\Gamma_O(S)$ and their respective actions on the product of trees are equivalent.

\begin{proposition}
The group $\pi_5^{-1}(A)$ is torsion-free.
\end{proposition}

\begin{proof}
The idea is similar to the proof of \Cref{prop:HamiltonOverQGamma3TorsionFree}.
One sees that $\pi_5^{-1}(A)$ is a residually $5$-group and is torsion-free by \Cref{lemma:FiniteOrderInBMWEq2}.
\end{proof}

\vspace{5mm}

\begin{longtable}{r|l}
$3$ & \begin{tabular}{l}
$-1-2i+j,\hspace{1mm} 1+i+2j,\hspace{1mm} 1-i-2j,\hspace{1mm} -1+2i-j$
\end{tabular}\\
\hline
$7$ & \begin{tabular}{l}
$-3+i+2j,\hspace{1mm} -3+2i-j,\hspace{1mm} -2+i-3j,\hspace{1mm} -2+3i+j,\hspace{1mm} -2-3i-j,$\\
$-2-i+3j,\hspace{1mm} -3-2i+j,\hspace{1mm} -3-i-2j$
\end{tabular}\\
\hline
$11$ & \begin{tabular}{l}
$-1+3i+j,\hspace{1mm} 3/2+3/2i+1/2j-5/2k,\hspace{1mm} -3/2-3/2i-1/2j-5/2k,$\\
$1+i-3j,\hspace{1mm} 3/2-1/2i+3/2j-5/2k,\hspace{1mm} -3/2-1/2i+3/2j+5/2k,$\\
$-3/2+1/2i-3/2j-5/2k,\hspace{1mm} 3/2+1/2i-3/2j+5/2k,\hspace{1mm} 1-i+3j,$\\
$-3/2+3/2i+1/2j+5/2k,\hspace{1mm} 3/2-3/2i-1/2j+5/2k,\hspace{1mm} -1-3i-j$
\end{tabular}\\
\hline
$13$ & \begin{tabular}{l}
$1+5j,\hspace{1mm} 1+5i,\hspace{1mm} -1-4i-3j,\hspace{1mm} 4-i+3j,\hspace{1mm} -1+5k,\hspace{1mm} -4-3i-j,$\\
$-1-3i+4j,\hspace{1mm} -1+3i-4j,\hspace{1mm} -4+3i+j,\hspace{1mm} -1-5k,\hspace{1mm} 4+i-3j,$\\
$-1+4i+3j,\hspace{1mm} 1-5i,\hspace{1mm} 1-5j$
\end{tabular}\\
\hline
$17$ & \begin{tabular}{l}
$3-5i,\hspace{1mm} 3+4i+3j,\hspace{1mm} -3+5j,\hspace{1mm} -2-i-2j-5k,\hspace{1mm} 3-3i+4j,$\\
$-2-2i+j+5k,\hspace{1mm} -2-i-2j+5k,\hspace{1mm} -2-2i+j-5k,\hspace{1mm} 3+5k,$\\
$3-5k,\hspace{1mm} -2+2i-j+5k,\hspace{1mm} -2+i+2j-5k,\hspace{1mm} -2+2i-j-5k,$\\
$3+3i-4j,\hspace{1mm} -2+i+2j+5k,\hspace{1mm} -3-5j,\hspace{1mm} 3-4i-3j,\hspace{1mm} 3+5i$
\end{tabular}
\end{longtable}

\begin{align*}
\Gamma_{3,7} =& \langle a_{1}, a_{2},\hspace{1mm} b_{1}, b_{2}, b_{3}, b_{4}\st a_{1}b_{1}a_{1}b_{3}^{-1}, \hspace{1mm} a_{1}b_{2}a_{1}^{-1}b_{2}^{-1}, \hspace{1mm} a_{1}b_{3}a_{2}^{-1}b_{3}, \hspace{1mm} a_{1}b_{4}a_{2}b_{4}, \hspace{1mm} a_{1}b_{4}^{-1}a_{1}b_{1}^{-1}, \hspace{1mm} \\*
 &\hspace{5mm} a_{2}b_{1}a_{2}^{-1}b_{1}^{-1}, \hspace{1mm} a_{2}b_{2}a_{2}b_{3}, \hspace{1mm} a_{2}b_{4}^{-1}a_{2}b_{2}^{-1}\rangle
\\
\Gamma_{3,11} =& \langle a_{1}, a_{2},\hspace{1mm} b_{1}, b_{2}, b_{3}, b_{4}, b_{5}, b_{6}\st a_{1}b_{1}a_{1}b_{4}^{-1}, \hspace{1mm} a_{1}b_{2}a_{1}^{-1}b_{3}^{-1}, \hspace{1mm} a_{1}b_{3}a_{2}b_{2}, \hspace{1mm} a_{1}b_{4}a_{2}^{-1}b_{5}, \hspace{1mm} a_{1}b_{5}a_{1}^{-1}b_{6}, \hspace{1mm} \\*
 &\hspace{5mm} a_{1}b_{6}a_{2}^{-1}b_{5}^{-1}, \hspace{1mm} a_{1}b_{6}^{-1}a_{2}^{-1}b_{4}, \hspace{1mm} a_{1}b_{3}^{-1}a_{2}b_{1}^{-1}, \hspace{1mm} a_{1}b_{1}^{-1}a_{2}b_{2}^{-1}, \hspace{1mm} a_{2}b_{1}a_{2}b_{4}, \hspace{1mm} a_{2}b_{3}a_{2}^{-1}b_{2}^{-1}, \hspace{1mm} \\*
 &\hspace{5mm} a_{2}b_{5}a_{2}^{-1}b_{6}\rangle
\\
\Gamma_{7,11} =& \langle a_{1}, a_{2}, a_{3}, a_{4},\hspace{1mm} b_{1}, b_{2}, b_{3}, b_{4}, b_{5}, b_{6}\st a_{1}b_{1}a_{2}b_{3}^{-1}, \hspace{1mm} a_{1}b_{2}a_{1}b_{6}^{-1}, \hspace{1mm} a_{1}b_{3}a_{1}b_{5}, \hspace{1mm} a_{1}b_{4}a_{2}^{-1}b_{5}^{-1}, \hspace{1mm} \\*
 &\hspace{5mm} a_{1}b_{6}a_{2}^{-1}b_{4}, \hspace{1mm} a_{1}b_{5}^{-1}a_{3}b_{1}^{-1}, \hspace{1mm} a_{1}b_{4}^{-1}a_{4}^{-1}b_{2}^{-1}, \hspace{1mm} a_{1}b_{3}^{-1}a_{4}^{-1}b_{4}^{-1}, \hspace{1mm} a_{1}b_{2}^{-1}a_{2}b_{1}, \hspace{1mm} a_{1}b_{1}^{-1}a_{3}b_{6}, \hspace{1mm} \\*
 &\hspace{5mm} a_{2}b_{2}a_{2}b_{5}^{-1}, \hspace{1mm} a_{2}b_{3}a_{2}b_{6}, \hspace{1mm} a_{2}b_{4}a_{4}^{-1}b_{3}^{-1}, \hspace{1mm} a_{2}b_{5}a_{3}^{-1}b_{1}^{-1}, \hspace{1mm} a_{2}b_{2}^{-1}a_{4}^{-1}b_{4}, \hspace{1mm} a_{2}b_{1}^{-1}a_{3}^{-1}b_{6}^{-1}, \hspace{1mm} \\*
 &\hspace{5mm} a_{3}b_{2}a_{4}^{-1}b_{6}^{-1}, \hspace{1mm} a_{3}b_{3}a_{3}^{-1}b_{2}^{-1}, \hspace{1mm} a_{3}b_{4}a_{3}^{-1}b_{4}^{-1}, \hspace{1mm} a_{3}b_{5}a_{4}^{-1}b_{3}, \hspace{1mm} a_{3}b_{6}^{-1}a_{4}b_{3}^{-1}, \hspace{1mm} a_{3}b_{2}^{-1}a_{4}b_{5}, \hspace{1mm} \\*
 &\hspace{5mm} a_{4}b_{1}a_{4}^{-1}b_{1}^{-1}, \hspace{1mm} a_{4}b_{6}a_{4}^{-1}b_{5}\rangle
\end{align*}

\paragraph{More families.}
The group $\PGL_2(\FF_p)$ has subgroups of index $12 = |\Gamma_O(\{\infty\})|$ or $24 = |\Gamma_O(\{\infty,2\})|$ for only finitely many primes $p$ (all less than $23$), so the arguments used for the afore-described families can only possibly apply to these primes.
For $p\in \{7,11,23\}$ one can indeed carry out these arguments and find other BMW families.

\subsection{\texorpdfstring{$(-2,-13 \st \QQ)$}{(-2,-13 | Q)}}\label{subsec:-2_-13_Q}

Take $K = \QQ$, $R = \ZZ$, $B = (-2,-13 \st \QQ)$.
We give this example since it was studied by \citep{Chiu92} and is easier to treat than the one in \Cref{subsec:-1_-1_Q}.
Indeed, the vertex stabilisers for the action of $\Gamma_O(S)$ are small.

One has $\Ram(B) = \{\infty, 13\}$, the quaternion algebra has type number $1$, and $\Idl(O) = \PIdl(O)$ for all maximal order $O$.
For example, take $O$ to be
\[\ZZ \oplus \ZZ i \oplus \ZZ \frac{2-i+k}{4} \oplus \ZZ \frac{1+i+j}{2}.\]
By \Cref{lemma:VertexStabGammaS} and \Cref{cor:GammaSTransIfClass1}, if $S$ is a finite set of places containing $V^K_\infty$, then $\Gamma_O(S)$ acts vertex-transitively on the product of trees with stabilisers conjugate to $\Gamma_O(S_0)$, where $S_0 = S\cap \Ram(B)$.

One checks that the group $\Gamma_O(V^K_\infty)$ is trivial.
It follows that if $S_0 = \{\infty\}$, then the group $\Gamma_O(S)$ already acts regularly on the vertices, giving a family of BMW groups.
We give here below examples of generators and group presentations for this family.

\begin{remark}
This is a generalisation of Chiu's lattice in \cite{Chiu92}, passing from a regular action on the $3$-regular Bruhat-Tits tree at the place $2$, to actions on products of trees.
Chiu's lattice is, however, not a free group as claimed in his paper: it is actually isomorphic to $\ZZ * C_2$.
\end{remark}

One the other hand, $\Gamma_O(\Ram(B))$ has order $2$ from the exact sequence in \Cref{prop:ExactSeqGammaScupRam}, and it is generated by $j$ since $v_{13}(j) = \frac{1}{2}v_{13}(\Nrd(j)) = \frac{1}{2} \notin \ZZ$.
If $S_0 = \Ram(B)$, take a prime $p$ not in $S$, such that $\Nrd(j) = 13$ is not a square modulo $p$, i.e.\ $p$ is odd and congruent to $2,5,6,7,8$ or $11$ modulo $13$ by the law of quadratic reciprocity.
Then $\pi_{p}^{-1}(\PSL_2(\FF_p))$ is also a BMW group, where $f_p : \Gamma_O(S) \rightarrow \PGL_1(O/p O) \simeq \PGL_2(\FF_p)$ is the reduction map, which is surjective by the hypothesis on $p$.

\vspace{5mm}

\begin{longtable}{r|l}
$2$ & \begin{tabular}{l}
$1/2-1/4i+1/4k,\hspace{1mm} i,\hspace{1mm} 1/2+1/4i-1/4k$
\end{tabular}\\
\hline
$3$ & \begin{tabular}{l}
$1/2-3/4i-1/4k,\hspace{1mm} 1-i,\hspace{1mm} 1+i,\hspace{1mm} 1/2+3/4i+1/4k$
\end{tabular}\\
\hline
$5$ & \begin{tabular}{l}
$-3/2+3/4i+1/4k,\hspace{1mm} -1/2-5/4i+1/4k,\hspace{1mm} 1/4i-1/2j+1/4k,$\\
$1/4i+1/2j+1/4k,\hspace{1mm} -1/2+5/4i-1/4k,\hspace{1mm} -3/2-3/4i-1/4k$
\end{tabular}\\
\hline
$7$ & \begin{tabular}{l}
$-1-3/4i-1/2j+1/4k,\hspace{1mm} -1-3/4i+1/2j+1/4k,\hspace{1mm} -3/2-5/4i+1/4k,$\\
$1/2i+1/2k,\hspace{1mm} 1/2i-1/2k,\hspace{1mm} -3/2+5/4i-1/4k,\hspace{1mm} -1+3/4i-1/2j-1/4k,$\\
$-1+3/4i+1/2j-1/4k$
\end{tabular}\\
\hline
$11$ & \begin{tabular}{l}
$-2+1/2i+1/2k,\hspace{1mm} 3-i,\hspace{1mm} -5/2-5/4i+1/4k,\hspace{1mm} -2+1/2i-1/2k,$\\
$3/2i-1/2k,\hspace{1mm} 3/2i+1/2k,\hspace{1mm} 7/4i-1/2j-1/4k,\hspace{1mm} 7/4i+1/2j-1/4k,$\\
$-2-1/2i+1/2k,\hspace{1mm} -5/2+5/4i-1/4k,\hspace{1mm} 3+i,\hspace{1mm} -2-1/2i-1/2k$
\end{tabular}\\
\hline
$17$ & \begin{tabular}{l}
$-1/2+11/4i+1/4k,\hspace{1mm} 3/2-1/4i-3/4k,\hspace{1mm} -3+5/4i+1/2j+1/4k,$\\
$-2-j,\hspace{1mm} -3+5/4i-1/2j+1/4k,\hspace{1mm} -3/2+1/4i-j-1/4k,$\\
$3-2i,\hspace{1mm} -7/2+5/4i-1/4k,\hspace{1mm} -3/2+1/4i+j-1/4k,\hspace{1mm} -3/2-1/4i-j+1/4k,$\\
$-7/2-5/4i+1/4k,\hspace{1mm} 3+2i,\hspace{1mm} -3/2-1/4i+j+1/4k,\hspace{1mm} -3-5/4i+1/2j-1/4k,$\\
$-2+j,\hspace{1mm} -3-5/4i-1/2j-1/4k,\hspace{1mm} 3/2+1/4i+3/4k,\hspace{1mm} -1/2-11/4i-1/4k$
\end{tabular}
\end{longtable}

\begin{align*}
\Gamma_{2,3} =& \langle a_{1}, a_{2},\hspace{1mm} b_{1}, b_{2}\st a_{2}^2,\hspace{1mm} a_{1}b_{1}a_{2}b_{1}, \hspace{1mm} a_{1}b_{2}a_{1}b_{2}, \hspace{1mm} a_{1}b_{2}^{-1}a_{1}b_{1}^{-1}, \hspace{1mm} a_{2}b_{2}a_{2}b_{2}^{-1}\rangle
\\
\Gamma_{2,5} =& \langle a_{1}, a_{2},\hspace{1mm} b_{1}, b_{2}, b_{3}, b_{4}\st a_{2}^2,\hspace{1mm} b_{3}^2,\hspace{1mm} b_{4}^2,\hspace{1mm} a_{1}b_{1}a_{1}b_{2}^{-1}, \hspace{1mm} a_{1}b_{2}a_{2}b_{2}, \hspace{1mm} a_{1}b_{3}a_{1}^{-1}b_{4}, \hspace{1mm} a_{1}b_{4}a_{2}b_{1}^{-1}, \hspace{1mm} \\*
 &\hspace{5mm} a_{1}b_{1}^{-1}a_{2}b_{3}\rangle
\\
\Gamma_{2,7} =& \langle a_{1}, a_{2},\hspace{1mm} b_{1}, b_{2}, b_{3}, b_{4}, b_{5}\st a_{2}^2,\hspace{1mm} b_{4}^2,\hspace{1mm} b_{5}^2,\hspace{1mm} a_{1}b_{1}a_{2}b_{3}^{-1}, \hspace{1mm} a_{1}b_{2}a_{1}^{-1}b_{1}^{-1}, \hspace{1mm} a_{1}b_{3}a_{1}b_{4}, \hspace{1mm} a_{1}b_{5}a_{1}^{-1}b_{5}, \hspace{1mm} \\*
 &\hspace{5mm} a_{1}b_{3}^{-1}a_{2}b_{2}, \hspace{1mm} a_{1}b_{1}^{-1}a_{2}b_{2}^{-1}, \hspace{1mm} a_{2}b_{4}a_{2}b_{5}\rangle
\\
\Gamma_{2,11} =& \langle a_{1}, a_{2},\hspace{1mm} b_{1}, b_{2}, b_{3}, b_{4}, b_{5}, b_{6}, b_{7}, b_{8}\st a_{2}^2,\hspace{1mm} b_{5}^2,\hspace{1mm} b_{6}^2,\hspace{1mm} b_{7}^2,\hspace{1mm} b_{8}^2,\hspace{1mm} a_{1}b_{1}a_{1}b_{5}, \hspace{1mm} a_{1}b_{2}a_{1}b_{3}^{-1}, \hspace{1mm} \\*
 &\hspace{5mm} a_{1}b_{3}a_{2}b_{8}, \hspace{1mm} a_{1}b_{4}a_{1}^{-1}b_{4}^{-1}, \hspace{1mm} a_{1}b_{6}a_{1}^{-1}b_{7}, \hspace{1mm} a_{1}b_{7}a_{2}b_{3}, \hspace{1mm} a_{1}b_{8}a_{1}^{-1}b_{6}, \hspace{1mm} a_{1}b_{2}^{-1}a_{1}b_{1}^{-1}, \hspace{1mm} \\*
 &\hspace{5mm} a_{2}b_{1}a_{2}b_{4}^{-1}, \hspace{1mm} a_{2}b_{2}a_{2}b_{2}^{-1}, \hspace{1mm} a_{2}b_{5}a_{2}b_{6}\rangle
\\
\Gamma_{3,5} =& \langle a_{1}, a_{2},\hspace{1mm} b_{1}, b_{2}, b_{3}, b_{4}\st b_{3}^2,\hspace{1mm} b_{4}^2,\hspace{1mm} a_{1}b_{1}a_{1}^{-1}b_{1}^{-1}, \hspace{1mm} a_{1}b_{2}a_{2}^{-1}b_{3}, \hspace{1mm} a_{1}b_{3}a_{1}^{-1}b_{4}, \hspace{1mm} a_{1}b_{4}a_{2}^{-1}b_{2}, \hspace{1mm} \\*
 &\hspace{5mm} a_{1}b_{2}^{-1}a_{1}b_{2}^{-1}, \hspace{1mm} a_{2}b_{1}a_{2}b_{1}, \hspace{1mm} a_{2}b_{2}a_{2}b_{1}^{-1}, \hspace{1mm} a_{2}b_{3}a_{2}^{-1}b_{4}\rangle
\\
\Gamma_{3,7} =& \langle a_{1}, a_{2},\hspace{1mm} b_{1}, b_{2}, b_{3}, b_{4}, b_{5}\st b_{4}^2,\hspace{1mm} b_{5}^2,\hspace{1mm} a_{1}b_{1}a_{1}b_{1}, \hspace{1mm} a_{1}b_{2}a_{1}b_{2}, \hspace{1mm} a_{1}b_{3}a_{2}b_{4}, \hspace{1mm} a_{1}b_{4}a_{2}b_{3}, \hspace{1mm} a_{1}b_{5}a_{2}^{-1}b_{1}^{-1}, \hspace{1mm} \\*
 &\hspace{5mm} a_{1}b_{3}^{-1}a_{2}b_{2}^{-1}, \hspace{1mm} a_{1}b_{2}^{-1}a_{2}^{-1}b_{5}, \hspace{1mm} a_{1}b_{1}^{-1}a_{2}b_{3}^{-1}, \hspace{1mm} a_{2}b_{1}a_{2}^{-1}b_{2}^{-1}\rangle
\\
\Gamma_{5,7} =& \langle a_{1}, a_{2}, a_{3}, a_{4},\hspace{1mm} b_{1}, b_{2}, b_{3}, b_{4}, b_{5}\st a_{3}^2,\hspace{1mm} a_{4}^2,\hspace{1mm} b_{4}^2,\hspace{1mm} b_{5}^2,\hspace{1mm} a_{1}b_{1}a_{1}^{-1}b_{2}^{-1}, \hspace{1mm} a_{1}b_{2}a_{3}b_{3}^{-1}, \hspace{1mm} a_{1}b_{3}a_{1}b_{5}, \hspace{1mm} \\*
 &\hspace{5mm} a_{1}b_{4}a_{2}b_{1}^{-1}, \hspace{1mm} a_{1}b_{3}^{-1}a_{4}b_{1}, \hspace{1mm} a_{1}b_{2}^{-1}a_{2}b_{4}, \hspace{1mm} a_{2}b_{1}a_{3}b_{5}, \hspace{1mm} a_{2}b_{2}a_{2}^{-1}b_{1}^{-1}, \hspace{1mm} a_{2}b_{3}a_{2}^{-1}b_{3}^{-1}, \hspace{1mm} \\*
 &\hspace{5mm} a_{2}b_{5}a_{4}b_{2}, \hspace{1mm} a_{3}b_{1}a_{4}b_{2}^{-1}, \hspace{1mm} a_{3}b_{3}a_{4}b_{3}^{-1}, \hspace{1mm} a_{3}b_{4}a_{4}b_{4}\rangle
\end{align*}

\subsection{\texorpdfstring{$(-13,-6 \st \QQ)$}{(-13,-6 | Q)}}\label{subsec:-13_-6_Q}

Take $K = \QQ$, $R = \ZZ$, $B = (-13,-6 \st \QQ)$.
As for \Cref{subsec:-2_-13_Q}, the vertex-stabilisers for the action of $\Gamma_O(S)$ are small.
Moreover, this is an example where the class groups of the maximal orders are not trivial.
One has $\Ram(B) = \{\infty, 2,3,13\}$, the quaternion algebra has type number $1$, and $\Idl(O)/\PIdl(O)$ has order $2$ for all maximal order $O$.
For example, take $O$ to be
\[\ZZ \oplus \ZZ i \oplus \ZZ \frac{1+i+j}{2} \oplus \ZZ \frac{1+i+k}{2}.\]
Let $P_2,P_3,P_{13}$ denote the prime two-sided $O$-ideals whose intersections with $R$ are the ideals $(2),(3),(13)$, respectively.
One sees that $P_{13}$ is principal generated by $i$, whereas $P_2$ and $P_3$ are not principal and those two individually generate $\Idl(O)/\PIdl(O)$, so that $\Idl(\localiz{O}{S_0}) = \PIdl(\localiz{O}{S_0})$ if $2\in S_0$ or $3\in S_0$.
By \Cref{lemma:VertexStabGammaS} and \Cref{cor:GammaSTransIfClass1}, if $2\in S$ or $3\in S$, then $\Gamma_O(S)$ acts vertex-transitively on the product of trees, with stabilisers conjugate to $\Gamma_O(S_0)$ where, as always, $S_0 = S\cap \Ram(B)$.
Moreover, $\Gamma_O(S_0)$ is trivial if $S_0 \cap V^K_f$ is equal to $\{2\}$ or $\{3\}$.
We get in this way two families of BMW groups.

\hspace{5mm}

First family with $S_0 \cap V^K_f = \{2\}$.

\begin{longtable}{r|l}
$5$ & \begin{tabular}{l}
$2+j,\hspace{1mm} 1/2-1/2i-1/2j,\hspace{1mm} 1/2-1/2i+1/2j,\hspace{1mm} 1/2+1/2i-1/2j,$\\
$1/2+1/2i+1/2j,\hspace{1mm} 2-j$
\end{tabular}\\
\hline
$7$ & \begin{tabular}{l}
$3/2-1/2i+1/2j,\hspace{1mm} 1+j,\hspace{1mm} 3/2-1/2i-1/2j,\hspace{1mm} 1-i,\hspace{1mm} 1+i,\hspace{1mm} 3/2+1/2i+1/2j,$\\
$1-j,\hspace{1mm} 3/2+1/2i-1/2j$
\end{tabular}\\
\hline
$11$ & \begin{tabular}{l}
$5/2+1/2i-1/2j,\hspace{1mm} 5/2-1/2i-1/2j,\hspace{1mm} 1-1/2j+1/2k,\hspace{1mm} 3+i,$\\
$4-j,\hspace{1mm} -1+1/2j+1/2k,\hspace{1mm} -1-1/2j-1/2k,\hspace{1mm} 4+j,\hspace{1mm} 3-i,\hspace{1mm} 1+1/2j-1/2k,$\\
$5/2+1/2i+1/2j,\hspace{1mm} 5/2-1/2i+1/2j$
\end{tabular}\\
\hline
$17$ & \begin{tabular}{l}
$1-3/2j-1/2k,\hspace{1mm} 1/2-1/2i-3/2j,\hspace{1mm} 1-3/2j+1/2k,\hspace{1mm} 7/2-1/2i+1/2j,$\\
$2+i,\hspace{1mm} 1/2+1/2i-3/2j,\hspace{1mm} 7/2+1/2i+1/2j,\hspace{1mm} i-1/2j+1/2k,$\\
$i-1/2j-1/2k,\hspace{1mm} i+1/2j+1/2k,\hspace{1mm} i+1/2j-1/2k,\hspace{1mm} 7/2-1/2i-1/2j,$\\
$1/2-1/2i+3/2j,\hspace{1mm} 2-i,\hspace{1mm} 7/2+1/2i-1/2j,\hspace{1mm} 1+3/2j-1/2k,$\\
$1/2+1/2i+3/2j,\hspace{1mm} 1+3/2j+1/2k$
\end{tabular}
\end{longtable}

\begin{align*}
\Gamma_{5,7} =& \langle a_{1}, a_{2}, a_{3},\hspace{1mm} b_{1}, b_{2}, b_{3}, b_{4}\st a_{1}b_{1}a_{1}b_{1}, \hspace{1mm} a_{1}b_{2}a_{1}^{-1}b_{2}^{-1}, \hspace{1mm} a_{1}b_{3}a_{3}b_{4}, \hspace{1mm} a_{1}b_{4}a_{3}b_{3}, \hspace{1mm} a_{1}b_{4}^{-1}a_{2}^{-1}b_{1}^{-1}, \hspace{1mm} \\*
 &\hspace{5mm} a_{1}b_{3}^{-1}a_{1}b_{3}^{-1}, \hspace{1mm} a_{1}b_{1}^{-1}a_{2}^{-1}b_{4}^{-1}, \hspace{1mm} a_{2}b_{2}a_{3}^{-1}b_{2}, \hspace{1mm} a_{2}b_{3}a_{2}^{-1}b_{3}^{-1}, \hspace{1mm} a_{2}b_{4}^{-1}a_{3}b_{4}^{-1}, \hspace{1mm} a_{2}b_{2}^{-1}a_{2}b_{1}^{-1}, \hspace{1mm} \\*
 &\hspace{5mm} a_{3}b_{1}a_{3}^{-1}b_{1}^{-1}, \hspace{1mm} a_{3}b_{2}a_{3}b_{3}^{-1}\rangle
\\
\Gamma_{5,11} =& \langle a_{1}, a_{2}, a_{3},\hspace{1mm} b_{1}, b_{2}, b_{3}, b_{4}, b_{5}, b_{6}\st a_{1}b_{1}a_{2}b_{4}, \hspace{1mm} a_{1}b_{2}a_{3}^{-1}b_{4}^{-1}, \hspace{1mm} a_{1}b_{3}a_{1}b_{6}, \hspace{1mm} a_{1}b_{4}a_{2}b_{1}, \hspace{1mm} a_{1}b_{5}a_{1}^{-1}b_{5}^{-1}, \hspace{1mm} \\*
 &\hspace{5mm} a_{1}b_{6}^{-1}a_{2}b_{1}^{-1}, \hspace{1mm} a_{1}b_{4}^{-1}a_{3}^{-1}b_{2}, \hspace{1mm} a_{1}b_{3}^{-1}a_{3}^{-1}b_{2}^{-1}, \hspace{1mm} a_{1}b_{2}^{-1}a_{3}^{-1}b_{6}^{-1}, \hspace{1mm} a_{1}b_{1}^{-1}a_{2}b_{3}^{-1}, \hspace{1mm} a_{2}b_{2}a_{2}^{-1}b_{2}^{-1}, \hspace{1mm} \\*
 &\hspace{5mm} a_{2}b_{3}a_{3}^{-1}b_{5}, \hspace{1mm} a_{2}b_{5}a_{3}^{-1}b_{6}, \hspace{1mm} a_{2}b_{6}a_{2}^{-1}b_{3}^{-1}, \hspace{1mm} a_{2}b_{5}^{-1}a_{2}b_{4}^{-1}, \hspace{1mm} a_{3}b_{1}a_{3}^{-1}b_{1}^{-1}, \hspace{1mm} a_{3}b_{5}a_{3}b_{4}^{-1}, \hspace{1mm} \\*
 &\hspace{5mm} a_{3}b_{6}a_{3}^{-1}b_{3}^{-1}\rangle
\\
\Gamma_{7,11} =& \langle a_{1}, a_{2}, a_{3}, a_{4},\hspace{1mm} b_{1}, b_{2}, b_{3}, b_{4}, b_{5}, b_{6}\st a_{1}b_{1}a_{1}^{-1}b_{1}^{-1}, \hspace{1mm} a_{1}b_{2}a_{4}^{-1}b_{5}, \hspace{1mm} a_{1}b_{3}a_{2}^{-1}b_{4}^{-1}, \hspace{1mm} a_{1}b_{4}a_{3}b_{4}, \hspace{1mm} \\*
 &\hspace{5mm} a_{1}b_{5}a_{4}^{-1}b_{2}, \hspace{1mm} a_{1}b_{6}a_{3}^{-1}b_{5}^{-1}, \hspace{1mm} a_{1}b_{6}^{-1}a_{1}b_{6}^{-1}, \hspace{1mm} a_{1}b_{5}^{-1}a_{3}^{-1}b_{3}, \hspace{1mm} a_{1}b_{4}^{-1}a_{2}^{-1}b_{6}, \hspace{1mm} a_{1}b_{3}^{-1}a_{1}b_{3}^{-1}, \hspace{1mm} \\*
 &\hspace{5mm} a_{1}b_{2}^{-1}a_{2}^{-1}b_{2}^{-1}, \hspace{1mm} a_{2}b_{1}a_{3}b_{1}, \hspace{1mm} a_{2}b_{3}a_{4}^{-1}b_{2}^{-1}, \hspace{1mm} a_{2}b_{5}a_{2}^{-1}b_{5}^{-1}, \hspace{1mm} a_{2}b_{6}a_{4}b_{1}^{-1}, \hspace{1mm} a_{2}b_{6}^{-1}a_{3}b_{4}^{-1}, \hspace{1mm} \\*
 &\hspace{5mm} a_{2}b_{4}^{-1}a_{3}b_{3}^{-1}, \hspace{1mm} a_{2}b_{2}^{-1}a_{4}^{-1}b_{6}, \hspace{1mm} a_{2}b_{1}^{-1}a_{4}b_{3}, \hspace{1mm} a_{3}b_{2}a_{3}^{-1}b_{2}^{-1}, \hspace{1mm} a_{3}b_{3}a_{3}b_{3}, \hspace{1mm} a_{3}b_{6}a_{3}b_{6}, \hspace{1mm} \\*
 &\hspace{5mm} a_{3}b_{5}^{-1}a_{4}^{-1}b_{1}^{-1}, \hspace{1mm} a_{3}b_{1}^{-1}a_{4}^{-1}b_{5}^{-1}, \hspace{1mm} a_{4}b_{4}a_{4}^{-1}b_{4}^{-1}, \hspace{1mm} a_{4}b_{6}a_{4}^{-1}b_{3}^{-1}\rangle
\end{align*}

\vspace{5mm}

Second family with $S_0 \cap V^K_f = \{3\}$.

\begin{longtable}{r|l}
$5$ & \begin{tabular}{l}
$1/2+1/2i-1/2j,\hspace{1mm} 1/2-1/2i-1/2j,\hspace{1mm} 3-j,\hspace{1mm} 3+j,\hspace{1mm} 1/2+1/2i+1/2j,$\\
$1/2-1/2i+1/2j$
\end{tabular}\\
\hline
$7$ & \begin{tabular}{l}
$3/2+1/2i+1/2j,\hspace{1mm} 3/2+1/2i-1/2j,\hspace{1mm} 1+j,\hspace{1mm} 1/2j+1/2k,\hspace{1mm} 1/2j-1/2k,$\\
$1-j,\hspace{1mm} 3/2-1/2i+1/2j,\hspace{1mm} 3/2-1/2i-1/2j$
\end{tabular}\\
\hline
$11$ & \begin{tabular}{l}
$3/2+3/2i-1/2j,\hspace{1mm} 3+2j,\hspace{1mm} 5/2-1/2i+1/2j,\hspace{1mm} 3/2+3/2i+1/2j,$\\
$5/2-1/2i-1/2j,\hspace{1mm} 3/2j-1/2k,\hspace{1mm} 3/2j+1/2k,\hspace{1mm} 5/2+1/2i+1/2j,$\\
$3/2-3/2i-1/2j,\hspace{1mm} 5/2+1/2i-1/2j,\hspace{1mm} 3-2j,\hspace{1mm} 3/2-3/2i+1/2j$
\end{tabular}\\
\hline
$17$ & \begin{tabular}{l}
$2+i,\hspace{1mm} 3/2-3/2i+1/2k,\hspace{1mm} 9/2+3/2i-1/2j,\hspace{1mm} 9/2-3/2i-1/2j,$\\
$7/2+1/2i+1/2j,\hspace{1mm} 1/2+1/2i-3/2j,\hspace{1mm} 1/2-1/2i-3/2j,\hspace{1mm} 3/2+3/2i+1/2k,$\\
$7/2+1/2i-1/2j,\hspace{1mm} 7/2-1/2i+1/2j,\hspace{1mm} 3/2-3/2i-1/2k,\hspace{1mm} 1/2+1/2i+3/2j,$\\
$1/2-1/2i+3/2j,\hspace{1mm} 7/2-1/2i-1/2j,\hspace{1mm} 9/2+3/2i+1/2j,\hspace{1mm} 9/2-3/2i+1/2j,$\\
$3/2+3/2i-1/2k,\hspace{1mm} 2-i$
\end{tabular}
\end{longtable}

\begin{align*}
\Gamma_{5,7} =& \langle a_{1}, a_{2}, a_{3},\hspace{1mm} b_{1}, b_{2}, b_{3}, b_{4}, b_{5}\st b_{4}^2,\hspace{1mm} b_{5}^2,\hspace{1mm} a_{1}b_{1}a_{3}^{-1}b_{5}, \hspace{1mm} a_{1}b_{2}a_{1}^{-1}b_{2}^{-1}, \hspace{1mm} a_{1}b_{3}a_{2}b_{3}, \hspace{1mm} a_{1}b_{4}a_{3}^{-1}b_{1}, \hspace{1mm} \\*
 &\hspace{5mm} a_{1}b_{5}a_{1}^{-1}b_{4}, \hspace{1mm} a_{1}b_{3}^{-1}a_{1}b_{1}^{-1}, \hspace{1mm} a_{2}b_{1}a_{2}^{-1}b_{1}^{-1}, \hspace{1mm} a_{2}b_{2}a_{2}b_{3}^{-1}, \hspace{1mm} a_{2}b_{4}a_{2}^{-1}b_{5}, \hspace{1mm} a_{2}b_{5}a_{3}^{-1}b_{2}^{-1}, \hspace{1mm} \\*
 &\hspace{5mm} a_{2}b_{2}^{-1}a_{3}^{-1}b_{4}, \hspace{1mm} a_{3}b_{1}a_{3}b_{2}^{-1}, \hspace{1mm} a_{3}b_{3}a_{3}^{-1}b_{3}^{-1}\rangle
\\
\Gamma_{5,11} =& \langle a_{1}, a_{2}, a_{3},\hspace{1mm} b_{1}, b_{2}, b_{3}, b_{4}, b_{5}, b_{6}, b_{7}\st b_{6}^2,\hspace{1mm} b_{7}^2,\hspace{1mm} a_{1}b_{1}a_{3}b_{5}^{-1}, \hspace{1mm} a_{1}b_{2}a_{1}b_{4}, \hspace{1mm} a_{1}b_{3}a_{1}^{-1}b_{3}^{-1}, \hspace{1mm} \\*
 &\hspace{5mm} a_{1}b_{5}a_{3}b_{7}, \hspace{1mm} a_{1}b_{6}a_{3}b_{5}, \hspace{1mm} a_{1}b_{7}a_{2}b_{2}^{-1}, \hspace{1mm} a_{1}b_{5}^{-1}a_{3}b_{1}, \hspace{1mm} a_{1}b_{4}^{-1}a_{2}^{-1}b_{1}^{-1}, \hspace{1mm} a_{1}b_{2}^{-1}a_{2}b_{6}, \hspace{1mm} \\*
 &\hspace{5mm} a_{1}b_{1}^{-1}a_{2}^{-1}b_{4}^{-1}, \hspace{1mm} a_{2}b_{2}a_{2}b_{1}^{-1}, \hspace{1mm} a_{2}b_{3}a_{3}b_{4}^{-1}, \hspace{1mm} a_{2}b_{5}a_{2}^{-1}b_{5}^{-1}, \hspace{1mm} a_{2}b_{7}a_{3}b_{3}^{-1}, \hspace{1mm} a_{2}b_{4}^{-1}a_{3}b_{3}, \hspace{1mm} \\*
 &\hspace{5mm} a_{2}b_{3}^{-1}a_{3}b_{6}, \hspace{1mm} a_{3}b_{2}a_{3}^{-1}b_{2}^{-1}, \hspace{1mm} a_{3}b_{4}a_{3}b_{1}^{-1}\rangle
\\
\Gamma_{7,11} =& \langle a_{1}, a_{2}, a_{3}, a_{4}, a_{5},\hspace{1mm} b_{1}, b_{2}, b_{3}, b_{4}, b_{5}, b_{6}, b_{7}\st a_{4}^2,\hspace{1mm} a_{5}^2,\hspace{1mm} b_{6}^2,\hspace{1mm} b_{7}^2,\hspace{1mm} a_{1}b_{1}a_{2}^{-1}b_{4}^{-1}, \hspace{1mm} a_{1}b_{2}a_{4}b_{3}^{-1}, \hspace{1mm} \\*
 &\hspace{5mm} a_{1}b_{3}a_{3}^{-1}b_{3}, \hspace{1mm} a_{1}b_{4}a_{3}^{-1}b_{6}, \hspace{1mm} a_{1}b_{5}a_{1}^{-1}b_{5}^{-1}, \hspace{1mm} a_{1}b_{6}a_{1}^{-1}b_{7}, \hspace{1mm} a_{1}b_{7}a_{3}^{-1}b_{4}, \hspace{1mm} a_{1}b_{4}^{-1}a_{2}^{-1}b_{1}, \hspace{1mm} \\*
 &\hspace{5mm} a_{1}b_{3}^{-1}a_{5}b_{2}, \hspace{1mm} a_{1}b_{2}^{-1}a_{1}b_{1}^{-1}, \hspace{1mm} a_{2}b_{1}a_{3}b_{6}, \hspace{1mm} a_{2}b_{2}a_{2}b_{4}^{-1}, \hspace{1mm} a_{2}b_{3}a_{2}^{-1}b_{3}^{-1}, \hspace{1mm} a_{2}b_{5}a_{3}b_{5}, \hspace{1mm} \\*
 &\hspace{5mm} a_{2}b_{6}a_{2}^{-1}b_{7}, \hspace{1mm} a_{2}b_{7}a_{3}b_{1}, \hspace{1mm} a_{2}b_{5}^{-1}a_{5}b_{2}^{-1}, \hspace{1mm} a_{2}b_{2}^{-1}a_{4}b_{5}^{-1}, \hspace{1mm} a_{3}b_{2}a_{3}^{-1}b_{2}^{-1}, \hspace{1mm} a_{3}b_{3}a_{5}b_{4}, \hspace{1mm} \\*
 &\hspace{5mm} a_{3}b_{4}a_{4}b_{3}, \hspace{1mm} a_{3}b_{5}^{-1}a_{4}b_{1}^{-1}, \hspace{1mm} a_{3}b_{1}^{-1}a_{5}b_{5}^{-1}, \hspace{1mm} a_{4}b_{1}a_{5}b_{1}^{-1}, \hspace{1mm} a_{4}b_{4}a_{5}b_{4}^{-1}, \hspace{1mm} a_{4}b_{6}a_{5}b_{7}\rangle
\end{align*}

\subsection{\texorpdfstring{$(-7-\sqrt{7},-8+3\sqrt{7} \st \QQ(\sqrt{7}))$}{(-7-sqrt7,-8+3sqrt7 | Q(sqrt7))}}\label{subsec:-7-al_-8+3al_Qsqrt7}

Take $K = \QQ(\sqrt{7})$, $R = \ZZ[\sqrt{7}]$ and $B = (-7-\sqrt{7},-8+3\sqrt{7} \st K)$.
This is our first example where the base field is not $\QQ$.
The class groups of the maximal orders are not trivial, as for \Cref{subsec:-13_-6_Q}.

\begin{remark}\label{rem:FirstExampleFound}
Interestingly, this example and the following one in \Cref{subsec:-1_-1_Qsqrt13} were the starting point of this project.
At the time, we already knew the BMW group $\Gamma_{2-\alpha,-2-\alpha}$ presented below and other ones from \Cref{subsec:-1_-1_Qsqrt13}, without knowing their arithmeticity.
We first came across them by enumerating all the possible BMW groups acting on a product of two $4$-regular trees.
We then discovered non-trivial quaternionic representations, for example $\Gamma_{2-\alpha,-2-\alpha} \rightarrow \PGL_1(B)$, by carefully analysing the $L_2$-quotients provided by the algorithm developed in \cite{Jambor15} by Jambor, which was implemented in Magma \cite{Magma}.
We now understand that these representations correspond to the description that we give here as subgroups of $\Gamma_O(S)\leq \PGL_1(B)$.
\end{remark}

The quaternion algebra $B$ has type number one, and $K$ has a trivial class group.
One has $\Ram(B) = \{\infty_1, \infty_2, \mathfrak{p}_2, \mathfrak{p}_7 \}$, where $\infty_1$ and $\infty_2$ are the two places at infinity and $\mathfrak{p}_2$ and $\mathfrak{p}_7$ are the prime ideals $(3+\sqrt{7})$ and $(\sqrt{7})$, respectively.
Note that $\mathfrak{p}_2^2 = (2)$ and $\mathfrak{p}_7^2 = (7)$.
For any maximal order $O$, the quotient $\Idl(O)/\PIdl(O)$ is not trivial but has order $4$, isomorphic to $C_2\times C_2$.
It is generated by the prime two-sided $O$-ideals of reduced norms $\mathfrak{p}_2$ and $\mathfrak{p}_7$, so that $\Idl(\localiz{O}{\Ram(B)}) = \PIdl(\localiz{O}{\Ram(B)})$.
For example, take $O$ to be the following maximal order:
\[R \oplus R i \oplus R \frac{1}{2}(\sqrt{7}+\sqrt{7}i+j) \oplus R \frac{2+\sqrt{7}}{6} ((1+\sqrt{7})+\sqrt{7}i + k).\]

Let $S$ be a finite set of places containing $\Ram(B)$, and $S_0 = S\cap \Ram(B) = \Ram(B)$.
By \Cref{lemma:VertexStabGammaS} and \Cref{cor:GammaSTransIfClass1}, $\Gamma_O(S)$ is vertex-transitive on the product of the Bruhat-Tits trees, with stabilisers conjugate to $\Gamma_O(S_0)$.
The group $\Gamma_O(S_0)$ has order $2$, generated by $j$.

Now let $\mathfrak{p}$ be a prime ideal not belonging to $S$ and such that $\Nrd(j) = 8-3\sqrt{7}$ is not a square modulo $\mathfrak{p}$.
Examples of such prime ideals are $(2-\sqrt{7})$ or $(3-2\sqrt{7})$.
The $R/\mathfrak{p}$-algebra $O/\mathfrak{p}O$ is isomorphic to $M_2(R/\mathfrak{p})$ since $\mathfrak{p}$ does not ramify in $B$.
Moreover, the reduction map $\pi_\mathfrak{p} : \Gamma_O(S) \rightarrow \PGL_1(O/\mathfrak{p}O)\simeq \PGL_2(R/\mathfrak{p})$ is surjective since $\Nrd(j)$ is not a square modulo $\mathfrak{p}$.
Note that the characteristic of $R/\mathfrak{p}$ is not $2$, so $\PSL_1(O/\mathfrak{p}O) \simeq \PSL_2(R/\mathfrak{p})$ has index $2$ in $\PGL_1(O/\mathfrak{p}O) \simeq \PGL_2(R/\mathfrak{p})$.
Take $A$ to be $\PSL_1(O/\mathfrak{p}O)$, and \Cref{lemma:Htransfromquotient} ensures that the group $\pi_\mathfrak{p}^{-1}(A)$ acts regularly on the vertices of the product of trees.

\vspace{5mm}

Family with $\mathfrak{p} = (3 - 2 \sqrt{7})$.
We write $\alpha$ for $\sqrt{7}$.

\begin{longtable}{r|l}
$2-\alpha$ & \begin{tabular}{l}
$1/2(3\alpha +7)+1/6(\alpha +2)i+1/6(14\alpha +37)k,\hspace{1mm} 1/6(-\alpha +1)i+1/6(\alpha +5)k,$\\
$(3\alpha +8)k,\hspace{1mm} 1/2(3\alpha +7)+1/6(-\alpha -2)i+1/6(-14\alpha -37)k$
\end{tabular}\\
\hline
$-2-\alpha$ & \begin{tabular}{l}
$1/2(-3\alpha -7)+1/6(-\alpha -2)i+1/6(4\alpha +11)k,$\\
$1/3(-\alpha -2)i+1/3(-5\alpha -13)k,\hspace{1mm} 1/6(5\alpha +13)i+1/6(-5\alpha -13)k,$\\
$1/2(-3\alpha -7)+1/6(\alpha +2)i+1/6(-4\alpha -11)k$
\end{tabular}
\end{longtable}

\begin{align*}
\Gamma_{2-\alpha,-2-\alpha} =& \langle a_{1}, a_{2}, a_{3},\hspace{1mm} b_{1}, b_{2}, b_{3}\st a_{2}^2,\hspace{1mm} a_{3}^2,\hspace{1mm} b_{2}^2,\hspace{1mm} b_{3}^2,\hspace{1mm} a_{1}b_{1}a_{1}b_{2}, \hspace{1mm} a_{1}b_{3}a_{1}b_{3}, \hspace{1mm} a_{1}b_{1}^{-1}a_{3}b_{1}^{-1}, \hspace{1mm} a_{2}b_{1}a_{2}b_{1}^{-1}, \hspace{1mm} \\*
 &\hspace{5mm} a_{2}b_{2}a_{3}b_{3}\rangle
\end{align*}

\subsection{\texorpdfstring{$(-1,-1 \st \QQ(\sqrt{13}))$}{(-1,-1 | Q(sqrt13))}}\label{subsec:-1_-1_Qsqrt13}

Take $K = \QQ(\sqrt{13})$, $R = \ZZ[\alpha]$ where $\alpha = (1+\sqrt{13})/2$, and $B = (-1,-1 \st K)$.
This is an easier example over a quadratic field than the one in \Cref{subsec:-7-al_-8+3al_Qsqrt7}.
It also has some importance in this project, as explained in \Cref{rem:FirstExampleFound}.

The field $K$ has a trivial class group, the quaternion algebra has type number $1$, one has $\Ram(B) = V^K_{\infty}$ and $\Idl(O) = \PIdl(O)$ for any maximal order $O$.
For example, take $O$ to be
\[R \oplus R i \oplus R \left(\frac{1-\sqrt{13}}{4}+ \frac{1+\sqrt{13}}{4} i+ \frac{1}{2} j \right) \oplus R \left(\frac{1+\sqrt{13}}{4}+ \frac{1-\sqrt{13}}{4} i+ \frac{1}{2} k \right).\]

Let $S$ be a finite set of places containing $\Ram(B)$.
In all cases, $S_0 = S\cap \Ram(B) = \Ram(B)$.
By \Cref{lemma:VertexStabGammaS} and \Cref{cor:GammaSTransIfClass1}, $\Gamma_O(S)$ is vertex-transitive on the product of the Bruhat-Tits trees, with stabilisers conjugate to $\Gamma_O(S_0)$.
One checks that 
\[\Gamma_O(S_0) = \{1,i,j,k,\frac{1\pm i\pm j \pm k}{2}\};\]
it is a group of order $12$, isomorphic to $\Alt(4)$.

\paragraph{First family.}
The ideal $I_2 = (2)$ is prime; take $S$ not containing it.
The $\localiz{R}{S}/I_2$-algebra $\localiz{O}{S}/I_2$ is isomorphic to $M_2(\FF_4)$, so $\PGL_1(\localiz{O}{S}/I_2) \simeq \PGL_2(\FF_4) \simeq \Alt(5)$, and 

The reduction map $\pi_{(2)}:\Gamma_O(S) \rightarrow \PGL_1(O/I_2) \simeq \PGL_2(\FF_4) \simeq \Alt(5)$ is surjective by the strong approximation theorem.
One checks that $\pi_{(2)}$ restricted to $\Gamma_O(S_0)$ is injective.
It follows by \Cref{lemma:Htransfromquotient} that if $A$ is any $5$-Sylow subgroup of $\PGL_1(O/I_2)$, the group $\pi_{(2)}^{-1}(A)$ acts regularly on the product of trees.

\vspace{5mm}

\begin{longtable}{r|l}
$\alpha$ & \begin{tabular}{l}
$1/2(-\alpha -2)-1/2i+1/2(\alpha +1)j,\hspace{1mm} 1/2(-\alpha -1)+1/2(-\alpha -2)i+1/2j,$\\
$1/2(-\alpha -1)+1/2(\alpha +2)i-1/2j,\hspace{1mm} 1/2(-\alpha -2)+1/2i+1/2(-\alpha -1)j$
\end{tabular}\\
\hline
$1-\alpha$ & \begin{tabular}{l}
$1/2+1/2(-\alpha -1)i-1/2\alpha j,\hspace{1mm} 1/2\alpha -1/2i+1/2(\alpha +1)j,$\\
$1/2\alpha +1/2i+1/2(-\alpha -1)j,\hspace{1mm} 2+1/2(\alpha +1)i+1/2\alpha j$
\end{tabular}\\
\hline
$-1 + 2 \alpha$ & \begin{tabular}{l}
$1/2(-\alpha -3)+1/2(\alpha +2)i-1/2j,\hspace{1mm} 1/2+1/2(-\alpha -3)i+1/2(\alpha +2)j,$\\
$1/2(-2\alpha -1)+1/2(-\alpha -1)i-1/2\alpha j,\hspace{1mm} 1/2(\alpha +1)+1/2\alpha i+1/2(2\alpha +1)j,$\\
$-i+(-\alpha -1)j,\hspace{1mm} i+(-\alpha -1)j,\hspace{1mm} j+(\alpha +1)k,\hspace{1mm} (-\alpha -1)i+k,$\\
$j+(-\alpha -1)k,\hspace{1mm} (-\alpha -1)i-k,\hspace{1mm} 1/2(\alpha +1)-1/2\alpha i+1/2(-2\alpha -1)j,$\\
$1/2(-2\alpha -1)+1/2(\alpha +1)i+1/2\alpha j,\hspace{1mm} 1/2+1/2(\alpha +3)i+1/2(-\alpha -2)j,$\\
$1/2(-\alpha -3)+1/2(-\alpha -2)i+1/2j$
\end{tabular}
\end{longtable}

\begin{align*}
\Gamma_{\alpha,1-\alpha} =& \langle a_{1}, a_{2},\hspace{1mm} b_{1}, b_{2}\st a_{1}b_{1}a_{1}b_{1}, \hspace{1mm} a_{1}b_{2}a_{1}^{-1}b_{2}^{-1}, \hspace{1mm} a_{1}b_{1}^{-1}a_{2}b_{1}^{-1}, \hspace{1mm} a_{2}b_{1}a_{2}b_{2}, \hspace{1mm} a_{2}b_{2}^{-1}a_{2}b_{2}^{-1}\rangle
\\
\Gamma_{\alpha,-1+2\alpha} =& \langle a_{1}, a_{2},\hspace{1mm} b_{1}, b_{2}, b_{3}, b_{4}, b_{5}, b_{6}, b_{7}, b_{8}, b_{9}, b_{10}\st b_{5}^2,\hspace{1mm} b_{6}^2,\hspace{1mm} b_{7}^2,\hspace{1mm} b_{8}^2,\hspace{1mm} b_{9}^2,\hspace{1mm} b_{10}^2,\hspace{1mm} a_{1}b_{1}a_{1}b_{4}, \hspace{1mm} a_{1}b_{2}a_{2}^{-1}b_{2}, \hspace{1mm} \\*
 &\hspace{5mm} a_{1}b_{3}a_{2}b_{10}, \hspace{1mm} a_{1}b_{5}a_{1}b_{3}^{-1}, \hspace{1mm} a_{1}b_{6}a_{1}^{-1}b_{6}, \hspace{1mm} a_{1}b_{7}a_{1}^{-1}b_{8}, \hspace{1mm} a_{1}b_{8}a_{2}b_{3}, \hspace{1mm} a_{1}b_{9}a_{2}^{-1}b_{4}^{-1}, \hspace{1mm} a_{1}b_{10}a_{1}^{-1}b_{9}, \hspace{1mm} \\*
 &\hspace{5mm} a_{1}b_{4}^{-1}a_{2}^{-1}b_{7}, \hspace{1mm} a_{1}b_{2}^{-1}a_{1}b_{1}^{-1}, \hspace{1mm} a_{2}b_{1}a_{2}^{-1}b_{1}^{-1}, \hspace{1mm} a_{2}b_{2}a_{2}b_{3}^{-1}, \hspace{1mm} a_{2}b_{5}a_{2}^{-1}b_{9}, \hspace{1mm} a_{2}b_{6}a_{2}b_{4}^{-1}, \hspace{1mm} \\*
 &\hspace{5mm} a_{2}b_{7}a_{2}^{-1}b_{5}, \hspace{1mm} a_{2}b_{8}a_{2}^{-1}b_{10}\rangle
\\
\Gamma_{1-\alpha,-1+2\alpha} =& \langle a_{1}, a_{2},\hspace{1mm} b_{1}, b_{2}, b_{3}, b_{4}, b_{5}, b_{6}, b_{7}, b_{8}, b_{9}, b_{10}\st b_{5}^2,\hspace{1mm} b_{6}^2,\hspace{1mm} b_{7}^2,\hspace{1mm} b_{8}^2,\hspace{1mm} b_{9}^2,\hspace{1mm} b_{10}^2,\hspace{1mm} a_{1}b_{1}a_{2}b_{7}, \hspace{1mm} a_{1}b_{2}a_{2}^{-1}b_{10}, \hspace{1mm} \\*
 &\hspace{5mm} a_{1}b_{3}a_{1}^{-1}b_{3}^{-1}, \hspace{1mm} a_{1}b_{4}a_{2}^{-1}b_{4}, \hspace{1mm} a_{1}b_{5}a_{1}^{-1}b_{8}, \hspace{1mm} a_{1}b_{6}a_{1}b_{2}^{-1}, \hspace{1mm} a_{1}b_{7}a_{1}^{-1}b_{9}, \hspace{1mm} a_{1}b_{8}a_{2}^{-1}b_{2}, \hspace{1mm} \\*
 &\hspace{5mm} a_{1}b_{9}a_{2}b_{1}, \hspace{1mm} a_{1}b_{10}a_{1}^{-1}b_{5}, \hspace{1mm} a_{1}b_{4}^{-1}a_{1}b_{1}^{-1}, \hspace{1mm} a_{2}b_{2}a_{2}b_{3}^{-1}, \hspace{1mm} a_{2}b_{3}a_{2}b_{4}, \hspace{1mm} a_{2}b_{5}a_{2}b_{1}^{-1}, \hspace{1mm} a_{2}b_{6}a_{2}^{-1}b_{6}, \hspace{1mm} \\*
 &\hspace{5mm} a_{2}b_{9}a_{2}^{-1}b_{8}, \hspace{1mm} a_{2}b_{10}a_{2}^{-1}b_{7}\rangle
\end{align*}

\paragraph{Second family.}
The $R$-ideal $(3)$ decomposes as a product of the two prime ideals $\mathfrak{p}_3 = (\alpha)$ and $\mathfrak{p}'_3 = (1-\alpha)$.
Take $S$ not containing $\mathfrak{p}_3$.
The reduction map $\pi_{\mathfrak{p}_3}:\Gamma_O(S) \rightarrow \PGL_1(O/\mathfrak{p}_3 O)\simeq \PGL_2(\FF_3)$ maps $\Gamma_O(S_0)$ bijectively onto its image, which is necessarily $\PSL_2(\FF_3)$.
If $A$ is any subgroup of $\PGL_2(\FF_3)$ of order $2$ intersecting trivially $\PSL_2(\FF_3)$, then $\pi_{\mathfrak{p}_3}^{-1}(A)$ is a BMW group by \Cref{lemma:Htransfromquotient}.

\vspace{5mm}

\begin{longtable}{r|l}
$2$ & \begin{tabular}{l}
$1/2\alpha +1/2i+1/2(\alpha -1)j,\hspace{1mm} 1/2(\alpha -1)i+1/2j+1/2\alpha k,$\\
$i-j,\hspace{1mm} 1/2(-\alpha +1)i-1/2j+1/2\alpha k,\hspace{1mm} 1/2\alpha -1/2i+1/2(-\alpha +1)j$
\end{tabular}\\
\hline
$1-\alpha$ & \begin{tabular}{l}
$1/2\alpha -1/2i+1/2(\alpha +1)j,\hspace{1mm} 1/2(\alpha +1)i-1/2j-1/2\alpha k,$\\
$1/2(-\alpha -1)i+1/2j-1/2\alpha k,\hspace{1mm} 1/2\alpha +1/2i+1/2(-\alpha -1)j$
\end{tabular}\\
\hline
$2\alpha-1$ & \begin{tabular}{l}
$1/2\alpha +1/2(2\alpha +1)i+1/2(-\alpha -1)j,$\\ 
$-1/2\alpha +1/2(\alpha +2)i+1/2(-\alpha -2)j-1/2\alpha k,$\\
$1/2\alpha +1/2(\alpha +2)i+1/2(-\alpha -2)j+1/2\alpha k,$\\
$\hspace{1mm} 1/2(-\alpha -3)+1/2(-\alpha -2)i-1/2j,$\\
$1/2(-\alpha -1)i+1/2(2\alpha +1)j+1/2\alpha k,\hspace{1mm} 1/2i+1/2(\alpha +2)j+1/2(-\alpha -3)k,$\\
$1/2i+1/2(\alpha +2)j+1/2(\alpha +3)k,\hspace{1mm} (-\alpha -1)i+j,$\\
$1/2(-\alpha -1)i+1/2(2\alpha +1)j-1/2\alpha k,\hspace{1mm} i+(-\alpha -1)j,$\\
$1/2(-\alpha -3)+1/2(\alpha +2)i+1/2j,$\\
$1/2\alpha +1/2(-\alpha -2)i+1/2(\alpha +2)j-1/2\alpha k,$\\
$-1/2\alpha +1/2(-\alpha -2)i+1/2(\alpha +2)j+1/2\alpha k,$\\
$1/2\alpha +1/2(-2\alpha -1)i+1/2(\alpha +1)j$
\end{tabular}
\end{longtable}

\begin{align*}
\Gamma_{2,1-\alpha} =& \langle a_{1}, a_{2}, a_{3}, a_{4},\hspace{1mm} b_{1}, b_{2}, b_{3}\st a_{2}^2,\hspace{1mm} a_{3}^2,\hspace{1mm} a_{4}^2,\hspace{1mm} b_{2}^2,\hspace{1mm} b_{3}^2,\hspace{1mm} a_{1}b_{1}a_{3}b_{1}, \hspace{1mm} a_{1}b_{2}a_{2}b_{1}^{-1}, \hspace{1mm} a_{1}b_{3}a_{1}^{-1}b_{2}, \hspace{1mm} \\*
 &\hspace{5mm} a_{1}b_{1}^{-1}a_{4}b_{3}, \hspace{1mm} a_{2}b_{1}a_{4}b_{1}^{-1}, \hspace{1mm} a_{2}b_{3}a_{3}b_{3}, \hspace{1mm} a_{3}b_{2}a_{4}b_{2}\rangle
\\
\Gamma_{2,-1+2\alpha} =& \langle a_{1}, a_{2}, a_{3}, a_{4},\hspace{1mm} b_{1}, b_{2}, b_{3}, b_{4}, b_{5}, b_{6}, b_{7}, b_{8}, b_{9}, b_{10}\st a_{2}^2,\hspace{1mm} a_{3}^2,\hspace{1mm} a_{4}^2,\hspace{1mm} b_{5}^2,\hspace{1mm} b_{6}^2,\hspace{1mm} b_{7}^2,\hspace{1mm} b_{8}^2,\hspace{1mm} b_{9}^2,\hspace{1mm} \\*
 &\hspace{5mm} b_{10}^2,\hspace{1mm} a_{1}b_{1}a_{1}b_{6}, \hspace{1mm} a_{1}b_{2}a_{3}b_{4}, \hspace{1mm} a_{1}b_{3}a_{2}b_{9}, \hspace{1mm} a_{1}b_{4}a_{3}b_{2}, \hspace{1mm} a_{1}b_{5}a_{1}b_{3}^{-1}, \hspace{1mm} a_{1}b_{7}a_{4}b_{1}^{-1}, \hspace{1mm} \\*
 &\hspace{5mm} a_{1}b_{8}a_{1}^{-1}b_{7}, \hspace{1mm} a_{1}b_{9}a_{1}^{-1}b_{8}, \hspace{1mm} a_{1}b_{10}a_{1}b_{4}^{-1}, \hspace{1mm} a_{1}b_{2}^{-1}a_{2}b_{3}, \hspace{1mm} a_{1}b_{1}^{-1}a_{4}b_{2}^{-1}, \hspace{1mm} a_{2}b_{1}a_{2}b_{4}, \hspace{1mm} \\*
 &\hspace{5mm} a_{2}b_{5}a_{2}b_{8}, \hspace{1mm} a_{2}b_{6}a_{3}b_{7}, \hspace{1mm} a_{2}b_{10}a_{4}b_{2}, \hspace{1mm} a_{3}b_{1}a_{3}b_{3}, \hspace{1mm} a_{3}b_{5}a_{4}b_{9}, \hspace{1mm} a_{3}b_{8}a_{3}b_{10}, \hspace{1mm} a_{4}b_{3}a_{4}b_{4}^{-1}, \hspace{1mm} \\*
 &\hspace{5mm} a_{4}b_{6}a_{4}b_{8}\rangle
\\
\Gamma_{1 -\alpha,-1+2\alpha} =& \langle a_{1}, a_{2}, a_{3},\hspace{1mm} b_{1}, b_{2}, b_{3}, b_{4}, b_{5}, b_{6}, b_{7}, b_{8}, b_{9}, b_{10}\st a_{2}^2,\hspace{1mm} a_{3}^2,\hspace{1mm} b_{5}^2,\hspace{1mm} b_{6}^2,\hspace{1mm} b_{7}^2,\hspace{1mm} b_{8}^2,\hspace{1mm} b_{9}^2,\hspace{1mm} b_{10}^2,\hspace{1mm} \\*
 &\hspace{5mm} a_{1}b_{1}a_{2}b_{1}, \hspace{1mm} a_{1}b_{2}a_{2}b_{3}, \hspace{1mm} a_{1}b_{3}a_{2}b_{7}, \hspace{1mm} a_{1}b_{4}a_{2}b_{6}, \hspace{1mm} a_{1}b_{5}a_{3}b_{4}, \hspace{1mm} a_{1}b_{6}a_{1}^{-1}b_{5}, \hspace{1mm} a_{1}b_{7}a_{1}^{-1}b_{10}, \hspace{1mm} \\*
 &\hspace{5mm} a_{1}b_{8}a_{1}^{-1}b_{8}, \hspace{1mm} a_{1}b_{9}a_{3}b_{1}^{-1}, \hspace{1mm} a_{1}b_{10}a_{1}^{-1}b_{9}, \hspace{1mm} a_{1}b_{4}^{-1}a_{1}b_{2}^{-1}, \hspace{1mm} a_{1}b_{3}^{-1}a_{3}b_{3}^{-1}, \hspace{1mm} a_{1}b_{1}^{-1}a_{3}b_{2}, \hspace{1mm} \\*
 &\hspace{5mm} a_{2}b_{2}a_{3}b_{8}, \hspace{1mm} a_{2}b_{4}a_{3}b_{4}^{-1}, \hspace{1mm} a_{2}b_{5}a_{2}b_{9}, \hspace{1mm} a_{2}b_{10}a_{3}b_{10}, \hspace{1mm} a_{3}b_{6}a_{3}b_{7}\rangle
\end{align*}

\subsection{\texorpdfstring{$(-1,-1 \st \QQ(\sqrt{21}))$}{(-1,-1 | Q(sqrt21))}}\label{subsec:-1_-1_Qsqrt21}

We now give an example where the type number is not $1$.
Take $K = \QQ(\sqrt{21})$, $R = \ZZ[\alpha]$ where $\alpha = (1+\sqrt{21})/2$ and $B = (-1,-1 \st K)$.
The number field $K$ has class number $1$, its narrow class group has order $2$, and $\Ram(B) = V^K_\infty$.
The quaternion algebra $B$ has type number $2$; in other words, there are two conjugacy classes of maximal orders in $B$, and representatives are given by those two:
\[O_1 = R \oplus R i \oplus R \frac{1}{2}((1+\alpha) + \alpha i + j) \oplus R \frac{1}{2}(\alpha + (1+\alpha) i + k),\]
\[O_2 = R \oplus R (1+\alpha) i \oplus R \frac{1}{2} ((1+\alpha) +\alpha i + j) \oplus  R \frac{1+\alpha}{6} (3+i+j+k).\]
Both $V^K(O_1)$ and $V^K(O_2)$ are non-empty, and they both consist of the totally positive principal prime $R$-ideals by \Cref{cor:VKOtypeGt1}.

One has $\Gamma_{O_1}(S_0)\simeq \Alt(4)$ whereas $\Gamma_{O_2}(S_0)\simeq D_6$.
Also, every two-sided $O_1$-ideal is principal, and the same is true for $O_2$.
By \Cref{lemma:VertexStabGammaS} and \Cref{cor:GammaSTransIfClass1}, the groups $\Gamma_{O_1}(S)$ and $\Gamma_{O_2}(S)$ both act vertex-transitively with stabilisers conjugate to $\Gamma_{O_1}(S_0)$ and $\Gamma_{O_2}(S_0)$, respectively, upon the condition that $S \cap V^K_f \subseteq V^K(O_1) = V^K(O_2)$.

\paragraph{First family.} 
Representatives of $\Gamma_{O_1}(S_0)$ are given by the set
\[\{1,i,j,k,\frac{\pm 1\pm i\pm j + k}{2}\}.\]
Write $\mathfrak{p}_3$ for the $R$-ideal generated by $1+\alpha$.
Note that $\mathfrak{p}_3^2 = (3)$, and $\mathfrak{p}_3$ is not totally positive so it does not belong to $S\cap V^K_f \subseteq V^K(O)$ by \Cref{cor:VKOtypeGt1}.
One has $\PGL_1(O_1/\mathfrak{p}_3 O_1)\simeq \PGL_2(3)$.
Moreover, the reduction map $\pi_{\mathfrak{p}_3} : \Gamma_{O_1}(S) \rightarrow \PGL_1(O_1/\mathfrak{p}_3 O_1)$ is injective when restricted to $\Gamma_{O_1}(S_0)$ and its image is necessarily $\PSL_1(O_1/\mathfrak{p}_3 O_1)$.
If $A$ is any subgroup of $\PGL_1(O_1/\mathfrak{p}_3 O_1)$ of order $2$ intersecting trivially $\PSL_1(O_1/\mathfrak{p}_3 O_1)$, then $f_{\mathfrak{p}_3}^{-1}(A)$ is a BMW group.

\vspace{5mm}

\begin{longtable}{r|l}
$2$ & \begin{tabular}{l}
$1/2(\alpha +1)+1/2(-\alpha -2)i+1/2j,\hspace{1mm} -i+j,$\\
$-1/2i+1/2(\alpha +2)j+1/2(-\alpha -1)k,$\\
$-1/2i+1/2(\alpha +2)j+1/2(\alpha +1)k,\hspace{1mm} 1/2(\alpha +1)+1/2(\alpha +2)i-1/2j$
\end{tabular}\\
\hline
$3+\alpha$ & \begin{tabular}{l}
$1/2(-\alpha -1)+1/2\alpha i+1/2j,\hspace{1mm} 1/2(-\alpha -1)+1/2(\alpha +2)i+1/2(-2\alpha -3)j,$\\
$-1/2i-1/2\alpha j+1/2(\alpha +1)k,\hspace{1mm} 1/2(2\alpha +3)i+1/2(-\alpha -2)j+1/2(\alpha +1)k,$\\
$1/2(2\alpha +3)i+1/2(-\alpha -2)j+1/2(-\alpha -1)k,$\\
$1/2i+1/2\alpha j+1/2(\alpha +1)k,\hspace{1mm} 1/2(-\alpha -1)+1/2(-\alpha -2)i+1/2(2\alpha +3)j,$\\
$1/2(-\alpha -1)-1/2\alpha i-1/2j$
\end{tabular}
\end{longtable}

\begin{align*}
\Gamma_{2,3+\alpha} =& \langle a_{1}, a_{2}, a_{3}, a_{4},\hspace{1mm} b_{1}, b_{2}, b_{3}, b_{4}, b_{5}, b_{6}\st a_{2}^2,\hspace{1mm} a_{3}^2,\hspace{1mm} a_{4}^2,\hspace{1mm} b_{3}^2,\hspace{1mm} b_{4}^2,\hspace{1mm} b_{5}^2,\hspace{1mm} b_{6}^2,\hspace{1mm} a_{1}b_{1}a_{3}b_{6}, \hspace{1mm} a_{1}b_{2}a_{2}b_{1}^{-1}, \hspace{1mm} \\*
 &\hspace{5mm} a_{1}b_{3}a_{4}b_{1}, \hspace{1mm} a_{1}b_{4}a_{3}b_{2}^{-1}, \hspace{1mm} a_{1}b_{5}a_{1}^{-1}b_{4}, \hspace{1mm} a_{1}b_{6}a_{1}^{-1}b_{3}, \hspace{1mm} a_{1}b_{2}^{-1}a_{4}b_{5}, \hspace{1mm} a_{1}b_{1}^{-1}a_{2}b_{2}, \hspace{1mm} a_{2}b_{3}a_{3}b_{5}, \hspace{1mm} \\*
 &\hspace{5mm} a_{2}b_{4}a_{4}b_{6}, \hspace{1mm} a_{3}b_{1}a_{4}b_{1}^{-1}, \hspace{1mm} a_{3}b_{2}a_{4}b_{2}^{-1}\rangle
\end{align*}

\paragraph{Second family.} 
The case of $O_2$ is trickier.
Indeed, the reduction map $\pi_{\mathfrak{p}_3} : \Gamma_{O_2}(S) \rightarrow \PGL_1(O_2/\mathfrak{p}_3 O_2) \simeq \PGL_2(\FF_3)$ is not injective when restricted to $\Gamma_{O_2}(S_0)$.
Nevertheless, one can go one step further and consider the reduction map $\pi_{(3)} : \Gamma_{O_2}(S) \rightarrow \PGL_1(O_2/(3) O_2)$, which is injective when restricted to $\Gamma_{O_2}(S_0)$.
Since $(3) = \mathfrak{p}_3^2$, the group $\PGL_1(O_2/(3) O_2)$ is isomorphic to $\PGL_2(\FF_3[X]/(X^2))$, which has order $648$.
By carefully analysing the group, one can find a complement (of order $54$) to the image of $\Gamma_{O_2}(S_0)$ in $\PGL_1(O_2/(3) O_2)$, giving us a BMW family.

\vspace{5mm}

\begin{longtable}{r|l}
$2$ & \begin{tabular}{l}
$1/2(\alpha -2)+1/2j+1/2(\alpha -3)k,$\\
$1/3(\alpha -2)i+1/3(-2\alpha +4)j+1/3(\alpha -2)k,$\\
$1/6(\alpha -5)i+1/6(-2\alpha +1)j+1/6(\alpha +4)k,$\\
$1/2(-\alpha +2)i+1/2(-\alpha +3)j-1/2k,$\\
$1/2(\alpha -2)-1/2j+1/2(-\alpha +3)k$
\end{tabular}\\
\hline
$3+\alpha$ & \begin{tabular}{l}
$1/2(-\alpha -1)+1/3(\alpha +1)i+1/6(-\alpha +2)j+1/6(2\alpha -1)k,$\\
$1/2(-\alpha +2)+1/6(\alpha +1)i+1/3(-\alpha +2)j+1/6(-2\alpha +1)k,$\\
$1/6(-\alpha -1)i+1/6(-4\alpha -1)j+1/6(-\alpha -4)k,$\\
$1/6(4\alpha -5)i+1/6(\alpha -2)j+1/6(\alpha -5)k,$\\
$1/6(\alpha +7)i+1/6(-2\alpha -5)j+1/6(\alpha -2)k,$\\
$1/6(-\alpha -1)i+1/6(2\alpha -7)j+1/6(-\alpha +8)k,$\\
$1/2(-\alpha +2)+1/6(-\alpha -1)i+1/3(\alpha -2)j+1/6(2\alpha -1)k,$\\
$1/2(-\alpha -1)+1/3(-\alpha -1)i+1/6(\alpha -2)j+1/6(-2\alpha +1)k$
\end{tabular}
\end{longtable}

\begin{align*}
\Gamma_{2,3+\alpha} =& \langle a_{1}, a_{2}, a_{3}, a_{4},\hspace{1mm} b_{1}, b_{2}, b_{3}, b_{4}, b_{5}, b_{6}\st a_{2}^2,\hspace{1mm} a_{3}^2,\hspace{1mm} a_{4}^2,\hspace{1mm} b_{3}^2,\hspace{1mm} b_{4}^2,\hspace{1mm} b_{5}^2,\hspace{1mm} b_{6}^2,\hspace{1mm} a_{1}b_{1}a_{2}b_{1}, \hspace{1mm} a_{1}b_{2}a_{2}b_{4}, \hspace{1mm} \\*
 &\hspace{5mm} a_{1}b_{3}a_{3}b_{2}, \hspace{1mm} a_{1}b_{4}a_{3}b_{1}^{-1}, \hspace{1mm} a_{1}b_{5}a_{1}^{-1}b_{3}, \hspace{1mm} a_{1}b_{6}a_{1}b_{2}^{-1}, \hspace{1mm} a_{1}b_{1}^{-1}a_{4}b_{5}, \hspace{1mm} a_{2}b_{2}a_{4}b_{2}^{-1}, \hspace{1mm} a_{2}b_{3}a_{4}b_{3}, \hspace{1mm} \\*
 &\hspace{5mm} a_{2}b_{5}a_{3}b_{6}, \hspace{1mm} a_{3}b_{1}a_{4}b_{2}, \hspace{1mm} a_{4}b_{4}a_{4}b_{6}\rangle
\end{align*}

\subsection{\texorpdfstring{$(-1,-1 \st \QQ(\sqrt{69}))$}{(-1,-1 | Q(sqrt69))}}\label{subsec:-1_-1_Qsqrt69}

Here, as in \Cref{subsec:-1_-1_Qsqrt21}, the type number of the quaternion algebra we consider is greater than $1$.
Moreover, this example is particular as there exist sets $S$ for which $\Gamma_O(S)$ is vertex-transitive but contains no BMW subgroup.

Take $K = \QQ(\sqrt{69})$, $R = \ZZ[\alpha]$ where $\alpha = (1+\sqrt{69})/2$, and $B = (-1,-1 \st K)$.
The class number of $K$ is $1$, its narrow class group has order $2$, and $\Ram(B) = V^K_\infty$.
The type number of $B$ is $5$, but there is only one conjugacy class of maximal orders such that $V^K(O)$ is non-empty, and a representative is given by
\[O = R \oplus R i \oplus R \frac{1}{2}(-1+(-5+\alpha) i + (4-\alpha) j) \oplus R \frac{3+\alpha}{10} (5+5i+7j+k).\]
One has $\Idl(O) = \PIdl(O)$, so that $\Gamma_{O}(S)$ acts vertex-transitively on the product of trees with stabilisers conjugate to $\Gamma_{O}(S_0)$ by \Cref{lemma:VertexStabGammaS} and \Cref{cor:GammaSTransIfClass1}, upon the condition that $S\cap V^K_f \subseteq V^K(O)$.
Recall \Cref{cor:VKOtypeGt1} where we showed that $V^K(O)$ is the set of totally positive principal prime ideals.

One can show that $\Gamma_{O}(S_0)$ has order $2$, generated by $i$.
The reduced norm of $i$ is $1$, so $\Gamma_{O}(S_0)$ embeds in $\PSL_1(O/\mathfrak{p}O)\simeq \PSL_2(R/\mathfrak{p})$ for all prime $R$-ideal $\mathfrak{p}$ not in $S$.
In all cases, there is no index-$2$ subgroup of $\PGL_1(\localiz{O}{S}/\mathfrak{p})$ not containing $i$.
Unlike the other examples presented in this section, one cannot find a BMW subgroup as a congruence subgroup.

By enumerating all the index-2 subgroups of $\Gamma_O(S)$ for $S = V^K_\infty \cup \{2, -5+\alpha\}$, one can also show that none of them is vertex-transitive.

\subsection{\texorpdfstring{$(-75-19\sqrt{15}, -4-\sqrt{15} \st \QQ(\sqrt{15}))$}{(-75-19sqrt15,-4-sqrt15 | Q(sqrt15))}}\label{subsec:-75-19al_-4-al_Qsqrt15}

Take $K = \QQ(\sqrt{15})$, $R = \ZZ[\sqrt{15}]$, $B = (-75-19\sqrt{15}, -4-\sqrt{15} \st K)$.
This is the only quaternion algebra with type number $1$ over a number field with a non-trivial class group: $\Cl(K)$ has order $2$.
One has $\Ram(B) = \{\infty_1, \infty_2, \mathfrak{p}_2, \mathfrak{p}_3\}$, where $\infty_1$ and $\infty_2$ are the two places at infinity, $\mathfrak{p}_2 = (2, 1+\sqrt{15})$ and $\mathfrak{p}_3 = (3,\sqrt{15})$.
Note that $\mathfrak{p}_2$ and $\mathfrak{p}_3$ both generate the class group of $K$, and $\mathfrak{p}_2^2 = (2)$ and $\mathfrak{p}_3^2 = (3)$.
Let $O$ be the following maximal order:
\[R \oplus R i \oplus R \frac{1+i+j}{2} \oplus R \frac{5+2\sqrt{15}}{70} ((15+\sqrt{15}) + 19 i + k).\]
The quotient $\Idl(O)/\PIdl(O)$ has order $8$ and is generated by the prime two-sided $O$-ideals of reduced norms $\mathfrak{p}_2$ and $\mathfrak{p}_3$.
It follows that $\Idl(\localiz{O}{\Ram(B)}) = \PIdl(\localiz{O}{\Ram(B)})$.
One also observes that $\PGL_1(\localiz{O}{\Ram(B)}) = \PGL_1(O)$ has order $2$, generated by $j$.
By \Cref{lemma:VertexStabGammaS} and \Cref{cor:GammaSTransIfClass1}, $\Gamma_{O}(S)$ acts vertex-transitively with stabilisers conjugate to $\Gamma_{O}(S_0)$ if $\Ram(B)\subseteq S$.

To find BMW families, take a prime ideal $\mathfrak{p}$ not belonging to $S$ such that $\Nrd(j) = 4+\sqrt{15}$ is not a square modulo $\mathfrak{p}$.
Take $\mathfrak{p}$ to be $(7, 1+\sqrt{15}), (7, 1-\sqrt{15}), (2+\sqrt{15})$ or $(2-\sqrt{15})$, for example.
Then $\pi_\mathfrak{p}^{-1}(\PSL_2(R/\mathfrak{p}))$ acts regularly on the product of trees, where $\pi_\mathfrak{p} : \Gamma_{O}(S) \rightarrow \PGL_1(O/\mathfrak{p}O) \simeq \PGL_2(R/\mathfrak{p})$ is the reduction map.

\vspace{5mm}

\begin{longtable}{r|l}
$(5,\alpha)$ & \begin{tabular}{l}
$1-j,\hspace{1mm} (\alpha +4)+1/70(\alpha -15)i+1/70(-11\alpha +25)k,$\\
$1/10(\alpha +5)i-j+1/10(-\alpha +5)k,\hspace{1mm} 1/10(-\alpha -5)i-j+1/10(\alpha -5)k,$\\
$(\alpha +4)+1/70(-\alpha +15)i+1/70(11\alpha -25)k,\hspace{1mm} 1+j$
\end{tabular}\\
\hline
$(7,1+\alpha)$ & \begin{tabular}{l}
$1/2(-\alpha +3)+1/70(-13\alpha +55)i+1/2(-\alpha +3)j+1/35(19\alpha -75)k,$\\
$1/2(-\alpha +3)+1/70(13\alpha -55)i+1/2(-\alpha +3)j+1/35(-19\alpha +75)k,$\\
$1/2(-\alpha -3)+1/70(2\alpha +5)i+1/70(-22\alpha +85)k,$\\
$1/10(\alpha -5)i+1/10(-\alpha +5)k,\hspace{1mm} 1/35(-9\alpha +30)i+1/35(-41\alpha +160)k,$\\
$1/2(-\alpha -3)+1/70(-2\alpha -5)i+1/70(22\alpha -85)k,$\\
$1/2(-\alpha +3)+1/70(-13\alpha +55)i+1/2(\alpha -3)j+1/35(19\alpha -75)k,$\\
$1/2(-\alpha +3)+1/70(13\alpha -55)i+1/2(\alpha -3)j+1/35(-19\alpha +75)k$
\end{tabular}\\
\hline
$(7,6+\alpha)$ & \begin{tabular}{l}
$1/2(-\alpha -3)+1/70(9\alpha +5)i+1/2(-\alpha -3)j+1/35(3\alpha -10)k,$\\
$1/2(-\alpha -3)+1/70(-9\alpha -5)i+1/2(-\alpha -3)j+1/35(-3\alpha +10)k,$\\
$1/70(-11\alpha +25)i+1/70(-19\alpha +75)k,$\\
$1/35(-3\alpha +10)i+1/35(-2\alpha -5)k,$\\
$1/70(-4\alpha +25)i+1/2(\alpha -3)j+1/70(9\alpha -30)k,$\\
$1/70(-4\alpha +25)i+1/2(-\alpha +3)j+1/70(9\alpha -30)k,$\\
$1/2(-\alpha -3)+1/70(9\alpha +5)i+1/2(\alpha +3)j+1/35(3\alpha -10)k,$\\
$1/2(-\alpha -3)+1/70(-9\alpha -5)i+1/2(\alpha +3)j+1/35(-3\alpha +10)k$
\end{tabular}
\end{longtable}

\begin{align*}
\Gamma_{(5,\alpha),(7,1+\alpha)} =& \langle a_{1}, a_{2}, a_{3}, a_{4},\hspace{1mm} b_{1}, b_{2}, b_{3}, b_{4}, b_{5}\st a_{3}^2,\hspace{1mm} a_{4}^2,\hspace{1mm} b_{4}^2,\hspace{1mm} b_{5}^2,\hspace{1mm} a_{1}b_{1}a_{1}b_{2}, \hspace{1mm} a_{1}b_{3}a_{2}b_{1}^{-1}, \hspace{1mm} a_{1}b_{4}a_{1}b_{4}, \hspace{1mm} \\*
 &\hspace{5mm} a_{1}b_{5}a_{1}b_{5}, \hspace{1mm} a_{1}b_{3}^{-1}a_{2}^{-1}b_{2}^{-1}, \hspace{1mm} a_{1}b_{2}^{-1}a_{4}b_{3}^{-1}, \hspace{1mm} a_{1}b_{1}^{-1}a_{3}b_{3}, \hspace{1mm} a_{2}b_{1}a_{3}b_{5}, \hspace{1mm} a_{2}b_{2}a_{2}^{-1}b_{1}, \hspace{1mm} \\*
 &\hspace{5mm} a_{2}b_{4}a_{2}b_{3}^{-1}, \hspace{1mm} a_{2}b_{5}a_{4}b_{2}^{-1}, \hspace{1mm} a_{3}b_{2}a_{3}b_{2}, \hspace{1mm} a_{3}b_{4}a_{3}b_{4}, \hspace{1mm} a_{3}b_{3}^{-1}a_{4}b_{3}, \hspace{1mm} a_{4}b_{1}a_{4}b_{1}, \hspace{1mm} a_{4}b_{4}a_{4}b_{4}\rangle
\\
\Gamma_{(5,\alpha),(7,6+\alpha)} =& \langle a_{1}, a_{2}, a_{3}, a_{4},\hspace{1mm} b_{1}, b_{2}, b_{3}, b_{4}, b_{5}, b_{6}\st a_{3}^2,\hspace{1mm} a_{4}^2,\hspace{1mm} b_{3}^2,\hspace{1mm} b_{4}^2,\hspace{1mm} b_{5}^2,\hspace{1mm} b_{6}^2,\hspace{1mm} a_{1}b_{1}a_{1}b_{2}, \hspace{1mm} a_{1}b_{3}a_{1}b_{3}, \hspace{1mm} \\*
 &\hspace{5mm} a_{1}b_{4}a_{1}b_{4}, \hspace{1mm} a_{1}b_{5}a_{3}b_{2}^{-1}, \hspace{1mm} a_{1}b_{6}a_{4}b_{1}^{-1}, \hspace{1mm} a_{1}b_{2}^{-1}a_{2}^{-1}b_{5}, \hspace{1mm} a_{1}b_{1}^{-1}a_{2}b_{6}, \hspace{1mm} a_{2}b_{1}a_{2}b_{1}, \hspace{1mm} a_{2}b_{3}a_{2}b_{3}, \hspace{1mm} \\*
 &\hspace{5mm} a_{2}b_{4}a_{4}b_{2}, \hspace{1mm} a_{2}b_{5}a_{2}^{-1}b_{6}, \hspace{1mm} a_{2}b_{2}^{-1}a_{2}b_{2}^{-1}, \hspace{1mm} a_{2}b_{1}^{-1}a_{3}b_{4}, \hspace{1mm} a_{3}b_{2}a_{4}b_{1}, \hspace{1mm} a_{3}b_{3}a_{3}b_{6}, \hspace{1mm} a_{4}b_{3}a_{4}b_{5}\rangle
\\
\Gamma_{(7,1+\alpha),(7,6+\alpha)} =& \langle a_{1}, a_{2}, a_{3}, a_{4}, a_{5},\hspace{1mm} b_{1}, b_{2}, b_{3}, b_{4}, b_{5}, b_{6}\st a_{4}^2,\hspace{1mm} a_{5}^2,\hspace{1mm} b_{3}^2,\hspace{1mm} b_{4}^2,\hspace{1mm} b_{5}^2,\hspace{1mm} b_{6}^2,\hspace{1mm} a_{1}b_{1}a_{3}^{-1}b_{3}, \hspace{1mm} a_{1}b_{2}a_{3}^{-1}b_{6}, \hspace{1mm} \\*
 &\hspace{5mm} a_{1}b_{3}a_{1}^{-1}b_{5}, \hspace{1mm} a_{1}b_{4}a_{1}b_{1}^{-1}, \hspace{1mm} a_{1}b_{5}a_{3}b_{2}, \hspace{1mm} a_{1}b_{6}a_{4}b_{1}, \hspace{1mm} a_{1}b_{2}^{-1}a_{5}b_{2}^{-1}, \hspace{1mm} a_{2}b_{1}a_{3}b_{5}, \hspace{1mm} a_{2}b_{2}a_{3}b_{3}, \hspace{1mm} \\*
 &\hspace{5mm} a_{2}b_{3}a_{2}^{-1}b_{6}, \hspace{1mm} a_{2}b_{4}a_{2}b_{2}^{-1}, \hspace{1mm} a_{2}b_{5}a_{4}b_{2}, \hspace{1mm} a_{2}b_{6}a_{3}^{-1}b_{1}, \hspace{1mm} a_{2}b_{1}^{-1}a_{5}b_{1}^{-1}, \hspace{1mm} a_{3}b_{1}a_{4}b_{2}^{-1}, \hspace{1mm} \\*
 &\hspace{5mm} a_{3}b_{4}a_{3}b_{4}, \hspace{1mm} a_{4}b_{3}a_{5}b_{4}, \hspace{1mm} a_{5}b_{5}a_{5}b_{5}, \hspace{1mm} a_{5}b_{6}a_{5}b_{6}\rangle
\end{align*}

\section{Examples over global function fields}\label{sec:ExFuncFields}

Let $q > 1$ be a prime power, let $\infty$ be the valuation of $K = \FF_q(t)$ equal to $-\deg$, so that $R = \FF_q[t]$, and let $f\in R$ be an irreducible monic polynomial of degree $d\in \{1,2\}$.
Let $B$ the quaternion algebra over $K$ with $\Ram(B) = \{\infty, f\}$ and let $O$ be a maximal order.

By \Cref{thm:poschartypenum}, we know that $B$ has type number $1$.
Moreover, by \cite{Gekeler90}[Théorème 1 (iii)], the mass of $O$ is 
\[\mass(O) = \frac{q^d-1}{q^2-1},\]
where the mass of $O$ is defined by
\[\mass(O) = \sum_{[I]\in \ClR(O)} [O_\mathrm{L}(I)^\times : R^\times]^{-1}.\]
One can say more: by \cite{Gekeler90}[Théorème 1 (i),(ii)], the right class set $\ClR(O)$ is a singleton, implying that $\mass(O) = [O^\times : R^\times]^{-1}$.
In other words, the group $\PGL_1(O)\simeq \Gamma_O(V^K_\infty)$ has order $q+1$ if $d = 1$, and order $1$ if $d=2$.

Since the right class set $\ClR(O)$ is trivial, the quotient $\Idl(O)/\PIdl(O)$ is trivial too.
By \Cref{lemma:VertexStabGammaS} and \Cref{cor:GammaSTransIfClass1}, the group $\Gamma_{O}(S)$ acts vertex-transitively on the product of trees with stabilisers conjugate to $\Gamma_{O}(S_0)$.
Note also that $\Gamma_{O}(V^K_\infty)$ has index $2$ in $\Gamma_{O}(\Ram(B))$ by the exact sequence in \Cref{prop:ExactSeqGammaScupRam}.

\subsection{\texorpdfstring{Quaternion algebra over $\FF_q(t)$ with $\Ram(B) = \{\infty, f\}$ and $f$ irreducible of degree $1$}{Quaternion algebra over Fq(t) with Ram(B) = \{infty, f\} and f irreducible of degree 1}} \label{subsec:ExFqtDeg1}

Let $L$ be the separable quadratic extension $\FF_{q^2}(t)$ of $K = \FF_q(t)$.
Let $i$ be an element of $\FF_{q^2}$ not in $\FF_q$ such that $\alpha = i^2$ belongs to $\FF_q$ if $q$ is odd, and such that $\alpha = i^2 + i$ belongs to $\FF_q$ if $q$ is even.
The quaternion algebra $B$ is isomorphic to $(L, f \st K)$, that is to \[(\alpha, f \st K)\] in the odd case, and \[[\alpha , f \st K)\] in the even case.

The order with basis $1,i,j,k$ is maximal.
One sees that $\Gamma_O(V^K_\infty) = \PGL_1(O) = \{ i \} \cup \{1 + \lambda i \st \lambda\in \FF_q \}$.
On the other hand, an element $x = a + b i + c j + d k$ of the algebra $O/(f)O$ is invertible if and only if its reduced norm is non-zero.
In the odd case, we have $\Nrd(x) = a^2 - \alpha b^2 - f c^2 + \alpha f d^2 \equiv a^2 - \alpha b^2$, and in the even case we have $\Nrd(x) = a^2 +  a b + \alpha b^2 + f c^2 + f c d + \alpha f d^2 \equiv a^2 +  a b + \alpha b^2 $. 
In both cases, by the definition of $\alpha$, there are $(q^2-1)q^2$ invertible elements in $O/(f)O$, so that $\PGL_1(O/(f))$ has order $(q+1)q^2$.
If $\pi : \Gamma_O(S) \rightarrow \PGL_1(O/(f)O)$  is the reduction map, assuming that $f\notin S$, it is clear that $\pi$ maps $\Gamma_O(S_0)$ injectively.
If $P$ is any $p$-Sylow subgroup of $\PGL_1(O/(f)O)$, where $q$ is a power of the prime $p$, then $\pi^{-1}(P)$ acts regularly on the product of trees by \Cref{lemma:Htransfromquotient}.

\begin{remark}
Rungtanapirom Stix and Vdovina constructed these lattices in their articles \cite{StixVdov17} and \cite{RungStixVdov19}.
This is also a generalisation of \cite{Morgenstern94}, passing from an action on one tree to an action on a product of trees.
\end{remark}

\paragraph{Family with $q = 3$.}
The parameters are 
\[q=3, \alpha = 2, f = t.\]

\begin{longtable}{r|l}
$t + 1$ & \begin{tabular}{l}
$2+2j+2k,\hspace{1mm} 1+j+2k,\hspace{1mm} 1+2j+k,\hspace{1mm} 2+j+k$
\end{tabular}\\
\hline
$t + 2$ & \begin{tabular}{l}
$1+2j,\hspace{1mm} 1+k,\hspace{1mm} 1+2k,\hspace{1mm} 1+j$
\end{tabular}
\end{longtable}

\begin{align*}
\Gamma_{t + 1,t + 2} =& \langle a_{1}, a_{2},\hspace{1mm} b_{1}, b_{2}\st a_{1}b_{1}a_{2}b_{1}, \hspace{1mm} a_{1}b_{2}a_{1}b_{1}^{-1}, \hspace{1mm} a_{1}b_{2}^{-1}a_{2}^{-1}b_{2}^{-1}, \hspace{1mm} a_{2}b_{2}^{-1}a_{2}b_{1}^{-1}\rangle
\end{align*}

\paragraph{Family with $q = 5$.}
The parameters are 
\[q=5, \alpha = 2, f = t.\]

\begin{longtable}{r|l}
$t + 1$ & \begin{tabular}{l}
$3+3j+3k,\hspace{1mm} 2+2j+3k,\hspace{1mm} 1+2j,\hspace{1mm} 1+3j,\hspace{1mm} 2+3j+2k,\hspace{1mm} 3+2j+2k$
\end{tabular}\\
\hline
$t + 2$ & \begin{tabular}{l}
$4+3j+k,\hspace{1mm} 3+k,\hspace{1mm} 3+j+3k,\hspace{1mm} 3+4j+2k,\hspace{1mm} 3+4k,\hspace{1mm} 4+2j+4k$
\end{tabular}\\
\hline
$t + 3$ & \begin{tabular}{l}
$1+k,\hspace{1mm} 1+j+3k,\hspace{1mm} 2+3j+k,\hspace{1mm} 2+2j+4k,\hspace{1mm} 1+4j+2k,\hspace{1mm} 1+4k$
\end{tabular}\\
\hline
$t + 4$ & \begin{tabular}{l}
$4+3j+3k,\hspace{1mm} 4+2j+3k,\hspace{1mm} 3+2j,\hspace{1mm} 3+3j,\hspace{1mm} 4+3j+2k,\hspace{1mm} 4+2j+2k$
\end{tabular}
\end{longtable}

\begin{align*}
\Gamma_{t + 1,t + 2} =& \langle a_{1}, a_{2}, a_{3},\hspace{1mm} b_{1}, b_{2}, b_{3}\st a_{1}b_{1}a_{2}^{-1}b_{3}^{-1}, \hspace{1mm} a_{1}b_{2}a_{1}b_{1}^{-1}, \hspace{1mm} a_{1}b_{3}a_{3}b_{2}^{-1}, \hspace{1mm} a_{1}b_{3}^{-1}a_{2}^{-1}b_{1}, \hspace{1mm} a_{1}b_{2}^{-1}a_{3}b_{3}, \hspace{1mm} \\*
 &\hspace{5mm} a_{2}b_{1}a_{3}b_{2}, \hspace{1mm} a_{2}b_{2}a_{3}b_{1}, \hspace{1mm} a_{2}b_{3}^{-1}a_{2}b_{2}^{-1}, \hspace{1mm} a_{3}b_{3}^{-1}a_{3}b_{1}^{-1}\rangle
\\
\Gamma_{t + 1,t + 4} =& \langle a_{1}, a_{2}, a_{3},\hspace{1mm} b_{1}, b_{2}, b_{3}\st a_{1}b_{1}a_{1}^{-1}b_{1}^{-1}, \hspace{1mm} a_{1}b_{2}a_{3}^{-1}b_{2}, \hspace{1mm} a_{1}b_{3}a_{2}b_{3}, \hspace{1mm} a_{1}b_{3}^{-1}a_{1}b_{2}^{-1}, \hspace{1mm} a_{2}b_{1}a_{2}b_{3}^{-1}, \hspace{1mm} \\*
 &\hspace{5mm} a_{2}b_{2}a_{2}^{-1}b_{2}^{-1}, \hspace{1mm} a_{2}b_{1}^{-1}a_{3}^{-1}b_{1}^{-1}, \hspace{1mm} a_{3}b_{2}a_{3}b_{1}^{-1}, \hspace{1mm} a_{3}b_{3}a_{3}^{-1}b_{3}^{-1}\rangle
\end{align*}

\paragraph{Family with $q = 7$.}
The parameters are 
\[q=7, \alpha = 3, f = t.\]

\begin{longtable}{r|l}
$t + 1$ & \begin{tabular}{l}
$4+5j+4k,\hspace{1mm} 4+2j+4k,\hspace{1mm} 1+2j+5k,\hspace{1mm} 5+3j+3k,\hspace{1mm} 5+4j+4k,$\\
$1+5j+2k,\hspace{1mm} 4+5j+3k,\hspace{1mm} 4+2j+3k$
\end{tabular}\\
\hline
$t + 2$ & \begin{tabular}{l}
$2+5j+4k,\hspace{1mm} 5+2j+4k,\hspace{1mm} 3+2j+5k,\hspace{1mm} 6+3j+3k,\hspace{1mm} 6+4j+4k,$\\
$3+5j+2k,\hspace{1mm} 5+5j+3k,\hspace{1mm} 2+2j+3k$
\end{tabular}\\
\hline
$t + 3$ & \begin{tabular}{l}
$3+4j+5k,\hspace{1mm} 3+k,\hspace{1mm} 5+2j+k,\hspace{1mm} 6+4j,\hspace{1mm} 6+3j,\hspace{1mm} 5+5j+6k,\hspace{1mm} 3+6k,$\\
$3+3j+2k$
\end{tabular}\\
\hline
$t + 4$ & \begin{tabular}{l}
$6+2j+4k,\hspace{1mm} 6+5j+4k,\hspace{1mm} 4+3j+3k,\hspace{1mm} 2+2j+5k,\hspace{1mm} 2+5j+2k,$\\
$4+4j+4k,\hspace{1mm} 6+2j+3k,\hspace{1mm} 6+5j+3k$
\end{tabular}\\
\hline
$t + 5$ & \begin{tabular}{l}
$1+k,\hspace{1mm} 4+2j+k,\hspace{1mm} 1+4j+5k,\hspace{1mm} 2+4j,\hspace{1mm} 2+3j,\hspace{1mm} 1+3j+2k,\hspace{1mm} 4+5j+6k,$\\
$1+6k$
\end{tabular}\\
\hline
$t + 6$ & \begin{tabular}{l}
$6+2j+k,\hspace{1mm} 2+4j+5k,\hspace{1mm} 4+4j,\hspace{1mm} 5+k,\hspace{1mm} 5+6k,\hspace{1mm} 4+3j,\hspace{1mm} 2+3j+2k,$\\
$6+5j+6k$
\end{tabular}
\end{longtable}

\begin{align*}
\Gamma_{t + 1,t + 2} =& \langle a_{1}, a_{2}, a_{3}, a_{4},\hspace{1mm} b_{1}, b_{2}, b_{3}, b_{4}\st a_{1}b_{1}a_{1}^{-1}b_{1}^{-1}, \hspace{1mm} a_{1}b_{2}a_{3}b_{4}^{-1}, \hspace{1mm} a_{1}b_{3}a_{4}^{-1}b_{3}, \hspace{1mm} a_{1}b_{4}a_{2}^{-1}b_{3}^{-1}, \hspace{1mm} a_{1}b_{4}^{-1}a_{3}b_{2}, \hspace{1mm} \\*
 &\hspace{5mm} a_{1}b_{3}^{-1}a_{2}^{-1}b_{4}, \hspace{1mm} a_{1}b_{2}^{-1}a_{4}b_{2}^{-1}, \hspace{1mm} a_{2}b_{1}a_{3}^{-1}b_{1}, \hspace{1mm} a_{2}b_{2}a_{2}^{-1}b_{2}^{-1}, \hspace{1mm} a_{2}b_{4}a_{3}b_{4}, \hspace{1mm} a_{2}b_{3}^{-1}a_{4}^{-1}b_{1}^{-1}, \hspace{1mm} \\*
 &\hspace{5mm} a_{2}b_{1}^{-1}a_{4}^{-1}b_{3}^{-1}, \hspace{1mm} a_{3}b_{1}a_{4}^{-1}b_{2}^{-1}, \hspace{1mm} a_{3}b_{3}a_{3}^{-1}b_{3}^{-1}, \hspace{1mm} a_{3}b_{2}^{-1}a_{4}^{-1}b_{1}, \hspace{1mm} a_{4}b_{4}a_{4}^{-1}b_{4}^{-1}\rangle
\\
\Gamma_{t + 1,t + 3} =& \langle a_{1}, a_{2}, a_{3}, a_{4},\hspace{1mm} b_{1}, b_{2}, b_{3}, b_{4}\st a_{1}b_{1}a_{4}b_{4}, \hspace{1mm} a_{1}b_{2}a_{4}^{-1}b_{3}, \hspace{1mm} a_{1}b_{3}a_{4}^{-1}b_{2}, \hspace{1mm} a_{1}b_{4}a_{4}b_{1}, \hspace{1mm} a_{1}b_{4}^{-1}a_{2}^{-1}b_{4}^{-1}, \hspace{1mm} \\*
 &\hspace{5mm} a_{1}b_{3}^{-1}a_{3}b_{3}^{-1}, \hspace{1mm} a_{1}b_{2}^{-1}a_{1}b_{1}^{-1}, \hspace{1mm} a_{2}b_{1}a_{4}^{-1}b_{1}, \hspace{1mm} a_{2}b_{2}a_{3}b_{1}^{-1}, \hspace{1mm} a_{2}b_{3}a_{2}b_{2}^{-1}, \hspace{1mm} a_{2}b_{4}^{-1}a_{3}^{-1}b_{3}^{-1}, \hspace{1mm} \\*
 &\hspace{5mm} a_{2}b_{3}^{-1}a_{3}^{-1}b_{4}^{-1}, \hspace{1mm} a_{2}b_{1}^{-1}a_{3}b_{2}, \hspace{1mm} a_{3}b_{1}a_{3}b_{4}^{-1}, \hspace{1mm} a_{3}b_{2}^{-1}a_{4}^{-1}b_{2}^{-1}, \hspace{1mm} a_{4}b_{3}a_{4}b_{4}^{-1}\rangle
\\
\Gamma_{t + 1,t + 6} =& \langle a_{1}, a_{2}, a_{3}, a_{4},\hspace{1mm} b_{1}, b_{2}, b_{3}, b_{4}\st a_{1}b_{1}a_{2}b_{2}, \hspace{1mm} a_{1}b_{2}a_{2}b_{1}, \hspace{1mm} a_{1}b_{3}a_{4}^{-1}b_{2}^{-1}, \hspace{1mm} a_{1}b_{4}a_{4}b_{1}^{-1}, \hspace{1mm} a_{1}b_{4}^{-1}a_{3}^{-1}b_{3}^{-1}, \hspace{1mm} \\*
 &\hspace{5mm} a_{1}b_{3}^{-1}a_{3}^{-1}b_{4}^{-1}, \hspace{1mm} a_{1}b_{2}^{-1}a_{4}^{-1}b_{3}, \hspace{1mm} a_{1}b_{1}^{-1}a_{4}b_{4}, \hspace{1mm} a_{2}b_{3}a_{4}b_{4}^{-1}, \hspace{1mm} a_{2}b_{4}a_{3}^{-1}b_{2}^{-1}, \hspace{1mm} a_{2}b_{4}^{-1}a_{4}b_{3}, \hspace{1mm} \\*
 &\hspace{5mm} a_{2}b_{3}^{-1}a_{3}b_{1}^{-1}, \hspace{1mm} a_{2}b_{2}^{-1}a_{3}^{-1}b_{4}, \hspace{1mm} a_{2}b_{1}^{-1}a_{3}b_{3}^{-1}, \hspace{1mm} a_{3}b_{1}a_{4}b_{2}^{-1}, \hspace{1mm} a_{3}b_{2}^{-1}a_{4}b_{1}\rangle
\end{align*}

\subsection{\texorpdfstring{Quaternion algebra over $\FF_q(t)$ with $\Ram(B) = \{\infty, f\}$ and $f$ irreducible of degree $2$}{Quaternion algebra over Fq(t) with Ram(B) = \{infty, f\} and f irreducible of degree 2}}\label{subsec:ExFqtDeg2}

The group $\Gamma_O(V^K_\infty) = \PGL_1(O)$ is trivial, so $\Gamma_O(S)$ acts regularly on the product of trees if $f\notin S$.
We give examples for this family below.

Alternatively, assume now that $f\in S$.
If $q$ is odd, one can find an irreducible polynomial $g\in R$ not in $S$ such that $f$ is not a square modulo $g$.
If $\pi_g : \Gamma_O(S) \rightarrow \PGL_1(O/(g)O) \simeq \PGL_2(R/(g))$ is the reduction map, it is surjective, and then $\pi_g^{-1}(\PSL_2(R/(g)))$ acts regularly on the product of trees by \Cref{lemma:Htransfromquotient}.

\begin{remark}
Rungtanapirom studied in \cite{Rungtanapirom18} the lattice $\Gamma = \Gamma_U(\{\infty,t,t+1\})$, where $U$ is a $\FF_2[(t+1)^{-1}]$-maximal order in the quaternion algebra $B$ over $\FF_2(t)$ which ramifies at $t+1$ and $t^2+t+1$.
He observed that $\Gamma$ acts regularly on the product of the Bruhat-Tits trees corresponding to the places $\infty$ and $t$.
Up to a field automorphism, this lattice is equivalent to the lattice $\Gamma_O(S)$ described just above, with $q=2$, $f = t^2+t+1$ and $S = \{\infty, t, t+1 \}$, acting regularly on the product of the Bruhat-Tits trees corresponding to the places $t$ and $t+1$.
\end{remark}

\paragraph{Family with $q = 3$.}
The parameters are 
\[q=3, f = t^2 + 1, B = (-t^2-1,t-1 \st K),\]
\[O = R\oplus R i \oplus R \frac{1}{t+2} (j+2k) \oplus R k\]

\begin{longtable}{r|l}
$t$ & \begin{tabular}{l}
$1+1/(t +2)j+2/(t +2)k,\hspace{1mm} 1+(t +1)/(t +2)j+1/(t +2)k,$\\
$1+(2t +2)/(t +2)j+2/(t +2)k,\hspace{1mm} 1+2/(t +2)j+1/(t +2)k$
\end{tabular}\\
\hline
$t + 1$ & \begin{tabular}{l}
$1+2j,\hspace{1mm} t /(t +2)j+2/(t +2)k,\hspace{1mm} 1/(t +2)j+2/(t +2)k,\hspace{1mm} 1+j$
\end{tabular}\\
\hline
$t + 2$ & \begin{tabular}{l}
$1+t /(t +2)j+2/(t +2)k,\hspace{1mm} (t +1)/(t +2)j+1/(t +2)k,\hspace{1mm} j,$\\
$1+2t /(t +2)j+1/(t +2)k$
\end{tabular}
\end{longtable}

\begin{align*}
\Gamma_{t,t + 1} =& \langle a_{1}, a_{2},\hspace{1mm} b_{1}, b_{2}, b_{3}\st b_{2}^2,\hspace{1mm} b_{3}^2,\hspace{1mm} a_{1}b_{1}a_{1}b_{1}, \hspace{1mm} a_{1}b_{2}a_{2}b_{1}^{-1}, \hspace{1mm} a_{1}b_{3}a_{1}^{-1}b_{3}, \hspace{1mm} a_{1}b_{1}^{-1}a_{2}b_{2}, \hspace{1mm} a_{2}b_{1}a_{2}b_{3}\rangle
\\
\Gamma_{t + 1,t + 2} =& \langle a_{1}, a_{2}, a_{3},\hspace{1mm} b_{1}, b_{2}, b_{3}\st a_{2}^2,\hspace{1mm} a_{3}^2,\hspace{1mm} b_{2}^2,\hspace{1mm} b_{3}^2,\hspace{1mm} a_{1}b_{1}a_{3}b_{1}, \hspace{1mm} a_{1}b_{2}a_{1}b_{1}^{-1}, \hspace{1mm} a_{1}b_{3}a_{1}^{-1}b_{3}, \hspace{1mm} a_{2}b_{1}a_{2}b_{1}^{-1}, \hspace{1mm} \\*
 &\hspace{5mm} a_{2}b_{2}a_{3}b_{3}\rangle
\end{align*}

\paragraph{Family with $q = 5$.}
The parameters are 
\[q=5, f = t^2 + 3, B = (3t^2+4,t+2 \st K),\]
\[O = R\oplus R i \oplus R \frac{1}{t+2} (j+4k) \oplus R k\]

\begin{longtable}{r|l}
$t$ & \begin{tabular}{l}
$1+(4t +4)/(t +2)j+4/(t +2)k,\hspace{1mm} 1+1/(t +2)j+4/(t +2)k,$\\
$(t +4)/(t +2)j+3/(t +2)k,\hspace{1mm} (t +3)/(t +2)j+4/(t +2)k,$\\
$1+4/(t +2)j+1/(t +2)k,\hspace{1mm} 1+(t +1)/(t +2)j+1/(t +2)k$
\end{tabular}\\
\hline
$t + 1$ & \begin{tabular}{l}
$3+1/(t +2)j+4/(t +2)k,\hspace{1mm} 1+(t +4)/(t +2)j+3/(t +2)k,$\\
$1+j,\hspace{1mm} 1+4j,\hspace{1mm} 1+(4t +1)/(t +2)j+2/(t +2)k,\hspace{1mm} 3+4/(t +2)j+1/(t +2)k$
\end{tabular}\\
\hline
$t + 2$ & \begin{tabular}{l}
$4+3t /(t +2)j+1/(t +2)k,\hspace{1mm} 4+(2t +1)/(t +2)j+3/(t +2)k,$\\
$j,\hspace{1mm} (t +1)/(t +2)j+1/(t +2)k,\hspace{1mm} 4+(3t +4)/(t +2)j+2/(t +2)k,$\\
$4+2t /(t +2)j+4/(t +2)k$
\end{tabular}\\
\hline
$t + 3$ & \begin{tabular}{l}
$1+(t +3)/(t +2)j+4/(t +2)k,\hspace{1mm} 4+2j,\hspace{1mm} 1/(t +2)j+4/(t +2)k,$\\
$t /(t +2)j+2/(t +2)k,\hspace{1mm} 4+3j,\hspace{1mm} 1+(4t +2)/(t +2)j+1/(t +2)k$
\end{tabular}\\
\hline
$t + 4$ & \begin{tabular}{l}
$1+4t /(t +2)j+3/(t +2)k,\hspace{1mm} 4+(3t +2)/(t +2)j+4/(t +2)k,$\\
$4+(3t +3)/(t +2)j+3/(t +2)k,\hspace{1mm} 4+(2t +2)/(t +2)j+2/(t +2)k,$\\
$4+(2t +3)/(t +2)j+1/(t +2)k,\hspace{1mm} 1+t /(t +2)j+2/(t +2)k$
\end{tabular}
\end{longtable}

\begin{align*}
\Gamma_{t,t + 1} =& \langle a_{1}, a_{2}, a_{3}, a_{4},\hspace{1mm} b_{1}, b_{2}, b_{3}\st a_{3}^2,\hspace{1mm} a_{4}^2,\hspace{1mm} a_{1}b_{1}a_{3}b_{3}^{-1}, \hspace{1mm} a_{1}b_{2}a_{1}b_{3}, \hspace{1mm} a_{1}b_{3}^{-1}a_{3}b_{1}, \hspace{1mm} a_{1}b_{2}^{-1}a_{2}^{-1}b_{2}^{-1}, \hspace{1mm} \\*
 &\hspace{5mm} a_{1}b_{1}^{-1}a_{4}b_{1}^{-1}, \hspace{1mm} a_{2}b_{1}a_{2}^{-1}b_{1}^{-1}, \hspace{1mm} a_{2}b_{3}a_{4}b_{2}^{-1}, \hspace{1mm} a_{2}b_{3}^{-1}a_{2}b_{3}^{-1}, \hspace{1mm} a_{2}b_{2}^{-1}a_{4}b_{3}, \hspace{1mm} a_{3}b_{2}a_{3}b_{2}^{-1}\rangle
\\
\Gamma_{t,t + 2} =& \langle a_{1}, a_{2}, a_{3}, a_{4},\hspace{1mm} b_{1}, b_{2}, b_{3}, b_{4}\st a_{3}^2,\hspace{1mm} a_{4}^2,\hspace{1mm} b_{3}^2,\hspace{1mm} b_{4}^2,\hspace{1mm} a_{1}b_{1}a_{3}b_{2}, \hspace{1mm} a_{1}b_{2}a_{3}b_{1}, \hspace{1mm} a_{1}b_{3}a_{2}b_{1}^{-1}, \hspace{1mm} \\*
 &\hspace{5mm} a_{1}b_{4}a_{1}^{-1}b_{4}, \hspace{1mm} a_{1}b_{2}^{-1}a_{2}^{-1}b_{2}^{-1}, \hspace{1mm} a_{1}b_{1}^{-1}a_{2}b_{3}, \hspace{1mm} a_{2}b_{1}a_{4}b_{1}, \hspace{1mm} a_{2}b_{4}a_{2}b_{2}^{-1}, \hspace{1mm} a_{3}b_{3}a_{4}b_{4}, \hspace{1mm} \\*
 &\hspace{5mm} a_{4}b_{2}a_{4}b_{2}^{-1}\rangle
\\
\Gamma_{t + 1,t + 4} =& \langle a_{1}, a_{2}, a_{3},\hspace{1mm} b_{1}, b_{2}, b_{3}\st a_{1}b_{1}a_{1}b_{2}^{-1}, \hspace{1mm} a_{1}b_{2}a_{3}b_{3}^{-1}, \hspace{1mm} a_{1}b_{3}a_{2}b_{1}^{-1}, \hspace{1mm} a_{1}b_{3}^{-1}a_{3}b_{2}, \hspace{1mm} a_{1}b_{1}^{-1}a_{2}b_{3}, \hspace{1mm} \\*
 &\hspace{5mm} a_{2}b_{1}a_{3}b_{1}, \hspace{1mm} a_{2}b_{2}a_{2}^{-1}b_{2}^{-1}, \hspace{1mm} a_{2}b_{3}^{-1}a_{3}^{-1}b_{3}^{-1}, \hspace{1mm} a_{3}b_{2}^{-1}a_{3}b_{1}^{-1}\rangle
\end{align*}

\paragraph{Family with $q = 7$.}
The parameters are
\[q=7, f = t^2 + 4, B = (6 t^2 + 3,t+1 \st K),\]
\[O = R\oplus R i \oplus R \frac{1}{t+1} (j+5k) \oplus R k\]

\begin{longtable}{r|l}
$t$ & \begin{tabular}{l}
$4+4j,\hspace{1mm} 4+(3t +1)/(t +1)j+4/(t +1)k,\hspace{1mm} 1+2t /(t +1)j+4/(t +1)k,$\\
$1+(2t +5)/(t +1)j+1/(t +1)k,\hspace{1mm} 1+(5t +2)/(t +1)j+6/(t +1)k,$\\
$1+5t /(t +1)j+3/(t +1)k,\hspace{1mm} 4+(4t +6)/(t +1)j+3/(t +1)k,$\\
$4+3j$
\end{tabular}\\
\hline
$t + 1$ & \begin{tabular}{l}
$2+(t +4)/(t +1)j+1/(t +1)k,\hspace{1mm} 2+t /(t +1)j+2/(t +1)k,$\\
$4+(3t +6)/(t +1)j+1/(t +1)k,\hspace{1mm} 4j,\hspace{1mm} (4t +5)/(t +1)j+5/(t +1)k,$\\
$4+(4t +1)/(t +1)j+6/(t +1)k,\hspace{1mm} 2+6t /(t +1)j+5/(t +1)k,$\\
$2+(6t +3)/(t +1)j+6/(t +1)k$
\end{tabular}\\
\hline
$t + 2$ & \begin{tabular}{l}
$2+(6t +2)/(t +1)j+1/(t +1)k,\hspace{1mm} 1+(2t +1)/(t +1)j+2/(t +1)k,$\\
$2+(6t +4)/(t +1)j+4/(t +1)k,\hspace{1mm} 3+1/(t +1)j+5/(t +1)k,$\\
$3+6/(t +1)j+2/(t +1)k,\hspace{1mm} 2+(t +3)/(t +1)j+3/(t +1)k,$\\
$1+(5t +6)/(t +1)j+5/(t +1)k,\hspace{1mm} 2+(t +5)/(t +1)j+6/(t +1)k$
\end{tabular}\\
\hline
$t + 3$ & \begin{tabular}{l}
$1+(2t +3)/(t +1)j+5/(t +1)k,\hspace{1mm} 4+3t /(t +1)j+6/(t +1)k,$\\
$4+(4t +5)/(t +1)j+5/(t +1)k,\hspace{1mm} (4t +2)/(t +1)j+4/(t +1)k,$\\
$(4t +3)/(t +1)j+2/(t +1)k,\hspace{1mm} 4+(3t +2)/(t +1)j+2/(t +1)k,$\\
$4+4t /(t +1)j+1/(t +1)k,\hspace{1mm} 1+(5t +4)/(t +1)j+2/(t +1)k$
\end{tabular}\\
\hline
$t + 4$ & \begin{tabular}{l}
$1+1/(t +1)j+5/(t +1)k,\hspace{1mm} 2+6j,\hspace{1mm} 2+(6t +1)/(t +1)j+3/(t +1)k,$\\
$(4t +6)/(t +1)j+3/(t +1)k,\hspace{1mm} (4t +1)/(t +1)j+6/(t +1)k,$\\
$2+(t +6)/(t +1)j+4/(t +1)k,\hspace{1mm} 2+j,\hspace{1mm} 1+6/(t +1)j+2/(t +1)k$
\end{tabular}\\
\hline
$t + 5$ & \begin{tabular}{l}
$1+(5t +1)/(t +1)j+1/(t +1)k,\hspace{1mm} 1+(5t +3)/(t +1)j+4/(t +1)k,$\\
$4+(3t +4)/(t +1)j+5/(t +1)k,\hspace{1mm} 5+1/(t +1)j+5/(t +1)k,$\\
$5+6/(t +1)j+2/(t +1)k,\hspace{1mm} 4+(4t +3)/(t +1)j+2/(t +1)k,$\\
$1+(2t +4)/(t +1)j+3/(t +1)k,\hspace{1mm} 1+(2t +6)/(t +1)j+6/(t +1)k$
\end{tabular}\\
\hline
$t + 6$ & \begin{tabular}{l}
$4+(3t +5)/(t +1)j+3/(t +1)k,\hspace{1mm} 2+(t +2)/(t +1)j+5/(t +1)k,$\\
$1+2j,\hspace{1mm} 4t /(t +1)j+1/(t +1)k,\hspace{1mm} 1/(t +1)j+5/(t +1)k,\hspace{1mm} 1+5j,$\\
$2+(6t +5)/(t +1)j+2/(t +1)k,\hspace{1mm} 4+(4t +2)/(t +1)j+4/(t +1)k$
\end{tabular}
\end{longtable}

\begin{align*}
\Gamma_{t,t + 1} =& \langle a_{1}, a_{2}, a_{3}, a_{4},\hspace{1mm} b_{1}, b_{2}, b_{3}, b_{4}, b_{5}\st b_{4}^2,\hspace{1mm} b_{5}^2,\hspace{1mm} a_{1}b_{1}a_{4}^{-1}b_{2}^{-1}, \hspace{1mm} a_{1}b_{2}a_{2}b_{5}, \hspace{1mm} a_{1}b_{3}a_{2}^{-1}b_{3}, \hspace{1mm} a_{1}b_{4}a_{1}^{-1}b_{4}, \hspace{1mm} \\*
 &\hspace{5mm} a_{1}b_{5}a_{2}b_{2}, \hspace{1mm} a_{1}b_{3}^{-1}a_{3}^{-1}b_{1}^{-1}, \hspace{1mm} a_{1}b_{2}^{-1}a_{4}^{-1}b_{1}, \hspace{1mm} a_{1}b_{1}^{-1}a_{3}^{-1}b_{3}^{-1}, \hspace{1mm} a_{2}b_{1}a_{2}b_{1}, \hspace{1mm} a_{2}b_{3}a_{4}^{-1}b_{1}^{-1}, \hspace{1mm} \\*
 &\hspace{5mm} a_{2}b_{4}a_{4}b_{2}^{-1}, \hspace{1mm} a_{2}b_{2}^{-1}a_{4}b_{4}, \hspace{1mm} a_{2}b_{1}^{-1}a_{4}^{-1}b_{3}, \hspace{1mm} a_{3}b_{2}a_{3}^{-1}b_{2}^{-1}, \hspace{1mm} a_{3}b_{4}a_{3}b_{1}^{-1}, \hspace{1mm} a_{3}b_{5}a_{4}^{-1}b_{3}^{-1}, \hspace{1mm} \\*
 &\hspace{5mm} a_{3}b_{3}^{-1}a_{4}^{-1}b_{5}\rangle
\\
\Gamma_{t,t + 2} =& \langle a_{1}, a_{2}, a_{3}, a_{4},\hspace{1mm} b_{1}, b_{2}, b_{3}, b_{4}\st a_{1}b_{1}a_{4}b_{3}^{-1}, \hspace{1mm} a_{1}b_{2}a_{2}b_{4}^{-1}, \hspace{1mm} a_{1}b_{3}a_{3}b_{4}, \hspace{1mm} a_{1}b_{4}a_{3}b_{3}, \hspace{1mm} a_{1}b_{4}^{-1}a_{2}b_{2}, \hspace{1mm} \\*
 &\hspace{5mm} a_{1}b_{3}^{-1}a_{4}b_{1}, \hspace{1mm} a_{1}b_{2}^{-1}a_{1}b_{1}^{-1}, \hspace{1mm} a_{2}b_{1}a_{2}^{-1}b_{1}^{-1}, \hspace{1mm} a_{2}b_{3}a_{3}^{-1}b_{2}^{-1}, \hspace{1mm} a_{2}b_{4}a_{4}^{-1}b_{3}^{-1}, \hspace{1mm} a_{2}b_{3}^{-1}a_{4}^{-1}b_{4}, \hspace{1mm} \\*
 &\hspace{5mm} a_{2}b_{2}^{-1}a_{3}^{-1}b_{3}, \hspace{1mm} a_{3}b_{1}a_{4}^{-1}b_{1}, \hspace{1mm} a_{3}b_{4}^{-1}a_{3}b_{1}^{-1}, \hspace{1mm} a_{3}b_{2}^{-1}a_{4}b_{2}^{-1}, \hspace{1mm} a_{4}b_{2}a_{4}b_{4}\rangle
\\
\Gamma_{t + 1,t + 3} =& \langle a_{1}, a_{2}, a_{3}, a_{4}, a_{5},\hspace{1mm} b_{1}, b_{2}, b_{3}, b_{4}, b_{5}\st a_{4}^2,\hspace{1mm} a_{5}^2,\hspace{1mm} b_{4}^2,\hspace{1mm} b_{5}^2,\hspace{1mm} a_{1}b_{1}a_{2}^{-1}b_{3}^{-1}, \hspace{1mm} a_{1}b_{2}a_{5}b_{1}^{-1}, \hspace{1mm} a_{1}b_{3}a_{4}b_{3}, \hspace{1mm} \\*
 &\hspace{5mm} a_{1}b_{4}a_{1}^{-1}b_{4}, \hspace{1mm} a_{1}b_{5}a_{1}b_{2}^{-1}, \hspace{1mm} a_{1}b_{3}^{-1}a_{2}^{-1}b_{1}, \hspace{1mm} a_{1}b_{1}^{-1}a_{5}b_{2}, \hspace{1mm} a_{2}b_{1}a_{3}^{-1}b_{4}, \hspace{1mm} a_{2}b_{2}a_{2}^{-1}b_{2}^{-1}, \hspace{1mm} \\*
 &\hspace{5mm} a_{2}b_{4}a_{3}^{-1}b_{1}, \hspace{1mm} a_{2}b_{5}a_{3}b_{3}^{-1}, \hspace{1mm} a_{2}b_{3}^{-1}a_{3}b_{5}, \hspace{1mm} a_{3}b_{1}a_{4}b_{2}, \hspace{1mm} a_{3}b_{2}a_{4}b_{1}, \hspace{1mm} a_{3}b_{3}a_{3}b_{2}^{-1}, \hspace{1mm} a_{4}b_{4}a_{5}b_{5}, \hspace{1mm} \\*
 &\hspace{5mm} a_{5}b_{3}a_{5}b_{3}^{-1}\rangle
\end{align*}

\bibliographystyle{alpha}
\bibliography{QuaternionicFamilies}

\appendix

\begin{landscape}
\section{Maximal orders of type number one in characteristic zero}\label{Appendix}

In the appendix, we give a table together with relevant information of the maximal orders $O$ of type number one in totally definite quaternion algebras over number fields, classified in \cite{KirschmerLorch16}.
See \Cref{subsec:TransLattType1} for more details.

\begin{itemize}
\item The first column is the quaternion algebra $B$.
Above each example, an irreducible polynomial $P \in \QQ[X]$ is given, and the base field of $B$ is the quotient $K = \QQ[X]/(P)$.
We denote by $\alpha$ the image of $X$ in $K$, which is a root of $P$.
\item The second column is the set $\Ram(B)\cap V^K_f$.
If the field has class number one, a place will be represented by an element in $K$ generating the prime ideal corresponding to the place.
Otherwise, we give the prime ideal.
\item The third column are the minimal sets of places $S_0\cap V^K_f$, such that $\Idl(\localiz{O}{S_0}) = \PIdl(\localiz{O}{S_0})$.
This column serves to apply \Cref{cor:GammaSTransIfClass1}: the lattice $\Gamma_O(S)$ acts vertex-transitively on the product of trees $\prod_{v\in S\backslash S_0} \mathcal{T}_v$ if and only if $S\backslash S_0 \subseteq V^K\backslash \Ram(B)$. 
\item The fourth column gives a group isomorphic to $\Gamma_O(\Ram(B))$.
By \Cref{lemma:VertexStabGammaS}, this is the stabiliser of the vertex $O$ for the action of $\Gamma_O(S)$ on the product of trees, if $S_0 (= S\cap \Ram(B)) = \Ram(B)$.
\item The fifth column gives a group isomorphic to $\Gamma_O(S_0)$, for any $S_0\cap V^K_f$ in the third column.
By \Cref{lemma:VertexStabGammaS}, this is the stabiliser of the vertex $O$ for the action of $\Gamma_O(S)$
\item The last column is the index $[\Idl(O) : \PIdl(O)]$.
\end{itemize}

\renewcommand{\arraystretch}{1.73}

\begin{longtable}{|c|c||c|c|c|c|c|}
\hline
\multicolumn{2}{|c||}{$X-1$} & \multicolumn{5}{c|}{} \\
\hline
1 & $\left(\frac{-1,-1}{\QQ}\right)$ & $\{2\}$ & $\emptyset$ & $\Sym(4)$ & $\Alt(4)$ & $1$  \\ 
\hline
2 & $\left(\frac{-3,-1}{\QQ}\right)$ & $\{3\}$ & $\emptyset$ & $D_6$ & $\Sym(3)$ & $1$ \\ 
\hline 
3 & $\left(\frac{-2,-5}{\QQ}\right)$ & $\{5\}$ & $\emptyset$ & $\Sym(3)$ & $C_3$ & $1$ \\ 
\hline 
4 & $\left(\frac{-3,-10}{\QQ}\right)$ & $\{2,3,5\}$ & $\{2\}, \{5\}$ & $D_6$ & $C_3$ & $2$ \\
\hline 
5 & $\left(\frac{-7,-1}{\QQ}\right)$ & $\{7\}$ & $\emptyset$ & $C_2^2$ & $C_2$ & $1$ \\
\hline 
6 & $\left(\frac{-21,-1}{\QQ}\right)$ & $\{2,3,7\}$ & $\{3\},\{7\}$ & $D_4$ & $C_2$ & $2$ \\
\hline 
7 & $\left(\frac{-2,-35}{\QQ}\right)$ & $\{2,5,7\}$ & $\{5\},\{7\}$ & $C_2^2$ & $1$ & $2$ \\
\hline 
8 & $\left(\frac{-2,-13}{\QQ}\right)$ & $\{13\}$ & $\emptyset$ & $C_2$ & $1$ & $1$ \\
\hline 
9 & $\left(\frac{-13,-6}{\QQ}\right)$ & $\{2,3,13\}$ & $\{2\},\{3\}$ & $C_2^2$ & $1$ & $2$ \\
\hline
\hline
\multicolumn{2}{|c||}{$X^2-X-1$} & \multicolumn{5}{c|}{}\\
\hline
10 & $\left(\frac{-1,-1}{K}\right)$ & $\emptyset$ & $\emptyset$ & $\Alt(5)$ & $\Alt(5)$ & $1$ \\
\hline 
11 & $\left(\frac{-4-2\alpha,-2}{K}\right)$ & $\{2, 1- 2 \alpha\}$ & $\emptyset$ & $D_{10}$ & $C_5$ & $1$ \\
\hline 
12 & $\left(\frac{-6,-3+\alpha}{K}\right)$ & $\{2, 3\}$ & $\{2\}, \{3\}$ & $D_5$ & $C_5$ & $2$ \\
\hline 
13 & $\left(\frac{-186+8\alpha,-25-7\alpha}{K}\right)$ & $\{2, 1 -3 \alpha\}$ & $\emptyset$ & $D_4$ & $C_2$ & $1$ \\
\hline 
14 & $\left(\frac{-133-19\alpha,-6+2\alpha}{K}\right)$ & $\{-1+2\alpha, -1+3\alpha\}$ & $\{-1+2\alpha\}, \{-1+3\alpha\}$ & $\Sym(3)$ & $C_3$ & $2$ \\
\hline 
15 & $\left(\frac{-178-8\alpha,-32+7\alpha}{K}\right)$ & $\{2, 2-3\alpha\}$ & $\emptyset$ & $D_4$ & $C_2$ & $1$ \\
\hline 
16 & $\left(\frac{-112+14\alpha,-6-2\alpha}{K}\right)$ & $\{-1+2\alpha, 2-3\alpha\}$ & $\{-1+2\alpha\}, \{2-3\alpha\}$ & $\Sym(3)$ & $C_3$ & $2$ \\
\hline 
17 & $\left(\frac{-432-36\alpha,-2}{K}\right)$ & $\{-1+2\alpha, 2-5\alpha\}$ & $\{-1+2\alpha\}, \{2-5\alpha\}$ & $C_2$ & $1$ & $2$ \\
\hline 
18 & $\left(\frac{-338+26\alpha,-2}{K}\right)$ & $\{-1+2\alpha, -3+5\alpha\}$ & $\{-1+2\alpha\}, \{-3+5\alpha\}$ & $C_2$ & $1$ & $2$ \\
\hline 
19 & $\left(\frac{-2410-60\alpha,-41+5\alpha}{K}\right)$ & $\{2, 6+\alpha\}$ & $\{2\}, \{6+\alpha\}$ & $C_2$ & $1$ & $2$ \\
\hline 
20 & $\left(\frac{-2470+60\alpha,-41-46\alpha}{K}\right)$ & $\{2, 7-\alpha\}$ & $\{2\}, \{7-\alpha\}$ & $C_2$ & $1$ & $2$ \\
\hline 
\hline
\multicolumn{2}{|c||}{$X^2-2$} & \multicolumn{5}{c|}{} \\
\hline
21 & $\left(\frac{-1,-1}{K}\right)$ & $\emptyset$ & $\emptyset$ & $\Sym(4)$ & $\Sym(4)$ & $1$ \\
\hline 
22 & $\left(\frac{-46+\alpha,-13-3\alpha}{K}\right)$ & $\{\alpha, -1+2\alpha\}$ & $\emptyset$ & $D_8$ & $C_4$ & $1$ \\
\hline 
23 & $\left(\frac{-46-\alpha,-13+3\alpha}{K}\right)$ & $\{\alpha, 1+2\alpha\}$ & $\emptyset$ & $D_8$ & $C_4$ & $1$ \\
\hline 
24 & $\left(\frac{-24+15\alpha,-3}{K}\right)$ & $\{\alpha, 3\}$ & $\emptyset$ & $D_6$ & $C_3$ & $1$ \\
\hline 
25 & $\left(\frac{-213+27\alpha,-25-10\alpha}{K}\right)$ & $\{-1+2\alpha, 3\}$ & $\{-1+2\alpha\}, \{3\}$ & $C_2$ & $1$ & $2$ \\
\hline 
26 & $\left(\frac{-213-27\alpha,-25+10\alpha}{K}\right)$ & $\{-1-2\alpha, 3\}$ & $\{-1-2\alpha\}, \{3\}$ & $C_2$ & $1$ & $2$ \\
\hline
27 & $\left(\frac{-246-7\alpha,-3}{K}\right)$ & $\{\alpha, 1+3\alpha\}$ & $\{\alpha\}, \{1+3\alpha\}$ & $\Sym(3)$ & $C_3$ & $2$ \\
\hline
28 & $\left(\frac{-246+7\alpha,-3}{K}\right)$ & $\{\alpha, 1-3\alpha\}$ & $\{\alpha\}, \{1-3\alpha\}$ & $\Sym(3)$ & $C_3$ & $2$ \\
\hline
29 & $\left(\frac{-70+45\alpha,-2+\alpha}{K}\right)$ & $\{\alpha, 5 \}$ & $\emptyset$ & $C_2^2$ & $1$ & $1$ \\
\hline
\hline
\multicolumn{2}{|c||}{$X^2-3$} & \multicolumn{5}{c|}{} \\
\hline
30 & $\left(\frac{-3-\alpha,-1}{K}\right)$ & $\{1+\alpha, \alpha\}$ & $\{1+\alpha\}, \{\alpha\}$ & $D_{12}$ & $C_{12}$ & $2$ \\
\hline
31 & $\left(\frac{-123-5\alpha,-2-\alpha}{K}\right)$ & $\{1+\alpha, 4-\alpha\}$ & $\{1+\alpha\}$ & $C_2^2$ & $C_2$ & $2$ \\
\hline
32 & $\left(\frac{-123-5\alpha,-21-11\alpha}{K}\right)$ & $\{\alpha, 4-\alpha\}$ & $\{\alpha\}$ & $C_2$ & $1$ & $2$ \\
\hline
33 & $\left(\frac{-123+5\alpha,-2+\alpha}{K}\right)$ & $\{1+\alpha, 4+\alpha\}$ & $\{1+\alpha\}$ & $C_2^2$ & $C_{2}$ & $2$ \\
\hline
34 & $\left(\frac{-30-14\alpha,-4-\alpha}{K}\right)$ & $\{\alpha, 4+\alpha\}$ & $\{\alpha\}$ & $C_2$ & $1$ & $2$ \\
\hline
35 & $\left(\frac{-75-5\alpha,-45-25\alpha}{K}\right)$ & $\{1-\alpha, 5\}$ & $\{1-\alpha\}$ & $C_2$ & $1$ & $2$ \\
\hline
\hline
\multicolumn{2}{|c||}{$X^2-X-3$} & \multicolumn{5}{c|}{} \\
\hline
36 & $\left(\frac{-1,-1}{K}\right)$ & $\emptyset$ & $\emptyset$ & $\Alt(4)$ & $\Alt(4)$ & $1$ \\
\hline
37 & $\left(\frac{-10+4\alpha,-1}{K}\right)$ & $\{1-\alpha,2\}$ & $\emptyset$ & $D_4$ & $C_2$ & $1$ \\
\hline
38 & $\left(\frac{-6-4\alpha,-1}{K}\right)$ & $\{\alpha,2\}$ & $\emptyset$ & $D_4$ & $C_2$ & $1$ \\
\hline
39 & $\left(\frac{-3,-1}{K}\right)$ & $\{\alpha,1-\alpha\}$ & $\{\alpha\}, \{1-\alpha\}$ & $D_6$ & $\Sym(3)$ & $2$ \\
\hline
40 & $\left(\frac{-19-14\alpha,-10-6\alpha}{K}\right)$ & $\{1-\alpha,1-2\alpha\}$ & $\{1-\alpha\},\{1-2\alpha\}$ & $C_2$ & $1$ & $2$ \\
\hline
41 & $\left(\frac{-6-\alpha,-10-6\alpha}{K}\right)$ & $\{\alpha,1-2\alpha\}$ & $\{\alpha\},\{1-2\alpha\}$ & $C_2$ & $1$ & $2$ \\
\hline
\hline
\multicolumn{2}{|c||}{$X^2-X-4$} & \multicolumn{5}{c|}{} \\
\hline
42 & $\left(\frac{-1,-18+7\alpha}{K}\right)$ & $\emptyset$ & $\emptyset$ & $\Sym(3)$ & $\Sym(3)$ & $1$ \\
\hline
43 & $\left(\frac{-2,-1}{K}\right)$ & $\{2-\alpha,1+\alpha\}$ & $\{2-\alpha\},\{1+\alpha\}$ & $\Sym(4)$ & $\Alt(4)$ & $2$ \\
\hline
44 & $\left(\frac{-282-7\alpha,-1317+514\alpha}{K}\right)$ & $\{2-\alpha,3-2\alpha\}$ & $\{2-\alpha\},\{3-2\alpha\}$ & $C_2$ & $1$ & $2$ \\
\hline
45 & $\left(\frac{-289+7\alpha,-979+382\alpha}{K}\right)$ & $\{1+\alpha,1+2\alpha\}$ & $\{1+\alpha\},\{1+2\alpha\}$ & $C_2$ & $1$ & $2$ \\
\hline
\hline
\multicolumn{2}{|c||}{$X^2-X-5$} & \multicolumn{5}{c|}{} \\
\hline
46 & $\left(\frac{-20-2\alpha,-11-5\alpha}{K}\right)$ & $\{2-\alpha,2\}$ & $\{2-\alpha\}$ & $C_4$ & $C_2$ & $2$ \\
\hline
47 & $\left(\frac{-10+2\alpha,-20-11\alpha}{K}\right)$ & $\{2,\alpha\}$ & $\{\alpha\}$ & $C_2$ & $1$ & $2$ \\
\hline
48 & $\left(\frac{-26-14\alpha,-13-7\alpha}{K}\right)$ & $\{2,1-\alpha\}$ & $\{1-\alpha\}$ & $C_2$ & $1$ & $2$ \\
\hline
49 & $\left(\frac{-10-\alpha,-1}{K}\right)$ & $\{1+\alpha,3+\alpha\}$ & $\{1+\alpha,3+\alpha\}$ & $C_2$ & $C_2$ & $4$ \\
\hline
\hline
\multicolumn{2}{|c||}{$X^2-6$} & \multicolumn{5}{c|}{} \\
\hline
50 & $\left(\frac{-18+7\alpha,-1}{K}\right)$ & $\{2+\alpha,3-\alpha\}$ & $\{2+\alpha\}$ & $D_4$ & $C_4$ & $2$ \\
\hline
51 & $\left(\frac{-18-7\alpha,-25+10\alpha}{K}\right)$ & $\{3+\alpha,1-\alpha\}$ & $\{-1+\alpha\}$ & $C_2$ & $1$ & $2$ \\
\hline
52 & $\left(\frac{-27-2\alpha,-25+10\alpha}{K}\right)$ & $\{3-\alpha,1+\alpha\}$ & $\{1+\alpha\}$ & $C_2$ & $1$ & $2$ \\
\hline
\hline
\multicolumn{2}{|c||}{$X^2-7$} & \multicolumn{5}{c|}{} \\
\hline
53 & $\left(\frac{-9+3\alpha,-3+\alpha}{K}\right)$ & $\{3-\alpha,2-\alpha\}$ & $\{2-\alpha\}$ & $\Sym(3)$ & $C_3$ & $2$ \\
\hline
54 & $\left(\frac{-9+3\alpha,-45+17\alpha}{K}\right)$ & $\{3+\alpha,2+\alpha\}$ & $\{2+\alpha\}$ & $\Sym(3)$ & $C_3$ & $2$ \\
\hline
55 & $\left(\frac{-7-\alpha,-8+3\alpha}{K}\right)$ & $\{3-\alpha,\alpha\}$ & $\{3-\alpha,\alpha\}$ & $C_2$ & $C_2$ & $4$ \\
\hline
\hline
\multicolumn{2}{|c||}{$X^2-X-8$} & \multicolumn{5}{c|}{} \\
\hline
56 & $\left(\frac{-10-4\alpha,-1}{K}\right)$ & $\{2+\alpha,-7+2\alpha\}$ & $\{2+\alpha\}$ & $C_2^2$ & $C_2$ & $2$ \\
\hline
57 & $\left(\frac{-10-4\alpha,-38-16\alpha}{K}\right)$ & $\{3-\alpha,5+2\alpha\}$ & $\{3-\alpha\}$ & $C_2^2$ & $C_2$ & $2$ \\
\hline
\hline
\multicolumn{2}{|c||}{$X^2-15$} & \multicolumn{5}{c|}{} \\
\hline
58 & $\left(\frac{-75-19\alpha,-4-\alpha}{K}\right)$ & $\{(2,1+\alpha), (3,\alpha)\}$ & $\{(2,1+\alpha), (3,\alpha)\}$ & $C_2$ & $C_2$ & $8$ \\
\hline
\hline
\multicolumn{2}{|c||}{$X^3-X^2-2X+1$} & \multicolumn{5}{c|}{} \\
\hline
59 & $\left(\frac{-2,-1}{K}\right)$ & $\{2\}$ & $\emptyset$ & $\Sym(4)$ & $\Alt(4)$ & $1$ \\
\hline
60 & $\left(\frac{-2-\alpha,-1}{K}\right)$ & $\{2-\alpha-\alpha^2\}$ & $\emptyset$ & $D_{14}$ & $D_7$ & $1$ \\
\hline
61 & $\left(\frac{-13,-2\alpha-4\alpha^2}{K}\right)$ & $\{-2-2\alpha+\alpha^2\}$ & $\emptyset$ & $D_7$ & $C_7$ & $1$ \\
\hline
62 & $\left(\frac{-13,6-10\alpha-18\alpha^2}{K}\right)$ & $\{3-\alpha-\alpha^2\}$ & $\emptyset$ & $D_7$ & $C_7$ & $1$ \\
\hline
63 & $\left(\frac{-13,-4+2\alpha-6\alpha^2}{K}\right)$ & $\{2+\alpha-2\alpha^2\}$ & $\emptyset$ & $D_7$ & $C_7$ & $1$ \\
\hline
64 & $\left(\frac{-29,-14+32\alpha-20\alpha^2}{K}\right)$ & $\{-3+2\alpha+\alpha^2\}$ & $\emptyset$ & $\Sym(3)$ & $C_3$ & $1$ \\
\hline
65 & $\left(\frac{-29,-12+34\alpha-26\alpha^2}{K}\right)$ & $\{4+\alpha-3\alpha^2\}$ & $\emptyset$ & $\Sym(3)$ & $C_3$ & $1$ \\
\hline
66 & $\left(\frac{-29,-4+6\alpha-4\alpha^2}{K}\right)$ & $\{3+3\alpha-2\alpha^2\}$ & $\emptyset$ & $\Sym(3)$ & $C_3$ & $1$ \\
\hline
67 & $\left(\frac{-43,-1-\alpha^2}{K}\right)$ & $\{5-\alpha-2\alpha^2\}$ & $\emptyset$ & $C_2^2$ & $C_2$ & $1$ \\
\hline
68 & $\left(\frac{-43,-1+2\alpha-2\alpha^2}{K}\right)$ & $\{2+3\alpha-\alpha^2\}$ & $\emptyset$ & $C_2^2$ & $C_2$ & $1$ \\
\hline
69 & $\left(\frac{-43,-\alpha-2\alpha^2}{K}\right)$ & $\{3+2\alpha-3\alpha^2\}$ & $\emptyset$ & $C_2^2$ & $C_2$ & $1$ \\
\hline
\hline
\multicolumn{2}{|c||}{$X^3-3X+1$} & \multicolumn{5}{c|}{} \\
\hline
70 & $\left(\frac{-2-\alpha,-1}{K}\right)$ & $\{1-\alpha^2\}$ & $\emptyset$ & $D_{18}$ & $D_9$ & $1$ \\
\hline
71 & $\left(\frac{-19,-5-7\alpha-4\alpha^2}{K}\right)$ & $\{5+2\alpha-2\alpha^2\}$ & $\emptyset$ & $C_2^2$ & $C_2$ & $1$ \\
\hline
72 & $\left(\frac{-19,-9-26\alpha-14\alpha^2}{K}\right)$ & $\{1-2\alpha\}$ & $\emptyset$ & $C_2^2$ & $C_2$ & $1$ \\
\hline
73 & $\left(\frac{-19,-10-15\alpha-6\alpha^2}{K}\right)$ & $\{3-2\alpha^2\}$ & $\emptyset$ & $C_2^2$ & $C_2$ & $1$ \\
\hline
74 & $\left(\frac{-37,-4-6\alpha-4\alpha^2}{K}\right)$ & $\{4+3\alpha-2\alpha^2\}$ & $\emptyset$ & $C_2$ & $1$ & $1$ \\
\hline
75 & $\left(\frac{-37,-52-128\alpha-62\alpha^2}{K}\right)$ & $\{6+\alpha-3\alpha^2\}$ & $\emptyset$ & $C_2$ & $1$ & $1$ \\
\hline
76 & $\left(\frac{-37,-56-170\alpha-88\alpha^2}{K}\right)$ & $\{2-2\alpha-\alpha^2\}$ & $\emptyset$ & $C_2$ & $1$ & $1$ \\
\hline
\hline
\multicolumn{2}{|c||}{$X^3-X^2-3X+1$} & \multicolumn{5}{c|}{} \\
\hline
77 & $\left(\frac{-\alpha-\alpha^2,-1}{K}\right)$ & $\{1-\alpha\}$ & $\emptyset$ & $\Sym(4)$ & $\Alt(4)$ & $1$ \\
\hline
78 & $\left(\frac{-5,-\alpha-\alpha^2}{K}\right)$ & $\{-1-\alpha+\alpha^2\}$ & $\emptyset$ & $\Sym(3)$ & $C_3$ & $1$ \\
\hline
79 & $\left(\frac{-13,-\alpha-\alpha^2}{K}\right)$ & $\{2+2\alpha-\alpha^2\}$ & $\emptyset$ & $C_2$ & $1$ & $1$ \\
\hline
\hline
\multicolumn{2}{|c||}{$X^3-X^2-4X-1$} & \multicolumn{5}{c|}{} \\
\hline
80 & $\left(\frac{-5,-2-10\alpha-12\alpha^2}{K}\right)$ & $\{3+2\alpha-\alpha^2\}$ & $\emptyset$ & $\Sym(3)$ & $C_3$ & $1$ \\
\hline
81 & $\left(\frac{-5,-2-4\alpha-4\alpha^2}{K}\right)$ & $\{1-\alpha\}$ & $\emptyset$ & $\Sym(3)$ & $C_3$ & $1$ \\
\hline
82 & $\left(\frac{-5,-2-8\alpha-6\alpha^2}{K}\right)$ & $\{2+\alpha-\alpha^2\}$ & $\emptyset$ & $\Sym(3)$ & $C_3$ & $1$ \\
\hline
83 & $\left(\frac{-4-\alpha,-3-3\alpha-3\alpha^2}{K}\right)$ & $\{3-\alpha^2\}$ & $\emptyset$ & $C_2$ & $1$ & $1$ \\
\hline
\hline
\multicolumn{2}{|c||}{$X^3-4X-1$} & \multicolumn{5}{c|}{} \\
\hline
84 & $\left(\frac{-2,-1}{K}\right)$ & $\{1+\alpha\}$ & $\{1+\alpha\}$ & $\Alt(4)$ & $\Alt(4)$ & $2$ \\
\hline
85 & $\left(\frac{-2,-2-\alpha}{K}\right)$ & $\{3+\alpha-\alpha^2\}$ & $\{3+\alpha-\alpha^2\}$ & $C_2\times C_2$ & $C_2\times C_2$ & $2$ \\
\hline
86 & $\left(\frac{-7,-1}{K}\right)$ & $\{2-\alpha^2\}$ & $\{2-\alpha^2\}$ & $C_2$ & $C_2$ & $2$ \\
\hline
\hline
\multicolumn{2}{|c||}{$X^3-X^2-4X+3$} & \multicolumn{5}{c|}{} \\
\hline
87 & $\left(\frac{-3,-1}{K}\right)$ & $\{\alpha\}$ & $\{\alpha\}$ & $\Sym(3)$ & $\Sym(3)$ & $2$ \\
\hline
88 & $\left(\frac{-5,-2}{K}\right)$ & $\{1+\alpha\}$ & $\{1+\alpha\}$ & $C_3$ & $C_3$ & $2$ \\
\hline
89 & $\left(\frac{-7,-1}{K}\right)$ & $\{2-\alpha^2\}$ & $\{2-\alpha^2\}$ & $C_2$ & $C_2$ & $2$ \\
\hline
\hline
\multicolumn{2}{|c||}{$X^3-X^2-4X+2$} & \multicolumn{5}{c|}{} \\
\hline
90 & $\left(\frac{-2+5\alpha-2\alpha^2,-9-\alpha+2\alpha^2}{K}\right)$ & $\{\alpha\}$ & $\emptyset$ & $\Sym(3)$ & $C_3$ & $1$ \\
\hline
\hline
\multicolumn{2}{|c||}{$X^3-X^2-4X+1$} & \multicolumn{5}{c|}{} \\
\hline
91 & $\left(\frac{-2-\alpha,-1}{K}\right)$ & $\{1-\alpha\}$ & $\emptyset$ & $C_2^2$ & $C_2$ & $1$ \\
\hline
\hline
\multicolumn{2}{|c||}{$X^4-X^3-3X^2+X+1$} & \multicolumn{5}{c|}{} \\
\hline
92 & $\left(\frac{-1,-1}{K}\right)$ & $\emptyset$ & $\emptyset$ & $\Alt(5)$ & $\Alt(5)$ & $1$ \\
\hline
\hline
\multicolumn{2}{|c||}{$X^4-X^3-4X^2+4X+1$} & \multicolumn{5}{c|}{} \\
\hline
93 & $\left(\frac{-2+20\alpha-22\alpha^2+6\alpha^3,-2}{K}\right)$ & $\{1+\alpha,2\}$ & $\{1+\alpha\}$ & $C_2$ & $1$ & $2$ \\
\hline
\hline
\multicolumn{2}{|c||}{$X^4-4X^2-X+1$} & \multicolumn{5}{c|}{} \\
\hline
94 & $\left(\frac{-1,-1}{K}\right)$ & $\emptyset$ & $\emptyset$ & $\Alt(4)$ & $\Alt(4)$ & $1$ \\
\hline
95 & $\left(\frac{-106-9\alpha-10\alpha^2+8\alpha^3,-3-4\alpha-2\alpha^2}{K}\right)$ & $
\{1+3\alpha-\alpha^3,2-\alpha^2\}$ & $\{1+3\alpha-\alpha^3\},\{2-\alpha^2\}$ & $C_2^2$ & $C_2$ & $2$ \\
\hline
\hline
\multicolumn{2}{|c||}{$X^4-5X^2+5$} & \multicolumn{5}{c|}{} \\
\hline
96 & $\left(\frac{-5-\alpha-\alpha^2,2-\alpha-2\alpha^2+\alpha^3}{K}\right)$ & $
\{-4-3\alpha+\alpha^2+\alpha^3,-3\alpha+\alpha^3\}$ & $\{-4-3\alpha+\alpha^2+\alpha^3\}$ & $C_2^2$ & $C_2$ & $2$ \\
\hline
\hline
\multicolumn{2}{|c||}{$X^4-4X^2+1$} & \multicolumn{5}{c|}{}  \\
\hline
97 & $\left(\frac{-4-7\alpha-\alpha^2+2\alpha^3,-2\alpha^2+\alpha^3}{K}\right)$ & $
\{1-4\alpha+\alpha^3,2-\alpha^2\}$ & $\{1-4\alpha+\alpha^3\}$ & $C_2^2$ & $C_2$ & $2$ \\
\hline
\hline
\multicolumn{2}{|c||}{$X^4-X^3-4X^2+X+2$} & \multicolumn{5}{c|}{}  \\
\hline
98 & $\left(\frac{-1,1-3\alpha^2+\alpha^3}{K}\right)$ & $
\emptyset$ & $\emptyset$ & $\Sym(3)$ & $\Sym(3)$ & $1$ \\
\hline
\hline
\multicolumn{2}{|c||}{$X^4-6X^2-4X+2$} & \multicolumn{5}{c|}{}  \\
\hline
99 & $\left(\frac{-15\alpha-12\alpha^2+\alpha^3,-2-\alpha}{K}\right)$ & $\{-5\alpha-\alpha^2+\alpha^3,1+5\alpha+\alpha^2-\alpha^3\}$ & $\{1+5\alpha+\alpha^2-\alpha^3\}$ & $C_2$ & $1$ & $2$ \\
\hline
100 & $\left(\frac{6-25\alpha-43\alpha^2-15\alpha^3,-\alpha-\alpha^2}{K}\right)$ & $\{\alpha,1-\alpha\}$ & $\{1-\alpha\}$ & $C_2$ & $1$ & $2$ \\
\hline
\hline
\multicolumn{2}{|c||}{$X^4-2X^3-3X^2+4X+1$} & \multicolumn{5}{c|}{}  \\
\hline
101 & $\left(\frac{-1+4\alpha-2\alpha^2-\alpha^3,-1+4\alpha-4\alpha^2+\alpha^3}{K}\right)$ & $\{1+\alpha,3-2\alpha-2\alpha^2+\alpha^3\}$ & $\{1+\alpha\}$ & $C_2$ & $1$ & $2$ \\
\hline
\hline
\multicolumn{2}{|c||}{$X^5-5X^3-X^2+3X+1$} & \multicolumn{5}{c|}{}  \\
\hline
102 & $\left(\frac{-5,-2}{K}\right)$ & $\{3+2\alpha-9\alpha^2-\alpha^3+2\alpha^4\}$ & $\emptyset$ & $\Sym(3)$ & $C_3$ & $1$ \\
\hline
\hline
\multicolumn{2}{|c||}{$X^5-2X^4-3X^3+5X^2+X-1$} & \multicolumn{5}{c|}{}  \\
\hline
103 & $\left(\frac{-3,-1}{K}\right)$ & $\{-1-4\alpha+3\alpha^2+2\alpha^3-\alpha^4\}$ & $\{-1-4\alpha+3\alpha^2+2\alpha^3-\alpha^4\}$ & $\Sym(3)$ & $\Sym(3)$ & $2$ \\
\hline
\hline
\multicolumn{2}{|c||}{$X^5-5X^3+4X-1$} & \multicolumn{5}{c|}{}  \\
\hline
104 & $\left(\frac{-7,-1}{K}\right)$ & $\{2-\alpha^2\}$ & $\{2-\alpha^2\}$ & $C_2$ & $C_2$ & $2$ \\
\hline
105 & $\left(\frac{-13,-2}{K}\right)$ & $\{-4+4\alpha+9\alpha^2-\alpha^3-2\alpha^4\}$ & $\{-4+4\alpha+9\alpha^2-\alpha^3-2\alpha^4\}$ & $1$ & $1$ & $2$ \\
\hline
\end{longtable}

\end{landscape}

\end{document}